\documentclass{article}
\RequirePackage{amsthm,amsmath,amssymb,amsthm,hyperref,xcolor}
\def\MR#1{\href{http://www.ams.org/mathscinet-getitem?mr=#1}{MR#1}}
\usepackage[utf8]{inputenc}

\numberwithin{equation}{section}
\theoremstyle{plain}
\newtheorem{thm}{Theorem}[section]

\newtheorem{pro}{Proposition}
\newtheorem*{hyp}{Hypothesis}
\newtheorem{lem}{Lemma}

\begin{document}

\title{Auxiliary information : the raking-ratio empirical process.}
\author{Mickael Albertus\footnote{mickael.albertus@math.univ-toulouse.fr}\  and Philippe Berthet\footnote{philippe.berthet@math.univ-toulouse.fr}}

\maketitle

\begin{center}
Institut de Mathématiques de Toulouse ; UMR5219 \\
Université de Toulouse ; CNRS
\end{center}

\begin{abstract}
We study the empirical measure associated to a sample of size $n$ and modified
by $N$ iterations of the raking-ratio method. This empirical measure is adjusted
to match the true probability of sets in a finite partition which changes
each step. We establish asymptotic properties of the raking-ratio empirical
process indexed by functions as $n\rightarrow +\infty$, for $N$ fixed. We study 
nonasymptotic properties by using a Gaussian approximation which yields uniform 
Berry-Esseen type bounds depending on $n, N$ and provides estimates of the uniform 
quadratic risk reduction. A closed-form expression of the limiting covariance matrices 
is derived as $N\rightarrow +\infty$. In the two-way contingency table case the limiting
process has a simple explicit formula.
\end{abstract}

\tableofcontents

\section{Introduction}

\subsection{The raking-ratio method}

\noindent In survey analysis, statistics, economics and computer sciences the 
raking-ratio iterative procedure aims to exploit the knowledge of one or several
marginals of a discrete multivariate distribution to fit the data after sampling. 
Despite many papers from the methodological and algorithmic viewpoint, and chapters 
in classical textbooks for statisticians, economists or engineers, no probabilistic
study is available to take into account that the entries of the algorithm are random
and the initial discrete measure is empirical. We intend to fill this gap. Let us
first describe the algorithm, usually considered with deterministic
entries, then recall the few known results and state the open question to be addressed.\smallskip

\noindent\textbf{The raking-ratio algorithm.} A sample is drawn from a population $P$ for which
$k\geqslant 2$ marginal finite discrete distributions are explicitly known. Initially, each
data point has a weight $1/n$. The \textit{ratio} step of the algorithm
consists in computing new weights in such a way that the modified empirical joint distribution
has the currently desired marginal. The \textit{raking} step
consists in iterating the correction according to another known marginal law,
changing again all the weights. The $k$ margin constraints are usually treated 
in a periodic order, only one being fulfilled at the same time. 
The raking-ratio method stops after $N$ iterations with the
implicit hope that the previous constraints are still almost satisfied.
See Section~\ref{RRExemple} for an elementary numerical example with $k=2$ and 
Section~\ref{InfNPart} for notation and mathematical definition of the algorithm.\smallskip

\noindent\textbf{The limit.} This algorithm was called
\textit{iterative proportions} by Deming and Stephan \cite{DemSte} who first
introduced it. They showed that the $k$ margins converge to the desired ones as 
$N\rightarrow+\infty$. They even claimed that if the frequencies of a multiway
contingency table are raked periodically as $N\rightarrow+\infty$ 
they converge to the frequencies minimizing the chi-square distance to
the initial frequencies, under the $k$ margin constraints. Two years later, Stephan \cite{Ste} 
observed that it is wrong and modified the algorithm accordingly to achieve the chi-square distance minimization.
Lewis \cite{Lewis} and Brown \cite{Brown} studied the case of Bernoulli marginals
from the Shannon entropy minimization viewpoint. When $k=2$ a two-way contingency table can be
viewed as a matrix. Sinkhorn \cite{Sinkh64,Sinkh67} proved that a unique doubly stochastic 
matrix can be obtained from each positive square matrix by alternately normalizing its rows 
and its columns, which shows that the algorithm converges in this special case. 
Finally Ireland and Kullback \cite{IrrKull} 
generalized previous arguments to rigorously justify that the raking-ratio
converges to the unique projection of the empirical measure in
Shannon-Kullback-Leibler relative entropy on the set of discrete
distributions satisfying the $k$ margin constraints. From a numerical viewpoint, the rate of convergence of the algorithm is geometric, see Franklin and Lorentz \cite{FranLor89}.\smallskip

\noindent\textit{Remark A.} When minimizing contrasts such as discrimination information, chi-square distance 
or likelihood, the minimizers are not explicit due to the nonlinearity of sums of ratios showing up in derivatives. 
This is why converging algorithms are used in practice. In the case of the iterative proportions algorithm 
each step is easily and fastly computed. What has been studied concerns the convergence 
when the iterations $N\rightarrow+\infty$, with $n$ fixed and initial empirical frequencies treated as 
deterministic entries. When the sample size $n\rightarrow+\infty$, these entries are close to $P$ itself, 
which satisfies the marginal constraints, hence one expects that the number $N$ of iterations necessary to 
converge is small. We shall study the $N_0$ first iterations in the statistical setting $n\rightarrow+\infty$.\smallskip

\noindent\textbf{Non explicit bias and variance.} The initial values being empirical frequencies 
the converged solution of the algorithm as $N\rightarrow+\infty$ is a joint distributions fulfilling the marginal 
requirements that still deviates from the true population distribution $P$, and moreover in a rather complicated way. 
The modified empirical distribution satisfying only the marginal constraint of the current iteration, 
there is a permanent bias with respect to other margins, and hence with $P$. The exact covariance matrix 
and bias vector of the random weights after $N$ iterations are tedious to compute. For instance, estimates
for the variance of cell probabilities in the case of a two-way contingency
table are given by Brackstone and Rao \cite{BraRao} for $N\leqslant4$, Konijn \cite{Kon}
or Choudhry and Lee \cite{ChouLee} for $N=2$. Bankier \cite{Bank} proposed a recursive linearization 
technique providing an estimator of the asymptotic variance of weights. In Binder and Th\'{e}berge
\cite{BinThe} the variance of the converged solution requires to calculate weights at each 
iteration.\smallskip

\noindent\textbf{Open question.} Since exact computations lead to intractable formulas for the bias 
and variance of frequencies and statistics as simple as means, an important open problem is to identify
leading terms when $n$ is large compared to $N$. We derive comprehensive explicit formulas as 
$n\rightarrow+\infty$ for $N\leqslant N_0$ and $N_0$ fixed, then for $N\rightarrow+\infty$. In order to further analyze the raking ratio 
method it is moreover desirable to control simultaneously large classes of statistics and hence 
to work at the empirical measure level rather than with the empirical weights or a single statistic only. This is the main motivation for 
the forthcoming general study of empirical measures indexed by functions and modified through auxiliary information 
given by partitions.\smallskip

\subsection{Statistical motivation}

\noindent\textbf{Representative sample.} In a survey analysis context, the
raking-ratio method modifies weights of a contingency table built from a
sample of size $n$ in order to fit exactly given marginals. Such a strict
margin correction is justified when a few properties of the finite population
under study are known, like the size of sub-populations. The modified sample
frequencies then reflect the marginal structure of the whole population.
If the population is large or infinite the information may come from
previous and independent statistical inference, from structural properties
of the model or from various experts.\smallskip

\noindent\textit{Remark B.} Making the sample representative of the population
is an ad hoc approach based on common sense. The mathematical impact is
twofold. On the one hand all statistics are affected by the new weights in terms
of bias, variance and limit law so that practitioners may very well be using
estimators, tests or confident bands that have lost their usual properties.
On the other hand, replacing marginal frequencies with the true ones may smooth 
sample fluctuations of statistics correlated to them while leaving
the uncorrelated ones rather unaffected. These statements will be quantified precisely
at Section~\ref{GenPro}.\smallskip

\noindent\textit{Remark C.} Fitting after sampling is a natural method that has 
been re-invented many times in various fields, and was probably used long time ago. 
Depending on the setting it may be viewed as stratification, calibration, iterating 
proportional fitting, matrix scaling and could be used to deal with missing data. 
Many fitting methods may be reduced to a raking-ratio type algorithm.
We initially called it auxiliary information of partitions as we re-invented
it as a special case of the nonparametric partial information problem stated 
in Section~\ref{AuxInfView}.\smallskip

\noindent\textit{Remark D.} An asymptotic approach is no more
relevant in survey analysis when the underlying population is rather small. In
the small population case, the way the sample is drawn has a so deep impact that it may even
become the main topic. A study of calibration methods for finite population
can be found in Deville and S\"{a}rndal \cite{DevSar92,DevSar93}. This is beyond
the scope of our work.\smallskip

\noindent\textbf{Quadratic risk reduction.} Modifying marginals frequencies of a sample may 
induce serious drawbacks. One should ask whether or not the estimation risk can be controlled. 
Typically, a statistic has more bias when sample weights are changed by
using raking, calibration or stratification methods after sampling. In the
spirit of Remark B, a variance reduction is expected if the 
statistic of interest is strongly correlated to the $k$ known
discrete marginal variables. Now, evaluating the quadratic risk of a specific statistic
requires tedious expansions for the bias, variance and correlations of weights, whence the 
very small $N$ studied in the literature. Likewise, no global risk reduction property has 
been established as $n\rightarrow+\infty$ and no multivariate or uniform central limit 
theorem. These results are established at Propositions~\ref{ProCVProc} to~\ref{Thm3}.\smallskip

\noindent\textbf{Contributions.} In this paper we consider classes of empirical means
raked $N$ times, sampled jointly from any population. We derive closed-form
expressions of their Gaussian limits and their limiting covariances as
$n\rightarrow+\infty$ then $N\rightarrow+\infty$. We also quantify
the uniform risk reduction phenomenon and provide sharp statistical estimation
tools such as uniform Berry-Esseen type bounds. In particular, a Donsker invariance 
principle for the raked empirical process provides joint limiting laws for additive 
statistics built from empirical means, and this can be extended to non linear estimators 
by applying the delta method, argmax theorems or plug-in approaches as 
in the classical setting -- see~\cite{VanWell,VanVaart}.\smallskip

\noindent\textbf{Organization of the paper.} In Section~\ref{AuxInfView} we
relate the raking-ratio problem to nonparametric auxiliary information.
The raking-ratio empirical process $\alpha_{n}^{(N)}(\mathcal{F})$ is defined in Section~\ref{InfNPart}.
Usual assumptions on an indexing class $\mathcal{F}$ of functions are
given in Section~\ref{RREmpGaussProc}. In Section~\ref{GenPro} we state our results for 
$\alpha_{n}^{(N)}(\mathcal{F})$ when the number $N$ of iterations is fixed. 
Our main theorem is a nonasymptotic strong approximation bound which yields the uniform
central limit theorem with rate as $n\rightarrow+\infty$, as well as an uniform
control of the bias and the covariances for fixed $n$. The approximating
Gaussian process is studied in Section~\ref{LimGau}, which establishes the uniform risk
reduction phenomenon provided the iterations are stopped properly. In Section~\ref{TwoMar}, in the two partitions case we characterize explicitly the limiting process as $N\rightarrow+\infty$.
All statements are proved in Sections~\ref{proofsGR} and \ref{proofsLP}.
The~\nameref{app} provides a few examples.

\subsection{An auxiliary information viewpoint} \label{AuxInfView}

Let $X_{1},...,X_{n}$ be
independent random variables with unknown law $P$ on some measurable space
$(\mathcal{X},\mathcal{A})$. Assumptions like separability or Haussdorf property are not necessary for this space. Let $\delta_x$ denote the Dirac mass at $x\in \mathcal{X}$ and consider the empirical measure $\mathbb{P}_{n}=n^{-1}\sum_{i=1}^{n} \delta_{X_{i}}$ on $\mathcal{A}$.\smallskip

\noindent\textbf{Auxiliary information.} Our interest for the raking-ratio method came while investigating 
how to exploit various kinds of partial information on $P$ to make $\mathbb{P}_{n}$ closer to $P$. 
The auxiliary information paradigm is as follows. Usually what is assumed on $P$ is formulated in terms of technical or regularity requirements.
Sometimes it is relevant to assume that $P$ satisfies simple properties that could be tested or estimated separately. Consider the following two extreme situations. First, a parametric model provides a tremendous amount
of information by specifying $P=P_{\theta}$ up to a finite dimensional
parameter $\theta$, so that $\mathbb{P}_{n}$ can be replaced with the most
likely $P_{\theta_{n}(X_{1},...,X_{n})}$ among the model. Notice that
$\mathbb{P}_{n}$ is used to minimize the empirical likelihood, but the
resulting $P_{\theta_{n}(X_{1},...,X_{n})}$ is of a very different nature, far
from the initial and discrete $\mathbb{P}_{n}$ thanks to the valuable parametric information.
On the opposite, in a purely nonparametric setting the information mainly comes from the sample itself, so
that only slight modifications of $\mathbb{P}_{n}$ are induced by weak
hypotheses on $P$ -- like support, regularity, symmetry, logconcavity, 
Bayesian model, semi-parametric model, etc.\ In between, we would like 
to formalize a notion of a priori auxiliary information on $P$ based on partial but concrete 
clues to be combined with the knowledge of $\mathbb{P}_{n}$. Such clues may come 
from experts, models, former inference, statistical learning or distributed data. A generic situation one can start with is when
the probabilities $P(A_{j})$ of a finite number of sets $A_{j}\in\mathcal{A}$ are known -- which
in a parametric setting already determines $\theta$ then $P$.\smallskip

\noindent\textbf{Information from one partition.} If the $A_{j}$ form a finite partition 
of $\mathcal{X}$ then the auxiliary information coincides with one discrete marginal distribution and a natural nonparametric redesign 
$\mathbb{P}_{n}^{(1)}$ of $\mathbb{P}_{n}$ is the following. Let $A_{1}^{(1)}%
,\dots,A_{m_{1}}^{(1)}\subset\mathcal{A}$ be a partition of $\mathcal{X}$ such
that $P(\mathcal{A}^{(1)})=(P(A_{1}^{(1)}),\dots,P(A_{m_{1}}^{(1)}))$ is
known. According to Proposition \ref{Kullback} below, the random measure%
\begin{equation}
\mathbb{P}_{n}^{(1)}=\frac{1}{n}%
{\displaystyle\sum\limits_{j=1}^{m_{1}}}
\frac{P(A_{j}^{(1)})}{\mathbb{P}_{n}(A_{j}^{(1)})}%
{\displaystyle\sum\limits_{X_{i}\in A_{j}^{(1)}}}
\delta_{X_{i}},\label{Pn1}%
\end{equation}
satisfies the auxiliary information $ \mathbb{P}_n^{(1)}(A_j^{(1)}) = P(A_j^{(1)}) $, for $ 1 \leqslant j \leqslant m_1 $, and is the relative entropy projection of $\mathbb{P}_{n}$ on these $m_1$ constraints. The random ratios in (\ref{Pn1}) induce a bias between $\mathbb{P}_{n}^{(1)}$ and $P$. We prove that the bias of $\alpha_n^{(1)}=\sqrt{n}(\mathbb{P}_{n}^{(1)}-P)$ vanishes uniformly and that the limiting Gaussian process of $\alpha_n^{(1)}$ has a smaller variance than the $P$-Brownian bridge.\smallskip

\noindent\textbf{Extension to }$N$\textbf{\ partitions.} If some among the sets $A_{j}$ are overlapping then the information comes from several marginal partitions. It is not obvious how to optimally combine these sources of information since there is no explicit modification
of $\mathbb{P}_{n}$ matching simultaneously several finite discrete marginals. In other words there is no closed form expression of the relative entropy projection of $\mathbb{P}_{n}$ on several margin constraints.
An alternative consists in recursively updating the current modification
$\mathbb{P}_{n}^{(N-1)}$ of $\mathbb{P}_{n}$ onto $\mathbb{P}_{n}^{(N)}$ according to the next known
marginal $P(\mathcal{A}^{(N)})=(P(A_{1}^{(N)}),\dots,P(A_{m_{N}}^{(N)}))$
exactly as in (\ref{Pn1}) for $\mathbb{P}_{n}^{(1)}$ from $\mathbb{P}%
_{n}^{(0)}=\mathbb{P}_{n}$ and $P(\mathcal{A}^{(1)})$. This coincides with the Deming and Stephan's
iterative procedure, that is the raking-ratio algorithm, as formalized in Section~\ref{InfNPart}.

\subsection{Information from N finite partitions} \label{InfNPart}

\textbf{The raking-ratio empirical measure.} For all $N\in\mathbb{N}_{\ast}$ let $m_{N}\geqslant2$ and
$\mathcal{A}^{(N)}=\{A_{1}^{(N)},\dots,A_{m_{N}}^{(N)}\}\subset\mathcal{A}$ be
a partition of $\mathcal{X}$ for which we are given the auxiliary information
$P(\mathcal{A}^{(N)})=(P(A_{1}^{(N)}),\dots,P(A_{m_{N}}^{(N)}))$ to be
exploited. Assume that
\begin{equation}
p_{N}=\min_{1\leqslant j\leqslant m_{N}}P(A_{j}^{(N)})>0,\quad
N\in\mathbb{N}_{\ast},\label{deltaN}%
\end{equation}
and $\mathcal{A}^{(N_{1})}\neq\mathcal{A}^{(N_{2})}$ if $\left\vert
N_{1}-N_{2}\right\vert =1$, otherwise $\mathcal{A}^{(N_{1})}=\mathcal{A}%
^{(N_{2})}$ is allowed. For $N=0$ there is no information and
$m_{0}=1$, $\mathcal{A}^{(0)}=\{\mathcal{X}\}$, $P(\mathcal{A}^{(0)})=\{1\}$, $p_{0}=1$. For any measurable real function $f$ write
$\mathbb{P}_{n}^{(0)}(f)=\mathbb{P}_{n}(f)=n^{-1}\sum_{i=1}^{n}f(X_{i})$,
$P(f)=\int_{\mathcal{X}}fdP$ and $\alpha_{n}^{(0)}(f)=\sqrt{n}(\mathbb{P}_{n}^{(0)}(f)-P(f))$. In (\ref{Pn1}) $\mathbb{P}_{n}^{(1)}$
allocates the random weight $P(A_{j}^{(1)})/n\mathbb{P}_{n}%
(A_{j}^{(1)})$ to each $X_{i}\in A_{j}^{(1)}$. Hence%
\begin{align*}
\mathbb{P}_{n}^{(1)}(f)  & =\sum_{i=1}^{n}\mathbb{P}_{n}^{(1)}(\left\{
X_{i}\right\}  )f(X_{i})=\sum_{j=1}^{m_{1}}%
{\displaystyle\sum\limits_{X_{i}\in A_{j}^{(1)}}}
\frac{P(A_{j}^{(1)})}{n\mathbb{P}_{n}(A_{j}^{(1)})}f(X_{i})\\
& =\sum_{j=1}^{m_{1}}\frac{P(A_{j}^{(1)})}{\mathbb{P}_{n}(A_{j}^{(1)})}\left(
\frac{1}{n}%
{\displaystyle\sum\limits_{X_{i}\in A_{j}^{(1)}}}
f(X_{i})\right)  =\sum_{j=1}^{m_{1}}\frac{P(A_{j}^{(1)})}{\mathbb{P}_{n}%
^{(0)}(A_{j}^{(1)})}\mathbb{P}_{n}^{(0)}(f1_{A_{j}^{(1)}}).
\end{align*}
Let define recursively, for $N\in\mathbb{N}_{\ast}$, the $N$-th raking-ratio
empirical measure%
\begin{equation}
\mathbb{P}_{n}^{(N)}(f)=\sum_{j=1}^{m_{N}}\frac{P(A_{j}^{(N)})}{\mathbb{P}%
_{n}^{(N-1)}(A_{j}^{(N)})}\mathbb{P}_{n}^{(N-1)}(f1_{A_{j}^{(N)}}),\label{PnN}%
\end{equation}
and the $N$-th raking-ratio empirical process%
\begin{equation}
\alpha_{n}^{(N)}(f)=\sqrt{n}(\mathbb{P}_{n}^{(N)}(f)-P(f)). \label{empnN}
\end{equation}
For $A \in \mathcal{A}$ we also write $\alpha_{n}^{(N)}(A)=\alpha_{n}^{(N)}(1_{A})$. By (\ref{PnN}) we have for all $N\in\mathbb{N}_{\ast}$%
\begin{equation}
\mathbb{P}_{n}^{(N)}(A_{j}^{(N)})=P(A_{j}^{(N)}),\quad\alpha_{n}%
^{(N)}(A_{j}^{(N)})=0,\quad1\leqslant j\leqslant m_{N},\label{annul}%
\end{equation}
as desired. Both weights and support $\{X_1,\dots,X_n\}$ of the discrete probability measure $\mathbb{P}%
_{n}^{(N)}$ are random since (\ref{PnN}) also reads%
\[
\mathbb{P}_{n}^{(N)}(\left\{  X_{i}\right\}  )=\mathbb{P}_{n}^{(N-1)}(\left\{
X_{i}\right\}  )\frac{P(A_{j}^{(N)})}{\mathbb{P}_{n}^{(N-1)}(A_{j}^{(N)}%
)},\quad\text{for }X_{i}\in A_{j}^{(N)}\text{.}%
\]
A\ few more formulas concerning $\alpha_{n}^{(N)}$
and $\mathbb{P}_{n}^{(N)}$ are derived in Section~\ref{RF}.\smallskip

\noindent\textbf{Iterated Kullback projections.} The random discrete measures
$\mathbb{P}_{n}^{(1)},...,\mathbb{P}_{n}^{(N)}$ are well defined provided that%
\begin{equation}
\min_{1\leqslant k\leqslant N}\min_{1\leqslant j\leqslant m_{k}}\mathbb{P}%
_{n}(A_{j}^{(k)})>0,\label{PnA}%
\end{equation}
which almost surely holds for all $n$ large enough and $N$ fixed, by
(\ref{deltaN}) and the law of large numbers. Given two
probability measures $Q_{n}$ and $Q$ supported by $\left\{  X_{1}%
,...,X_{n}\right\}  $ we define the relative
entropy of $ Q_n $ and $ Q $ -- see for instance~\cite{CovThom} -- to be %
\[
d_{K}(Q_{n}\mid\mid Q)=\sum_{i=1}^{n}Q_{n}(\left\{  X_{i}\right\}
)\log\left(  \frac{Q_{n}(\left\{  X_{i}\right\}  )}{Q(\left\{  X_{i}\right\}
)}\right)  .
\]

\begin{pro}
\label{Kullback}If (\ref{PnA}) holds then%
\[
\mathbb{P}_{n}^{(N)}=\arg\min\left\{  d_{K}(\mathbb{P}_{n}^{(N-1)}\mid\mid
Q):Q(\mathcal{A}^{(N)})=P(\mathcal{A}^{(N)}),\text{ supp}(Q)=\left\{
X_{1},...,X_{n}\right\}  \right\}  .
\]

\end{pro}

\noindent As a consequence, the formula (\ref{PnN}) means that the $N$-th iteration $\mathbb{P}_{n}^{(N)}$ is the Shannon-Kullback-Leibler projection of $\mathbb{P}_{n}^{(N-1)}$ under the constraint $P(\mathcal{A}^{(N)})$. Therefore the raking-ratio method is an iterated maximum likelihood procedure.\smallskip

\noindent\textbf{A mixture of conditional empirical processes.} By introducing,
for $A\in\mathcal{A}$ such that $P(A)>0$ and $\mathbb{P}_{n}^{(N)}(A)>0$,
the conditional expectations%
\begin{equation}
\mathbb{E}_{n}^{(N)}(f|A)=\frac{\mathbb{P}_{n}^{(N)}(f1_{A})}{\mathbb{P}%
_{n}^{(N)}(A)},\quad\mathbb{E}(f|A)=\frac{P(f1_{A})}{P(A)},\label{EfA}%
\end{equation}
we see that (\ref{PnN}) further reads%
\[
\mathbb{P}_{n}^{(N)}(f)=\sum_{j=1}^{m_{N}}P(A_{j}^{(N)})\mathbb{E}_{n}%
^{(N-1)}(f|A_{j}^{(N)}),\quad P(f)=\sum_{j=1}^{m_{N}}P(A_{j}^{(N)}%
)\mathbb{E}(f|A_{j}^{(N)}).
\]
Therefore (\ref{empnN}) can also be formulated into%
\begin{align}
\alpha_{n}^{(N)}(f)&=\sum_{j=1}^{m_{N}}P(A_{j}^{(N)})\alpha_{n,j}%
^{(N-1)}(f), \notag\\
\alpha_{n,j}^{(N-1)}(f)&=\sqrt{n}\left(  \mathbb{E}%
_{n}^{(N-1)}(f|A_{j}^{(N)})-\mathbb{E}(f|A_{j}^{(N)})\right)
.\label{condalphanN}%
\end{align}
Each $\alpha_{n,j}^{(N-1)}$ is the conditional empirical process of
$\mathbb{P}_{n}^{(N-1)}$ on a set $A_{j}^{(N)}$ of the new partition
$\mathcal{A}^{(N)}$. Their mixture with weights $P(\mathcal{A}^{(N)})$ is
$\alpha_{n}^{(N)}$. In view of (\ref{EfA}) and (\ref{condalphanN}) we have to
study the consequences of (\ref{annul}) on $\mathbb{P}_{n}^{(N-1)}%
(f1_{A_{j}^{(N)}})$ and $\mathbb{E}_{n}^{(N-1)}(f|A_{j}^{(N)})$ as
$n\rightarrow+\infty$, for $f\neq1_{A_{j}^{(N)}}$.\smallskip

\noindent\textbf{Bias and variance problem.} The processes $\alpha_{n,j}%
^{(N-1)}$ from (\ref{condalphanN}) are not centered due to the factors $1/\mathbb{P}_{n}^{(N-1)}%
(A_{j}^{(N)})$ in (\ref{PnN}) and
(\ref{EfA}). In general it holds
\[
\mathbb{E}\left(  \mathbb{E}_{n}^{(N-1)}(f|A)-\mathbb{E}(f|A)\right)
=\mathbb{E}\left(  \mathbb{P}_{n}^{(N-1)}(f1_{A})\left(  \frac{1}%
{\mathbb{P}_{n}^{(N-1)}(A)}-\frac{1}{P(A)}\right)  \right)  \neq0,
\]
except for $(A,f)=(A_{j}^{(N-1)},1_{A_{j}^{(N-1)}})$ hence
$\alpha_{n}^{(N)}$
is no more centered if $N\geqslant1$. This unavoidable bias is induced by
(\ref{annul}) to globally compensate for the local cancellation of the
variance of $\mathbb{P}_{n}^{(N)}(\mathcal{A}^{(N)})=P(\mathcal{A}^{(N)})$. The bias tends to spread
through (\ref{PnN}) since the information $P(\mathcal{A}^{(N)})$ is
applied to the biased $\mathbb{P}_{n}^{(N-1)}$ instead of the unbiased
$\mathbb{P}_{n}$. The variance of $\mathbb{P}_{n}^{(N)}(f)$ for the step functions 
$f=1_{A}$ being null if $A\in\mathcal{A}^{(N)}$ one expects that 
$\mathbb{V}(\mathbb{\alpha}_{n}^{(N)}(f))\leqslant\mathbb{V}(\mathbb{\alpha}_{n}^{(0)}(f))$ 
for many more functions $f$. Our results show that, uniformly over a large class of functions,
the bias vanishes asymptotically and the variance decreases, as well as the quadratic risk, thus $\mathbb{E}((\mathbb{P}_{n}^{(N)}(f)-P(f))^{2})\leqslant \mathbb{E}%
((\mathbb{P}_{n}^{(0)}(f)-P(f))^{2})$ for $n$ large. \smallskip

\section{Main results}

\subsection{The raking-ratio empirical and Gaussian processes} \label{RREmpGaussProc}

Let $\mathcal{M}$ denote the set of measurable real valued functions on
$(\mathcal{X},\mathcal{A})$. Consider a class $\mathcal{F}\subset\mathcal{M}$
such that $\sup_{f\in\mathcal{F}}|f|\leqslant M<+\infty$ and satisfying the
pointwise measurability condition often used to avoid measurability problems.
Namely, $\lim_{m\rightarrow+\infty}f_{m}(x)=f(x)$ for all $x\in\mathcal{X}$
and $f\in\mathcal{F}$ where $\left\{  f_{m}\right\}  \subset\mathcal{F}_{*}$
depends on $f$ and $\mathcal{F}_{*}\subset\mathcal{F}$ is countable. With no
loss of generality also assume that
\begin{equation}
1_{A}f\in\mathcal{F},\quad A\in\mathcal{A}_{\cup}^{(N)}=\mathcal{A}%
^{(1)}\cup...\cup\mathcal{A}^{(N)},\quad f\in\mathcal{F}.\label{1A}%
\end{equation}
In addition $\mathcal{F}$ is assumed to have either a small
uniform entropy, like Vapnik-Chervonenkis classes -- for short, VC-classes -- or a small $P$-bracketing entropy, like many classes of smooth functions. These entropy conditions are defined below and named (VC) and (BR) respectively. For a probability measure $Q$ 
on $(\mathcal{X},\mathcal{A})$ and
$f,g\in\mathcal{M}$ define $d_{Q}^{2}(f,g)=\int_{\mathcal{X}}(f-g)^{2}dQ$. Let
$N(\mathcal{F},\varepsilon,d_{Q})$ be the minimum number of balls having
$d_{Q}$-radius $\varepsilon$ needed to cover $\mathcal{F}$. Let $N_{[\ ]}%
(\mathcal{F},\varepsilon,d_{P})$ be the least number of $\varepsilon$-brackets
necessary to cover $\mathcal{F}$, of the form
$\left[  g_{-},g_{+}\right]  =\left\{  f:g_{-}\leqslant f\leqslant
g_{+}\right\}  $ with $d_{P}(g_{-},g_{+})<\varepsilon$.

\begin{hyp}
[VC]\label{VC}For $c_{0}>0$, $\nu_{0}>0$ it holds $\sup
_{Q}N\left(  \mathcal{F},\varepsilon,d_{Q}\right)  \leqslant c_{0}%
/\varepsilon^{\nu_{0}} $ where the supremum is taken over all discrete
probability measures $Q $ on $(\mathcal{X},\mathcal{A})$.
\end{hyp}

\begin{hyp}
[BR]\label{BR}For $b_{0}>0$, $r_{0}\in\left(0,1\right)$ it holds
$N_{[\ ]}(\mathcal{F},\varepsilon,d_{P})\leqslant\exp(b_{0}%
^{2}/\varepsilon^{2r_{0}})  $.
\end{hyp}

\noindent If one modifies a class $ \mathcal{F} $ satisfying (VC) or (BR) by adding functions necessary to also satisfy the condition (\ref{1A}) then (VC) or (BR) still holds with a new constant $c_{0}$ or $b_{0}$ respectively. Many properties and examples of VC-classes or classes satisfying (BR) can be found in Pollard~\cite{Poll84}, Van der Vaart and Wellner~\cite{VanWell} or Dudley~\cite{DudUnif}. Uniform boundedness is the less crucial assumption and could be replaced by a moment condition allowing some truncation arguments, however adding technicalities.\smallskip

\noindent Let $\ell^{\infty}(\mathcal{F})$ denote the set of real-valued
functions bounded on $\mathcal{F}$, endowed with the supremum norm $\left\Vert
\mathbb{\cdot}\right\Vert _{\mathcal{F}}$. The raking-ratio empirical process
$\alpha_{n}^{(N)}$ defined at (\ref{empnN}) is now denoted %
$\alpha_{n}^{(N)}(\mathcal{F})=\{\alpha_{n}^{(N)}(f):f\in\mathcal{F} \}$.
Under $(VC)$ or $(BR)$ $\mathcal{F}$ is a $P$-Donsker class -- see Sections 2.5.1 and 2.5.2 of \cite{VanWell}. Thus $\alpha
_{n}^{(0)}(\mathcal{F})$ converges weakly in $\ell^{\infty}(\mathcal{F})$ to
the $P$-Brownian bridge $\mathbb{G}$ indexed by $\mathcal{F}$, that we denote
$\mathbb{G}(\mathcal{F})=\left\{  \mathbb{G}(f):f\in\mathcal{F}\right\}  $.
Hence $\mathbb{G}(\mathcal{F})$ is a Gaussian process such that $ f \mapsto \mathbb{G}(f) $ is linear and, for any $f,g\in\mathcal{F}$,%
\begin{equation}
\mathbb{E}\left(  \mathbb{G}(f)\right)  =0,\quad \mathrm{Cov}(\mathbb{G}(f),\mathbb{G}%
(g))=P(fg)-P(f)P(g).\label{GF}%
\end{equation}
As for $\alpha_{n}^{(N)}$ we write $\mathbb{G}^{(0)}(\mathcal{F})=\mathbb{G}(\mathcal{F})$ and, for short, $\mathbb{G}(A)=\mathbb{G}(1_{A})$ if
$A\in\mathcal{A}$. Remind (\ref{EfA}). Let us introduce a new centered Gaussian process $\mathbb{G}^{(N)}(\mathcal{F})$ indexed by $\mathcal{F}$ that we call the $N$-th raking-ratio $P$-Brownian bridge and that is defined recursively, for any $N\in \mathbb{N}_{\ast}$ and $f\in\mathcal{F}$, by%
\begin{equation}
\mathbb{G}^{(N)}(f)=\mathbb{G}^{(N-1)}(f)-\sum_{j=1}^{m_{N}}\mathbb{E(}%
f|A_{j}^{(N)})\mathbb{G}^{(N-1)}(A_{j}^{(N)})\label{GNf}.%
\end{equation}
The distribution of $ \mathbb{G}^{(N)} $ is given in Proposition \ref{Thm2}. Lastly, the following notation will be useful,
\begin{equation}
\sigma^{2}_{f}=\mathbb{V}(f(X))=P(f^{2})-P(f)^{2}, \quad\sigma_{\mathcal{F}}^{2}=\sup_{f\in\mathcal{F}}\sigma^{2}_{f}.\label{sigma}
\end{equation}
Notice that $\sigma^{2}_{f}=\mathbb{V}(\alpha_{n}^{(0)}(f))=\mathbb{V}(\mathbb{G}^{(0)}
(f))$.

\subsection{General properties}\label{GenPro}

We now state asymptotic and nonasymptotic properties that always hold after raking $N_{0}$ times.
The i.i.d. sequence $\left\{  X_{n}\right\}  $ is defined on a
probability space $(\Omega,\mathcal{T},\mathbb{P})$ so that $\mathbb{P}$
implicitly leads all convergences when $n\rightarrow+\infty$ and
$(\mathcal{X},\mathcal{A})$ is endowed
with $P=\mathbb{P}^{X_{1}}$. For all $N\leqslant N_{0}$ the information
$P(\mathcal{A}^{(N)})$ satisfies (\ref{deltaN}). Most of the subsequent
constants can be bounded by using only $N_{0}$ and%
\begin{equation}
p_{(N_{0})}=\min_{0\leqslant N\leqslant N_{0}}p_{N}=\min_{0\leqslant
N\leqslant N_{0}}\min_{1\leqslant j\leqslant m_{N}}P(A_{j}^{(N)}%
)>0.\label{deltaN0}%
\end{equation}
Write $L(x)=\log(\max(e,x))$ and define $\kappa_{N_{0}}=\prod_{N=1}^{N_{0}}(1+Mm_{N})$.

\begin{pro}
\label{ProLLIMes}If $\mathcal{F}$ satisfies (VC) or (BR) then for all
$N_{0}\in\mathbb{N}$ it holds
\[
\underset{n\rightarrow+\infty}{\lim\sup}\frac{1}{\sqrt{2L\circ L(n)}}%
\sup_{0\leqslant N\leqslant N_{0}}\left\Vert \alpha_{n}^{(N)}\right\Vert
_{\mathcal{F}}\leqslant\kappa_{N_{0}}\sigma_{\mathcal{F}}\quad a.s.
\]
\end{pro}

\noindent\textit{Remark E.} The limiting constant 
$\kappa_{N_{0}}\leqslant\left(  1+M/p_{(N_{0})}\right)^{N_{0}}$
is large, and possibly largely suboptimal, except for $N_{0}=0$ where $\kappa_{0}=1$
coincides with the classical law of the iterated logarithm -- from 
which the proposition follows.\smallskip

\noindent The next result shows that the nonasymptotic deviation probability for $\Vert \alpha_{n}^{(N)}\Vert _{\mathcal{F}}$ can be controlled by the deviation probability of $\Vert \alpha_{n}^{(0)}\Vert _{\mathcal{F}}$ which in turn can be bounded by using Talagrand \cite{Tal}, van der Vaart and Wellner ~\cite{VanWell} or more recent bounds from empirical processes theory. However, since the partition changes at each
step the constants are penalized by factors similar to $\kappa_{N_{0}}$ above, involving
\begin{equation}
P_{N_{0}}=\prod_{N=1}^{N_{0}}p_{N},\quad M_{N_{0}}=
\prod_{N=1}^{N_{0}}m_{N},\quad S_{N_{0}}=\sum_{N=1}^{N_{0}}m_{N}.\label{deltaMN}%
\end{equation}

\begin{pro}
\label{ProConcentr}If $\mathcal{F}$ is pointwise measurable, bounded by $M$
then for any $N_{0}\in\mathbb{N}$, any $n\in\mathbb{N}_{\ast} $ and any $ \lambda > 0 $ we have%
\begin{align*}
\mathbb{P}\left(  \sup_{0\leqslant N\leqslant N_{0}}\left\Vert \alpha
_{n}^{(N)}\right\Vert _{\mathcal{F}}\geqslant\lambda\right)   &
\leqslant2^{N_{0}} N_0 M_{N_{0}}\mathbb{P}\left(  \left\Vert \alpha_{n}%
^{(0)}\right\Vert _{\mathcal{F}}\geqslant\frac{\lambda P_{N_{0}}%
}{(1+M+\lambda/\sqrt{n})^{N_{0}}}\right)  \\
&  \quad+S_{N_{0}}\left(  1-p_{(N_{0})}\right)  ^{n}.
\end{align*}
Under (BR) it holds, for $n>n_{0}$ and $\lambda_{0}<\lambda<D_{0}\sqrt{n}$,
\begin{align*}
 \mathbb{P}\left(  \sup_{0\leqslant N\leqslant N_{0}}\left\Vert \alpha
_{n}^{(N)}\right\Vert _{\mathcal{F}}\geqslant\lambda\right)
\leqslant D_1\exp(-D_2 \lambda^2)+S_{N_{0}}\left(  1-p_{(N_{0})}\right)  ^{n},%
 \end{align*}
where the positive constants $D_{0}$, $D_{1}$, $D_{2}$, $n_0$, $\lambda_0$ are defined at (\ref{D012}).
Under (VC) it holds, for $n>n_{0}$ and $\lambda
_{0}<\lambda<2M\sqrt{n}$,%
\begin{align*}
\mathbb{P}\left(  \sup_{0\leqslant N\leqslant N_{0}}\left\Vert \alpha
_{n}^{(N)}\right\Vert _{\mathcal{F}}\geqslant\lambda\right) 
\leqslant D_{3}\lambda^{v_{0}}\exp(-D_{4}\lambda^{2})+S_{N_{0}}\left(  1-p_{(N_{0})}\right)^{n},%
\end{align*}
where the positive constants $D_{3}$, $D_{4}$, $n_0$, $\lambda_0$ are defined at (\ref{D34}).
\end{pro}

\noindent\textit{Remark F.} Clearly, to avoid drawbacks $N_{0}$ should be fixed as $n$ increases, and
$\mathcal{F}$ limited to the bare necessities for the actual statistical problem. In this case, Proposition \ref{ProConcentr} shows that $\Vert
\alpha_{n}^{(N_{0})}\Vert _{\mathcal{F}}$ is of order $C\sqrt{\log n}$
with probability less than $ 1/n^2 $ and $C>0$. Concentration of measure type probability bounds for $\Vert \alpha_{n}^{(N_{0})}\Vert _{\mathcal{F}}-\mathbb{E}(\Vert \alpha_{n}^{(N_{0})}\Vert _{\mathcal{F}})$  are more difficult to handle due to the mixture (\ref{condalphanN}) of processes 
$\alpha_{n,j}^{(N-1)}$ involving unbounded random coefficients.\smallskip

\noindent Our main result is that the raking-ratio empirical processes $\alpha_{n}^{(0)},...,\alpha_{n}^{(N_0)}$
jointly converge weakly at some explicit rate to the
raking-ratio $P$-Brownian bridges $\mathbb{G}^{(0)},..., \mathbb{G}^{(N_0)}$
defined at (\ref{GNf}) and studied in Section~\ref{LimGau}. 
The $\mathbb{R}^{N_{0}+1}$-valued version can be stated as follows.

\begin{pro}
\label{ProCVProc}If $\mathcal{F}$ satisfies (VC) or (BR) then for all
$N_{0}\in\mathbb{N}$, as $n\rightarrow+\infty$ the sequence 
$(\alpha_{n}^{(0)}(\mathcal{F}),...,\alpha_{n}^{(N_{0})}(\mathcal{F}))$ 
converges weakly to 
$(\mathbb{G}^{(0)}(\mathcal{F}),...,\mathbb{G}^{(N_{0})}(\mathcal{F}))$
on $\ell_{\infty}(\mathcal{F}\rightarrow\mathbb{R}^{N_{0}+1})$.
\end{pro}

\noindent By using Berthet and Mason \cite{BerMas06} we further obtain the following upper
bound for the speed of Gaussian approximation of $ \alpha_n^{(N)} $ in $\left\Vert
\mathbb{\cdot}\right\Vert _{\mathcal{F}}$ distance. The powers provided at their Propositions 1 and 2 are $\alpha=1/(2+5\nu_0)$, $\beta=(4+5\nu_0)/(4+10\nu_0)$ and $\gamma=(1-r_0)/2r_0$ -- they could be slightly improved.

\begin{thm}
\label{ProCVSpeed} Let $\theta_{0}>0$. If $\mathcal{F}$ satisfies (VC) then write $v_{n}%
=(\log n)^{\beta}/n^{\alpha}$. If $\mathcal{F}$ satisfies (BR) then write
$v_{n}=1/(\log n)^{\gamma}$. In both cases, one can define on the same
probability space $(\Omega,\mathcal{T},\mathbb{P})$ a sequence $\left\{
X_{n}\right\}  $ of independent random variables with law $P$ and a sequence
$\left\{  \mathbb{G}_{n}\right\}  $ of versions of $\mathbb{G}$ satisfying the
following property. For any $N_{0}\geqslant0$ there exists $n_{0}\in
\mathbb{N}$ and $d_{0}>0$ such that we have, for all $n\geqslant n_{0}$,
\[
\mathbb{P}\left(  \sup_{0\leqslant N\leqslant N_{0}}\left\Vert \alpha
_{n}^{(N)}-\mathbb{G}_{n}^{(N)}\right\Vert _{\mathcal{F}}\geqslant d_{0}%
v_{n}\right)  <\frac{1}{n^{\theta_{0}}},
\]
where $\mathbb{G}_{n}^{(N)}$ is the version of $\mathbb{G}^{(N)}$ derived from
$\mathbb{G}_{n}^{(0)}=\mathbb{G}_{n}$ through (\ref{GNf}).
\end{thm}

\noindent\textit{Remark G.} Applied with $\theta_{0}>1$, Theorem \ref{ProCVSpeed} makes the 
study of weak convergence of functions of $\alpha_{n}^{(N)}(\mathcal{F})$ 
easier by substituting $\mathbb{G}_{n}^{(N)}$ to $\alpha_{n}^{(N)}$ through%
\[
\underset{n\rightarrow+\infty}{\lim\sup}\frac{1}{v_{n}}\sup_{0\leqslant
N\leqslant N_{0}}\left\Vert \alpha_{n}^{(N)}-\mathbb{G}_{n}^{(N)}\right\Vert
_{\mathcal{F}}\leqslant d_{0}<+\infty\quad a.s.
\]
then exploiting the properties induced by (\ref{GNf}) as in Section \ref{LimGau}. For instance
the finite dimensional laws of $\mathbb{G}^{(N)}$ are computed explicitly 
at Proposition \ref{Thm2}. For nonasymptotic applications, given a 
class $\mathcal{F}$ of interest it is possible to compute crude bounds 
for $n_{0}$ and $d_{0}$ since most constants are left explicit in our proofs as well as in \cite{BerMas06}. Indeed $d_{0}$ depends on  $p_{N_{0}}$  from (\ref{deltaN0}), on $P_{N_{0}}, M_{N_{0}}, S_{N_{0}}$ from (\ref{deltaMN}), on $\nu_0, c_0, r_0, b_0$ from (VC) or (BR), on $N_{0}, M,\theta_{0}$ and on some universal constants from the literature.\smallskip

\noindent Clearly, Theorem \ref{ProCVSpeed} implies that the speed of weak convergence of $\alpha_n^{(N)}$ to $\mathbb{G}^{(N)}$ in L{\'e}vy-Prokhorov distance $ d_{LP} $ is at least $ d_0 v_n $ -- see (\ref{Levy}) and Section 11.3 of \cite{Dud} for a definition of this metric. More deeply, from Theorem \ref{ProCVSpeed}\ we derive the following rates of uniform convergence for the bias and the variance.

\begin{pro}
\label{ProBias}If $\mathcal{F}$ satisfies (VC) or (BR) then for $N_{0}%
\in\mathbb{N}$ it holds%
\[
\limsup_{n\rightarrow+\infty}\frac{\sqrt{n}}{v_{n}}\max_{0\leqslant N\leqslant
N_{0}}\sup_{f\in\mathcal{F}}\left\vert \mathbb{E}\left(  \mathbb{P}_{n}%
^{(N)}(f)\right)  -P(f)\right\vert  \leqslant d_{0},
\]
where $v_{n}\rightarrow 0$ and $d_{0}$ are the same as in Theorem~\ref{ProCVSpeed}, and
\begin{align*}
& \limsup_{n\rightarrow+\infty}\frac{n}{v_{n}}\sup_{f,g\in\mathcal{F}}\left\vert
\mathbb{E}\left( (\mathbb{P}_{n}^{(N)}(f)-P(f))(\mathbb{P}_{n}^{(N)}(g)-P(g)) \right)
-\frac{1}{n}\mathrm{Cov}\left(  \mathbb{G}^{(N)}(f),\mathbb{G}^{(N)}(g)\right)\right\vert \\
& = \limsup_{n\rightarrow+\infty}\frac{n}{v_{n}}\sup_{f,g\in\mathcal{F}}\left\vert
\mathrm{Cov}\left(  \mathbb{P}_{n}^{(N)}(f),\mathbb{P}_{n}^{(N)}(g)\right)  -\frac{1}%
{n}\mathrm{Cov}\left(  \mathbb{G}^{(N)}(f),\mathbb{G}^{(N)}(g)\right)  \right\vert
\leqslant\sqrt{\frac{8}{\pi}}d_{0}\sigma_{\mathcal{F}}.
\end{align*}

\end{pro}

\noindent By Proposition~\ref{ProBias}, the bias process $
\mathbb{E}(\alpha_{n}^{(N)}(f))=\sqrt{n}(  \mathbb{E(P}_{n}%
^{(N)}(f))-Pf) $ vanishes at the uniform rate $1/\sqrt{n}$. The covariance of $\mathbb{G}^{(N)}$ is computed in Section~\ref{GenPro} and the quadratic risk is estimated at Remark I. 

\noindent A second consequence of Theorem \ref{ProCVSpeed} is uniform
Berry-Esseen type bounds. Let $\Phi$ denote the distribution function of the
centered standardized normal law.

\begin{pro}
\label{ProBerryEsseen}Assume that $\mathcal{F}$ satisfies (VC) or (BR), fix
$N_{0}\in\mathbb{N}$ and let $d_{0}>0$, $v_{n}\rightarrow0$ be defined as
in Theorem \ref{ProCVSpeed}. If $\mathcal{F}_{0}\subset\mathcal{F}$ is such that $$\sigma_{0}%
^{2}=\inf\left\{  \mathbb{V}\left(  \mathbb{G}^{(N)}(f)\right)  :f\in
\mathcal{F}_{0},0\leqslant N\leqslant N_{0}\right\}  >0,$$ then for any
$d_{1}>d_{0}$ there exists $n_{1}\in\mathbb{N}$ such that for all $n\geqslant
n_{1}$,
\begin{align}\label{BE_1st_stat}
    \max_{0\leqslant N\leqslant N_{0}}\sup_{f\in\mathcal{F}_{0}}\sup
_{x\in\mathbb{R}}\left\vert \mathbb{P}\left( \sqrt{n}\frac{\mathbb{P}%
_{n}^{(N)}(f)-P(f)}{\sqrt{\mathbb{V}\left(\mathbb{G}^{(N)}(f)\right)}%
}\leqslant x\right)  -\Phi(x)\right\vert \leqslant\frac{d_{1}}{\sqrt{2\pi
}\sigma_{0}}v_{n}.
\end{align}
Let $\mathcal{L}$ be a collection of real valued Lipschitz functions $\varphi$
defined on $\ell_{\infty}(\mathcal{F})$ with Lipschitz constant bounded by
$C_{1}<+\infty$ and such that $\varphi(\mathbb{G}^{(N)})$ has a density
bounded by $C_{2}<+\infty$ for all $0\leqslant N\leqslant N_{0}$. Then for all $\varphi\in\mathcal{L}$,
$n\geqslant n_{1}$,%
\begin{align}\label{BE_2nd_stat}
    \max_{0\leqslant N\leqslant N_{0}}\sup_{\varphi\in\mathcal{L}}\sup
_{x\in\mathbb{R}}\left\vert \mathbb{P}\left(  \varphi(\alpha_{n}%
^{(N)})\leqslant x\right)  -\mathbb{P}\left(  \varphi(\mathbb{G}%
^{(N)})\leqslant x\right)  \right\vert \leqslant d_{1}C_{1}C_{2}v_{n}.
\end{align}
\end{pro}

\noindent\textit{Remark H.} The formula (\ref{BE_1st_stat}) is a special case of the second one (\ref{BE_2nd_stat}) and reads
\[
\max_{0\leqslant N\leqslant N_{0}}\sup_{f\in\mathcal{F}_{0}}\sup
_{x\in\mathbb{R}}\left\vert \mathbb{P}\left(  \alpha_{n}^{(N)}(f)\leqslant
x\right)  -\mathbb{P}\left(  \mathbb{G}^{(N)}(f)\leqslant x\right)
\right\vert \leqslant\frac{d_{1}}{\sqrt{2\pi}\sigma_{0}}v_{n}.
\]
The functions $f\in\mathcal{F}$ overdetermined by the knowledge of
$P(\mathcal{A}^{(N)})$ have a small $\mathbb{V}(\mathbb{G}^{(N)}(f))$
and are excluded from $\mathcal{F}_{0}$. Proposition 
\ref{ProBerryEsseen} is especially useful under (VC) since $v_{n}
$ is then polynomialy decreasing, thus allowing larger $C_{1}C_{2}$ and
$\mathcal{L}$. An example is given in Section~\ref{AppEx}. Whenever the class $\mathcal{F}$ is finite, the density of the transform $ \varphi(\mathbb{G}(\mathcal{F})) $ of the finite dimensional Gaussian vector $\mathbb{G}(\mathcal{F})$ is easily computed. The conditions for (\ref{BE_2nd_stat}) of Proposition \ref{ProBerryEsseen} are fulfilled if, for example, all random variables $ \varphi(\mathbb{G}(\mathcal{F})) $ can be controlled by discretizing the small entropy class $\mathcal{F}$, by bounding their densities then by taking limits accordingly.\smallskip

\subsection{Limiting variance and risk reduction}\label{LimGau}

In this section we study the covariance structure of $\mathbb{G}%
^{(N)}(\mathcal{F})$ from \eqref{GNf}, for $N$ fixed. The following matrix
notation is introduced to shorten formulas. The brackets $[\cdot]$ refer to
column vectors built from the partition $\mathcal{A}^{(k)}$ appearing inside. Let
$V^{t}$ denote the transpose of a vector $V$. For $k\leqslant N$ write%
\[
\mathbb{E}\left[f|\mathcal{A}^{(k)}\right]=\left(  \mathbb{E(}f|A_{1}^{(k)}%
),\dots,\mathbb{E}(f|A_{m_{k}}^{(k)})\right)  ^{t},\quad\mathbb{G}%
\left[\mathcal{A}^{(k)}\right]=\left(  \mathbb{G}(A_{1}^{(k)}),\dots,\mathbb{G}(A_{m_{k}%
}^{(k)})\right)  ^{t},%
\]
and, for $l\leqslant k\leqslant N$ define the matrix $\mathbf{P}_{\mathcal{A}^{(k)}%
|\mathcal{A}^{(l)}}$ to be%
\[
\left(  \mathbf{P}_{\mathcal{A}^{(k)}|\mathcal{A}^{(l)}}\right)
_{i,j}=P(A_{j}^{(k)}|A_{i}^{(l)})=\frac{P(A_{j}^{(k)}\cap A_{i}^{(l)}%
)}{P(A_{i}^{(l)})},\quad 1\leqslant i\leqslant m_{l}, 1\leqslant
j\leqslant m_{k}.
\]

\noindent Write $\mathrm{Id}_k$ the identity matrix $k\times k$. Remind that
$\mathbb{V}(\mathbb{G}(f))=P(f^{2})-(P(f))^{2}$, $P(\mathcal{A}^{(k)}%
)^{t}=P[\mathcal{A}^{(k)}]$ and
$P(A_{i}^{(k)}\cap A_{j}^{(k)})=0$ if $i\neq j$. The covariance matrix of the
Gaussian vector $\mathbb{G}[\mathcal{A}^{(k)}]$ is $\mathbb{V}(\mathbb{G}%
[\mathcal{A}^{(k)}])=diag(P(\mathcal{A}^{(k)}))-P(\mathcal{A}^{(k)}%
)^{t}P(\mathcal{A}^{(k)})$. Let $\cdot$ denote a product between a square matrix and a vector. Finally define%
\begin{align}
& \Phi_{k}^{(N)}(f) =\mathbb{E}\left[f|\mathcal{A}^{(k)}\right]+ \nonumber \\
& \sum_{\substack{1\leqslant
L\leqslant N-k\\k<l_{1}<l_{2}<...<l_{L}\leqslant N}}(-1)^{L}\mathbf{P}%
_{\mathcal{A}^{(l_{1})}|\mathcal{A}^{(i)}}\mathbf{P}_{\mathcal{A}%
^{(l_{2})}|\mathcal{A}^{(l_{1})}}\dots\mathbf{P}_{\mathcal{A}%
^{(l_{L})}|\mathcal{A}^{(l_{L-1})}}\cdot\mathbb{E}\left[f|\mathcal{A}^{(l_{L})}\right].\label{phikNf}
\end{align}
An explicit expression for $ \mathbb{G}^{(N)} $ is given in Lemma~\ref{Thm1} and the closed form for the covariance function of $\mathbb{G}%
^{(N)}(\mathcal{F})$ is as follows.

\begin{pro}
\label{Thm2}For all $N\in\mathbb{N}$ the process $\mathbb{G}^{(N)}(\mathcal{F})$ is Gaussian, centered and linear
with covariance function defined to be, for $(f,g)\in\mathcal{F}^{2}$,%
\begin{align*}
\mathrm{Cov}\left(\mathbb{G}^{(N)}(f),\mathbb{G}^{(N)}(g)\right)  &  =\mathrm{Cov}%
\left(\mathbb{G}(f),\mathbb{G}(g)\right)-\sum_{k=1}^{N}\Phi_{k}^{(N)}(f)^{t}\cdot\mathbb{V}\left(  \mathbb{G}%
[\mathcal{A}^{(k)}]\right)  \cdot\Phi_{k}^{(N)}(g).
\end{align*}

\end{pro}

\noindent Proposition \ref{Thm2} implies the following variance reduction 
phenomenon.

\begin{pro}
\label{Thm2b}For any $\left\{  f_{1},...,f_{m}\right\}
\subset\mathcal{F}$ and $N\in \mathbb{N}$ the covariance matrices $\Sigma_{m}^{(N)}=\mathbb{V}%
((\mathbb{G}^{(N)}(f_{1}),...,\mathbb{G}^{(N)}(f_{k})))$ are such that
$\Sigma_{m}^{(0)}-\Sigma_{m}^{(N)}$ is positive definite. \ 
\end{pro}

\noindent\textit{Remark I.} In particular we have 
$\mathbb{V}(\mathbb{G}^{(N)}(f))\leqslant\mathbb{V}(\mathbb{G}^{(0)}(f))=\sigma_{f}^{2}$, $f\in\mathcal{F}$. The asymptotic risk 
reduction after raking is quantified by combining Propositions \ref{ProBias} and
\ref{Thm2b}. Given $\varepsilon_{0}>0$ and $0<\sigma_{0}<\sigma_{\mathcal{F}}$ there exists
some $n_{0}=n_{0}(\varepsilon_{0},\mathcal{F})$ such that if $n>n_{0}$ then any $f\in\mathcal{F}$ with initial quadratic risk $\sigma_{f}^{2}/n>\sigma_{0}/n$ has a new risk, after raking $N$ times, equal to%
\[
\mathbb{E}\left(  (\mathbb{P}_{n}^{(N)}(f)-P(f))^{2}\right)
=\frac{\sigma_{f}^{2}}{n}(\Delta(f)+e(f)v_{n}),
\]
where $v_{n}\rightarrow 0$ and $d_{0}$ are as in Theorem~\ref{ProCVSpeed} and
\begin{align}
\Delta(f) & =\frac{\mathbb{V}(\mathbb{G}^{(N)}(f))}{\sigma_{f}^{2}}%
\in\left[  0,1\right], \nonumber \\
\sup_{f\in\mathcal{F},\text{ }\sigma_{f}\geqslant\sigma_{0}}\left\vert
e(f)\right\vert  & <(1+\varepsilon_{0})\sqrt{\frac{8}{\pi}}d_{0}\frac
{\sigma_{\mathcal{F}}}{\sigma_{0}}, \nonumber \\
\mathbb{V}\left(\mathbb{G}^{(N)}(f)\right)  &  =\sigma_{f}^{2}-\sum_{k=1}%
^{N}\Phi_{k}^{(N)}(f)^{t}\cdot\mathbb{V}\left(\mathbb{G}[\mathcal{A}^{(k)}%
]\right)\cdot\Phi_{k}^{(N)}(f)\label{VGNf},
\end{align}
so that the risk is reduced whenever $\Delta(f)<1$ and $n$ is large
enough.\smallskip

\noindent When $N_{1}>N_{0}>0$ it is not automatically true that the covariance structure of
$\alpha_{n}^{(N_{1})}(\mathcal{F})$ decreases compared to that of $\alpha_{n}^{(N_{0}%
)}(\mathcal{F})$. According to the next statement, a simple sufficient condition is to rake two
times along the same cycle of partitions.

\begin{pro}
\label{Thm3} Let $N_{0},N_{1}\in\mathbb{N}$ be such that $N_{1}\geqslant
2N_{0}$ and
\[
\mathcal{A}^{(N_{0}-k)}=\mathcal{A}^{(N_{1}-k)},\text{ for }0\leqslant
k<N_{0}.
\]
Then it holds $\mathbb{V}(\mathbb{G}^{(N_{1})}(f))\leqslant\mathbb{V}%
(\mathbb{G}^{(N_{0})}(f))$ for all $f\in\mathcal{F}$ and $\Sigma_{m}^{(N_{0}%
)}-\Sigma_{m}^{(N_{1})}$ is positive definite for all $\left\{  f_{1}%
,...,f_{m}\right\}  \subset\mathcal{F}$.
\end{pro}

\noindent\textit{Remark J.} In Appendix~\ref{AppCounterEx} a counter-example with $N_{1}=N_{0}+1$
shows that the variance does not decrease for all functions
at each iteration. This case is excluded from Proposition \ref{Thm3} 
since $N_{1}=N_{0}+1<2N_{0}$ if $N_{0}>1$ and, whenever $N_{0}=1$ and $N_{1}=2$
the requirement $\mathcal{A}^{(N_{0})}=\mathcal{A}^{(N_{1})}$ is not allowed.

\subsection{The case of two marginals}\label{TwoMar}

We now consider the original method where $k$ partitions
are raked in a periodic order. Let us focus on the case $k=2$ of the two-way
contingency table. The Deming and Stephan algorithm coincides with
the Sinkhorn-Knopp algorithm for matrix scaling~\cite{SinkhKnopp}. Denote $\mathcal{A}=\mathcal{A}%
^{(1)}=\left\{  A_{1},...,A_{m_{1}}\right\}  $ and $\mathcal{B}=\mathcal{A}%
^{(2)}=\left\{  B_{1},...,B_{m_{2}}\right\}  $ the two known margins, thus 
$\mathcal{A}^{(2m+1)}=\mathcal{A}$ and $\mathcal{A}^{(2m)}=\mathcal{B}$. Likewise
for $1\leqslant i\leqslant m_{1}$ and $1\leqslant j\leqslant m_{2}$ rewrite
$(\mathbf{P}_{\mathcal{A}|\mathcal{B}})_{i,j}=P(A_{j}|B_{i})$, $(\mathbf{P}%
_{\mathcal{B}|\mathcal{A}})_{i,j}=P(B_{j}|A_{i})$ and
\begin{align*}
\mathbb{G}[\mathcal{A}]=(\mathbb{G}(A_{1}),\dots,\mathbb{G}(A_{m_{1}}))^{t},
&  \qquad\mathbb{E}[f|\mathcal{A}]=\left(  \mathbb{E}(f|A_{1}),\dots
,\mathbb{E}(f|A_{m_{1}})\right)^{t},\\
\mathbb{G}[\mathcal{B}]=(\mathbb{G}(B_{1}),\dots,\mathbb{G}(B_{m_{2}}))^{t},
&  \qquad\mathbb{E}[f|\mathcal{B}]=\left(  \mathbb{E}(f|B_{1}),\dots
,\mathbb{E}(f|B_{m_{2}})\right)  ^{t}.
\end{align*}
The matrix $\mathbf{P}_{\mathcal{A}|\mathcal{B}}\mathbf{P}_{\mathcal{B}%
|\mathcal{A}}$ is $m_{1}\times m_{1}$ and $\mathbf{P}_{\mathcal{B}%
|\mathcal{A}}\mathbf{P}_{\mathcal{A}|\mathcal{B}}$ is $m_{2}\times m_{2}$. A
sum with a negative upper index is null, a matrix with a negative power is also
null, and a square matrix with power zero is the identity matrix. For $N\in
\mathbb{N}_{\ast}$ define%
\begin{align}
S_{1,\text{even}}^{(N)}(f) &  =\sum_{k=0}^{N}\left(\mathbf{P}_{\mathcal{B}%
|\mathcal{A}}\mathbf{P}_{\mathcal{A}|\mathcal{B}}\right)^{k}\cdot\left(\mathbb{E}%
[f|\mathcal{A}]-\mathbf{P}_{\mathcal{B}|\mathcal{A}}\cdot\mathbb{E}[f|\mathcal{B}%
]\right) \text{ is } m_{1}\times 1, \label{RakingRatioS1p} \\
S_{2,\text{odd}}^{(N)}(f) &  =\sum_{k=0}^{N}\left(\mathbf{P}_{\mathcal{A}%
|\mathcal{B}}\mathbf{P}_{\mathcal{B}|\mathcal{A}}\right)^{k}\cdot\left(\mathbb{E}%
[f|\mathcal{B}]-\mathbf{P}_{\mathcal{A}|\mathcal{B}}\cdot\mathbb{E}[f|\mathcal{A}%
]\right) \text{ is } m_{2}\times 1\label{RakingRatioS2i},\\
S_{1,\text{odd}}^{(N)}(f) &  =S_{1,\text{even}}^{(N)}+\left(\mathbf{P}%
_{\mathcal{B}|\mathcal{A}}\mathbf{P}_{\mathcal{A}|\mathcal{B}}\right)^{N+1}%
\cdot\mathbb{E}[f|\mathcal{A}] \text{ is } m_{1}\times 1\label{RakingRatioS1i},\\
S_{2,\text{even}}^{(N)}(f) &  =S_{2,\text{odd}}^{(N)}+\left(\mathbf{P}%
_{\mathcal{A}|\mathcal{B}}\mathbf{P}_{\mathcal{B}|\mathcal{A}}\right)^{N+1}%
\cdot\mathbb{E}[f|\mathcal{B}] \text{ is } m_{2}\times 1\label{RakingRatioS2p}.%
\end{align}

\begin{pro}
\label{Pro3} Let $m\in\mathbb{N}$. We have%
\begin{align}
\mathbb{G}^{(2m)}(f) &  =\mathbb{G}(f)- S_{1,\text{even}}^{(m-1)}(f)^{t}
\cdot\mathbb{G}[\mathcal{A}]- S_{2,\text{even}}^{(m-2)}(f)^{t}
\cdot\mathbb{G}[\mathcal{B}],\label{G2m}\\
\mathbb{G}^{(2m+1)}(f) & =\mathbb{G}(f)- S_{1,\text{odd}}^{(m-1)}(f)^{t}
\cdot\mathbb{G}[\mathcal{A}]- S_{2,\text{odd}}^{(m-1)}(f)^{t}
\cdot\mathbb{G}[\mathcal{B}].\label{G2m1}%
\end{align}

\end{pro}

\noindent\textit{Remark K.} The limiting process $\mathbb{G}^{(N)}$ evaluated at 
$f$ is then simply $\mathbb{G}(f)$ with a correction depending on the Gaussian vectors
$\mathbb{G}[\mathcal{A}]$ and $\mathbb{G}[\mathcal{B}]$ through the two
deterministic matrices $\mathbf{P}_{\mathcal{A}|\mathcal{B}}$ and
$\mathbf{P}_{\mathcal{B}|\mathcal{A}}$ carrying the information and operating
on the conditional expectation vectors $\mathbb{E}[f|\mathcal{A}]$ and
$\mathbb{E}[f|\mathcal{B}]$.\smallskip

\noindent The following assumption simplifies the limits and ensures a geometric rate of convergence for matrices $S_{i,\text{even}}^{(N)}$ and $S_{1,\text{odd}}^{(N)}$ as $N\rightarrow+\infty$.

\begin{hyp}
[ER]The matrices $\mathbf{P}_{\mathcal{A}|\mathcal{B}}\mathbf{P}_{\mathcal{B}%
|\mathcal{A}}$ and $\mathbf{P}_{\mathcal{B}|\mathcal{A}}\mathbf{P}%
_{\mathcal{A}|\mathcal{B}}$ are ergodic.
\end{hyp}

\noindent\textit{Remark L.} Notice that (ER) holds whenever the matrices
have strictly positive coefficients.
This is true for $\mathbf{P}_{\mathcal{A}|\mathcal{B}}\mathbf{P}_{\mathcal{B}%
|\mathcal{A}}$ if $\sum_{j=1}^{m_{2}}P(A\cap B_{j})P(B_{j}\cap A^{\prime})>0$
for all $A,A^{\prime}\in\mathcal{A}$ hence if each pair
$A,A^{\prime}\in\mathcal{A}$ is intersected by some $B\in\mathcal{B}$ with
positive probability. The latter requirement is for instance met if
$\mathcal{X}=\mathbb{R}^{d}$, $P$ has a positive density and the partitions concern
two distinct coordinates. 

\begin{pro}
\label{Pro4}Under (ER) the matrices $S_{l,\text{even}}^{(N)}(f)$ and $S_{l,\text{odd}}%
^{(N)}(f)$ for $l=1,2$ converge uniformly on $\mathcal{F}$ to
$S_{l,\text{even}}(f)$ and $S_{l,\text{odd}}(f)$ satisfying
\[
S_{1,\text{odd}}(f)=S_{1,\text{even}}(f)+P_{1}[f],\quad
S_{2,\text{even}}(f)=S_{2,\text{odd}}(f)+P_{2}[f],\]
where $P_{l}[f]=(P(f),\dots,P(f))^t$ are $m_{l}\times 1$ vectors. 
More precisely, given any vector norms $\lVert \cdot \rVert _{m_l}$ for
$l=1,2$, there exists $c_{l}>0$ and $0<\lambda_{l}<1$ such that
\begin{align*}
\sup_{f \in \mathcal{F}} \left\Vert S_{l,\text{even}}^{(N)}(f)-
S_{l,\text{even}}(f)\right\Vert _{m_l} & \leqslant c_{l}\lambda_{l}^{N}, \\
\sup_{f \in \mathcal{F}} \left\Vert S_{l,\text{odd}}^{(N)}(f)-
S_{l,\text{odd}}(f)\right\Vert _{m_l} & \leqslant 
c_{l}\lambda_{l}^{N}.
\end{align*}
\end{pro}

\noindent The main result of this section is the simple expression of the limiting process for a two partitions raking procedure. Let $d_{LP}$ denote the L\'{e}vy-Prokhorov distance. The matrices $S_{1,\text{even}}(f), S_{2,\text{odd}}(f)$ and scalars $\lambda_{1},\lambda_{2}$ are as in Proposition \ref{Pro4}.

\begin{thm}
\label{Thm4} Under (ER) the sequence $\left\{\mathbb{G}^{(N)}(\mathcal{F})\right\}$
defined at (\ref{GNf}) converges almost surely
to the centered Gaussian process $\mathbb{G}^{(\infty)}(\mathcal{F})$ defined to be
\[
\mathbb{G}^{(\infty)}(f)=\mathbb{G}(f)-S_{1,\text{even}}(f)^{t}\cdot
\mathbb{G}[\mathcal{A}]-S_{2,\text{odd}}(f)^{t}\cdot\mathbb{G}[\mathcal{B}%
],\quad f\in\mathcal{F}\text{.}%
\]
Moreover we have, for all $N$ large and $c_{3}>0$ depending on $\lambda_{1},\lambda_{2},P(\mathcal{A}),P(\mathcal{B})$, 
\[
d_{LP}(\mathbb{G}^{(N)},\mathbb{G}^{(\infty)})\leqslant
c_{3}\sqrt{N}\max(\lambda_{1},\lambda_{2})^{N/2}.
\]
\end{thm}

\noindent Theorem \ref{Thm4} may be viewed as a stochastic counterpart of the deterministic rate obtained by Franklin and Lorentz \cite{FranLor89} for the Sinkhorn algorithm. 
Mixing both approaches could strengthen the following two remarks.\smallskip

\noindent\textit{Remark M.} The matrices $\mathbf{P}_{\mathcal{A}|\mathcal{B}}$, $\mathbf{P}_{\mathcal{B}|\mathcal{A}}$ and the vectors $\mathbb{E}[f|\mathcal{A}]$, $\mathbb{E}[f|\mathcal{B}]$ are not known without additional information. They can be
estimated uniformly over $\mathcal{F}$ as $n\rightarrow+\infty$ to evaluate the distribution of $\mathbb{G}^{(N)}$ and $\mathbb{G}^{(\infty)}$, thus
giving access to adaptative tests or estimators. 
Since $\lambda_{1}$, $\lambda_{2}$ and $c_3$ are related to eigenvalues of $\mathbf{P}_{\mathcal{A}|\mathcal{B}}\mathbf{P}_{\mathcal{B}%
|\mathcal{A}}$ and $\mathbf{P}_{\mathcal{B}|\mathcal{A}}\mathbf{P}%
_{\mathcal{A}|\mathcal{B}}$ they can be estimated 
adaptively at rate $1/\sqrt{n}$ in probability. This in turn provides an evaluation of $d_{LP}(\mathbb{G}^{(N)},\mathbb{G}^{(\infty)})$. \smallskip

\noindent\textit{Remark N.} In the case of an auxiliary information reduced to $P(A), P(B)$ 
one should use $\mathcal{A}=\{ A,A^{c}\}$, $\mathcal{B}=\{ B,B^{c}\}$, estimate the missing $P(A\cap B)$ in $\mathbf{P}_{\mathcal{A}|\mathcal{B}}$, $\mathbf{P}_{\mathcal{B}|\mathcal{A}}$ and the conditional expectations on the four sets, then $S_{1,\text{even}}$, $S_{2,\text{odd}}$. 
If the probabilities of more overlapping sets are known the above characterization of the limiting process $\mathbb{G}^{(\infty)}(\mathcal{F})$ can be generalized to a recursive raking among $k$ partitions $\{A_j,A_j^{c}\}$ in the same order.


\section{Proofs of general results}\label{proofsGR}

\subsection{Raking formulas}\label{RF}

Write%
\begin{equation}
B_{n,N_{0}}=\left\{  \min_{0\leqslant N\leqslant N_{0}}\min_{1\leqslant
j\leqslant m_{N}}\mathbb{P}_{n}\left(  A_{j}^{(N)}\right)  >0\right\}
\label{biendef},%
\end{equation}
and $B_{n,N_{0}}^{c}=\Omega\setminus B_{n,N_{0}}$.
By (\ref{deltaN}), the probability that $\alpha_{n}^{(N_{0})}$ is undefined is
\[
P\left(  B_{n,N_{0}}^{c}\right)
\leqslant\ \sum_{N=1}^{N_{0}}m_{N}\left(  1-p_{N}\right)  ^{n}\leqslant
S_{N_{0}}\left(  1-p_{(N_{0})}\right)  ^{n}.
\]
On $B_{n,N_{0}}$ we have, by (\ref{PnN}) and since $\mathcal{A}^{(N)}$ is a partition,%
\begin{align}
&  \alpha_{n}^{(N)}(f)\nonumber\\
&  =\sqrt{n}(\mathbb{P}_{n}^{(N)}(f)-P(f))\nonumber\\
&  =\sqrt{n}\left(  \sum_{j=1}^{m_{N}}\frac{P(A_{j}^{(N)})}{\mathbb{P}%
_{n}^{(N-1)}(A_{j}^{(N)})}\mathbb{P}_{n}^{(N-1)}(f1_{A_{j}^{(N)}})-\sum
_{j=1}^{m_{N}}P\left(  f1_{A_{j}^{(N)}}\right)  \right)  \nonumber\\
&  =\sum_{j=1}^{m_{N}}\left(  \frac{P(A_{j}^{(N)})}{\mathbb{P}_{n}%
^{(N-1)}(A_{j}^{(N)})}\alpha_{n}^{(N-1)}(f1_{A_{j}^{(N)}})-\frac{P\left(
f1_{A_{j}^{(N)}}\right)  }{\mathbb{P}_{n}^{(N-1)}(A_{j}^{(N)})}\alpha
_{n}^{(N-1)}(A_{j}^{(N)})\right)  \nonumber\\
&  =\sum_{j=1}^{m_{N}}\frac{P(A_{j}^{(N)})}{\mathbb{P}_{n}^{(N-1)}(A_{j}%
^{(N)})}\alpha_{n}^{(N-1)}\left(  (f-\mathbb{E}(f|A_{j}^{(N)}))1_{A_{j}^{(N)}%
}\right)  .\label{alphanN}%
\end{align}
By (\ref{EfA}) this implies (\ref{condalphanN}) since for any $A\in\mathcal{A}$,%
\[
\alpha_{n}^{(N-1)}(f1_{A})-\mathbb{E}(f|A)\alpha_{n}^{(N-1)}(A)=\mathbb{P}%
_{n}^{(N-1)}(f1_{A})-\mathbb{E}(f|A)\mathbb{P}_{n}^{(N-1)}(A).
\]
Define $\mathcal{A}_{\cap}^{(0)}=\{\Omega\}$ and, for $N\geqslant1$,%
\begin{equation}
\mathcal{A}_{\cap}^{(N)}
=\left\{  A:A=A_{j_{1}}^{(1)}\cap A_{j_{2}}^{(2)}\cap...\cap A_{j_{N}%
}^{(N)},j_{k}\leqslant m_{k},k\leqslant N\right\}.\label{Ainter}%
\end{equation}
For any $A\in\mathcal{A}_{\cap}^{(N)}$, $\mathbb{P}_{n}%
^{(N)}$ associates to each $X_{i}\in A$ the weight%
\begin{equation}
\omega_{n}^{(N)}(A)=\frac{1}{n}%
{\displaystyle\prod\limits_{k=1}^{N}}
\frac{P(A_{j_{k}}^{(k)})}{\mathbb{P}_{n}^{(k-1)}(A_{j_{k}}^{(k)}%
)}.\label{exact}%
\end{equation}
The case $N=1$ yields (\ref{Pn1}). By induction on (\ref{PnN}), (\ref{exact}) and
since $\mathcal{A}_{\cap}^{(N)}$ is a refined finite partition of
$\mathcal{X}$, we get%
\begin{align*}
\mathbb{P}_{n}^{(N+1)}(f)  & =\sum_{j=1}^{m_{N+1}}\frac{P(A_{j}^{(N+1)}%
)}{\mathbb{P}_{n}^{(N)}(A_{j}^{(N+1)})}\mathbb{P}_{n}^{(N)}(f1_{A_{j}^{(N+1)}%
})\\
& =\sum_{j=1}^{m_{N+1}}\frac{P(A_{j}^{(N+1)})}{\mathbb{P}_{n}^{(N)}%
(A_{j}^{(N+1)})}\sum_{A\in\mathcal{A}_{\cap}^{(N)}}\sum_{i=1}^{n}\omega
_{n}^{(N)}(A)f(X_{i})1_{A\cap A_{j}^{(N+1)}}(X_{i})\\
& =\sum_{i=1}^{n}f(X_{i})\sum_{A\in\mathcal{A}_{\cap}^{(N)}}\sum
_{j=1}^{m_{N+1}}1_{A\cap A_{j}^{(N+1)}}(X_{i})\omega_{n}^{(N+1)}(A\cap
A_{j}^{(N+1)})\\
& =\sum_{i=1}^{n}f(X_{i})\sum_{A\in\mathcal{A}_{\cap}^{(N+1)}}1_{A}%
(X_{i})\omega_{n}^{(N+1)}(A).
\end{align*}

\subsection{Proof of Proposition~\ref{Kullback}.}

The partition $\mathcal{A}_{\cap}^{(N)}$ is defined at (\ref{Ainter}).
By using (\ref{exact}) it holds, for $N\geqslant1$,%
\begin{align}
& d_{K}\left(  \mathbb{P}_{n}^{(N-1)}\mid\mid\mathbb{P}_{n}^{(N)}\right) \nonumber\\
& =\sum_{i=1}^{n}\mathbb{P}_{n}^{(N-1)}(\left\{  X_{i}\right\}  )\log\left(
\frac{\mathbb{P}_{n}^{(N-1)}(\left\{  X_{i}\right\}  )}{\mathbb{P}_{n}%
^{(N)}(\left\{  X_{i}\right\}  )}\right) \nonumber\\
& =\sum_{i=1}^{n}\mathbb{P}_{n}^{(N-1)}(\left\{  X_{i}\right\}  )\sum
_{A\in\mathcal{A}_{\cap}^{(N-1)}}\sum_{j=1}^{m_{N}}1_{A\cap A_{j}^{(N)}}%
(X_{i})\log\left(  \frac{\omega_{n}^{(N-1)}(A)}{\omega_{n}^{(N)}(A\cap
A_{j}^{(N)})}\right) \nonumber\\
& =\sum_{j=1}^{m_{N}}\log\left(  \frac{\mathbb{P}_{n}^{(N-1)}(A_{j}^{(N)}%
)}{P(A_{j}^{(N)})}\right)  \sum_{i=1}^{n}\mathbb{P}_{n}^{(N-1)}(\left\{
X_{i}\right\}  )\sum_{A\in\mathcal{A}_{\cap}^{(N-1)}}1_{A\cap A_{j}^{(N)}%
}(X_{i})\nonumber\\
& =\sum_{j=1}^{m_{N}}\mathbb{P}_{n}^{(N-1)}(A_{j}^{(N)})\log\left(
\frac{\mathbb{P}_{n}^{(N-1)}(A_{j}^{(N)})}{P(A_{j}^{(N)})}\right)\label{DKA},
\end{align}
since $\mathcal{A}_{\cap}^{(N-1)}$ is a partition of $\mathcal{X}$.
Hence the contrast between $\mathbb{P}_{n}^{(N-1)}$ and $\mathbb{P}_{n}^{(N)}$
is the same on $\left\{  X_{1},...,X_{n}\right\}  $ as on $\mathcal{A}^{(N)}$. 
Now, by convexity of $-\log(x)$ it follows, for any probability distribution $Q$ supported by
$\left\{  X_{1},...,X_{n}\right\}  $,%
\begin{align*}
& d_{K}\left(  \mathbb{P}_{n}^{(N-1)}\mid\mid Q\right) \\
& =-\sum_{i=1}^{n}%
\log\left(  \frac{Q(\left\{  X_{i}\right\}  )}{\mathbb{P}_{n}^{(N-1)}(\left\{
X_{i}\right\}  )}\right)  \sum_{j=1}^{m_{N}}\mathbb{P}_{n}^{(N-1)}(\left\{
X_{i}\right\}  )1_{A_{j}^{(N)}}(X_{i})\\
& =-\sum_{j=1}^{m_{N}}\mathbb{P}_{n}^{(N-1)}(A_{j}^{(N)})\sum_{i=1}^{n}%
\frac{\mathbb{P}_{n}^{(N-1)}(\left\{  X_{i}\right\}  )}{\mathbb{P}_{n}%
^{(N-1)}(A_{j}^{(N)})}1_{A_{j}^{(N)}}(X_{i})\log\left(  \frac{Q(\left\{
X_{i}\right\}  )}{\mathbb{P}_{n}^{(N-1)}(\left\{  X_{i}\right\}  )}\right) \\
& \geqslant-\sum_{j=1}^{m_{N}}\mathbb{P}_{n}^{(N-1)}(A_{j}^{(N)})\log\left(
\sum_{i=1}^{n}\frac{Q(\left\{  X_{i}\right\}  )}{\mathbb{P}_{n}^{(N-1)}%
(A_{j}^{(N)})}1_{A_{j}^{(N)}}(X_{i})\right) \\
& =\sum_{j=1}^{m_{N}}\mathbb{P}_{n}^{(N-1)}(A_{j}^{(N)})\log\left(
\frac{\mathbb{P}_{n}^{(N-1)}(A_{j}^{(N)})}{P(A_{j}^{(N)})}\right)
=d_{K}\left(  \mathbb{P}_{n}^{(N-1)}\mid\mid\mathbb{P}_{n}^{(N)}\right),
\end{align*}
where the final identification relies on (\ref{DKA}) and $\mathbb{P}_{n}^{(N)}=P$ on $\mathcal{A}^{(N)}$.

\subsection{Proof of Proposition~\ref{ProLLIMes}.}

The classical law of the iterated logarithm holds for the empirical process
$\alpha_{n}^{(0)}$ indexed by $\mathcal{F}$ under (VC) and (BR). See Alexander
\cite{Alex84}, in particular Theorem 2.12 for (BR) and Theorem~2.13
based on Theorem 2.8 that uses in its proof the consequence of Lemma 2.7, which is indeed (VC). Namely, for any $\varepsilon>0$,
with probability one there exists $n(\omega)$ such that, for all $n>n(\omega
)$,%
\begin{equation}
u_{n}||\mathbb{P}_{n}^{(0)}-P||_{\mathcal{F}}<1+\varepsilon,\quad u_{n}%
=\sqrt{\frac{n}{2\sigma_{\mathcal{F}}^{2}L\circ L(n)}}\label{loglog},%
\end{equation}
where $\sigma_{\mathcal{F}}^{2}=\sup_{\mathcal{F}}\mathbb{V}(f)\leqslant
M^{2}$. Let $1\leqslant N\leqslant N_{0}$. By (\ref{PnN}) it holds%
\begin{align*}
\mathbb{P}_{n}^{(N)}(f)-\mathbb{P}_{n}^{(N-1)}(f)  & =\sum_{j=1}^{m_{N}%
}\mathbb{P}_{n}^{(N-1)}(f1_{A_{j}^{(N)}})\frac{P(A_{j}^{(N)})}{\mathbb{P}%
_{n}^{(N-1)}(A_{j}^{(N)})}-\sum_{j=1}^{m_{N}}\mathbb{P}_{n}^{(N-1)}%
(f1_{A_{j}^{(N)}})\\
& =\sum_{j=1}^{m_{N}}\mathbb{E}_{n}^{(N-1)}(f\;|\;A_{j}^{(N)})\left(
P(A_{j}^{(N)})-\mathbb{P}_{n}^{(N-1)}(A_{j}^{(N)})\right)  .
\end{align*}
Since $\mathbb{P}_{n}^{(N)}$ is a probability measure we have $\Vert
\mathbb{P}_{n}^{(N)}(f1_{A})\Vert _{\mathcal{F}}\leqslant M\mathbb{P}%
_{n}^{(N)}(A)$ hence $\Vert \mathbb{E}_{n}^{(N-1)}(f\;|\;A_{j}%
^{(N)})\Vert _{\mathcal{F}}\leqslant M$ and $\vert
\mathbb{P}_{n}^{(N)}(f)-\mathbb{P}_{n}^{(N-1)}(f)\vert \leqslant
Mm_{N}\Vert \mathbb{P}_{n}^{(N-1)}-P\Vert _{\mathcal{F}}$. Also
observe that (\ref{deltaN}) combined with the fact that $\mathcal{A}^{(N)}$ is a
partition implies $m_{N}\leqslant1/p_{N}$ and $p_{N}\geqslant
p_{(N)}\geqslant p_{(N_{0})}$. Therefore%
\begin{align*}
u_{n}\left\Vert \mathbb{P}_{n}^{(N)}-P\right\Vert _{\mathcal{F}}  & \leqslant
u_{n}\left\Vert \mathbb{P}_{n}^{(N-1)}-P\right\Vert _{\mathcal{F}}%
+u_{n}\left\Vert \mathbb{P}_{n}^{(N)}-\mathbb{P}_{n}^{(N-1)}\right\Vert
_{\mathcal{F}}\\
& \leqslant u_{n}\left(  1+Mm_{N}\right)  \left\Vert \mathbb{P}_{n}%
^{(N-1)}-P\right\Vert _{\mathcal{F}}\\
& \leqslant u_{n}\kappa_{N}\left\Vert \mathbb{P}_{n}^{(0)}-P\right\Vert
_{\mathcal{F}},%
\end{align*}
where, for $N>0$,
\[
\kappa_{N}=\prod_{k=1}^{N}\left(  1+Mm_{N}\right)  \leqslant\prod_{k=1}%
^{N}\left(  1+\frac{M}{p_{k}}\right)  \leqslant\left(  1+\frac{M}%
{p_{(N_{0})}}\right)  ^{N_{0}},%
\]
which by (\ref{loglog}) remains true for $N=N_{0}=0$. This proves that, given
$N_{0}\in\mathbb{N}$ and for all $\varepsilon>0$,%
\[
\underset{n\rightarrow+\infty}{\lim\sup}\sqrt{\frac{n}{2L\circ L(n)}}%
\sup_{0\leqslant N\leqslant N_{0}}\left\Vert \mathbb{P}_{n}^{(N)}-P\right\Vert
_{\mathcal{F}}\leqslant(1+\varepsilon)\kappa_{N}\sigma_{\mathcal{F}}\quad a.s.
\]
and Proposition \ref{ProLLIMes} follows.

\subsection{Proof of Proposition \ref{ProConcentr}}

\noindent\textbf{Step 1.} We work on $B_{n,N_{0}}$ from (\ref{biendef}), which means that all the probabilities considered below concern events that we implicitly assume to be intersected with $B_{n,N_{0}}$. By (\ref{condalphanN}) and
(\ref{alphanN}) we have, for $N\geqslant1$,
\[
\alpha_{n,j}^{(N-1)}(f)=\frac{1}{\mathbb{P}_{n}^{(N-1)}(A_{j}^{(N)})}%
\alpha_{n}^{(N-1)}\left(  (f-\mathbb{E}(f|A_{j}^{(N)}))1_{A_{j}^{(N)}}\right),
\]
with $\left\vert \mathbb{E}(f|A_{j}^{(N)})\right\vert \leqslant M$, and%
\begin{align*}
\mathbb{P}\left(  \left\Vert \alpha_{n}^{(N)}\right\Vert _{\mathcal{F}%
}\geqslant\lambda\right)   &  \leqslant\mathbb{P}\left(  \sum_{j=1}^{m_{N}%
}P(A_{j}^{(N)})\left\Vert \alpha_{n,j}^{(N-1)}\right\Vert _{\mathcal{F}%
}\geqslant\sum_{j=1}^{m_{N}}P(A_{j}^{(N)})\lambda\right)  \\
&  \leqslant\sum_{j=1}^{m_{N}}\mathbb{P}\left(  \left\Vert \alpha
_{n,j}^{(N-1)}\right\Vert _{\mathcal{F}}\geqslant\lambda\right)  .
\end{align*}
Each term in the latter sum satisfies, for any positive numbers $K\leqslant P(A_{j}^{(N)})$ and
$K^{\prime}\leqslant P(A_{j}^{(N)})-K$,%
\begin{align}
&  \mathbb{P}\left(  \frac{1}{\mathbb{P}_{n}^{(N-1)}(A_{j}^{(N)})}\left\Vert
\alpha_{n}^{(N-1)}(f)-\mathbb{E}(f|A_{j}^{(N)})\alpha_{n}^{(N-1)}(A_{j}%
^{(N)})\right\Vert _{\mathcal{F}}\geqslant\lambda\right)  \nonumber\\
&  \leqslant\mathbb{P}\left(  (1+M)\left\Vert \alpha_{n}^{(N-1)}\right\Vert
_{\mathcal{F}}\geqslant K\lambda\right)  +\mathbb{P}\left(  \mathbb{P}%
_{n}^{(N-1)}(A_{j}^{(N)})\leqslant K\right)  \nonumber\\
&  \leqslant\mathbb{P}\left(  \left\Vert \alpha_{n}^{(N-1)}\right\Vert
_{\mathcal{F}}\geqslant\frac{\lambda K}{1+M}\right)  +\mathbb{P}\left(
\alpha_{n}^{(N-1)}(A_{j}^{(N)})\leqslant-K^{\prime}\sqrt{n}\right)
\nonumber\\
&  \leqslant2\mathbb{P}\left(  \left\Vert \alpha_{n}^{(N-1)}\right\Vert
_{\mathcal{F}}\geqslant\frac{\lambda K}{1+M}\right),  \label{lambdaK}%
\end{align}
where the last bound holds provided that $K^{\prime}\sqrt{n}\geqslant\lambda
K/(1+M)$. Define%
\[
\beta=\frac{1}{1+\lambda/(1+M)\sqrt{n}}\in\left(  0,1\right), 
\quad K=\beta p_{N}, \quad K^{\prime}=p_{N}(1-\beta),
\]
where $p_{N}$ is defined in (\ref{deltaN}). Since
$p_{N}\leqslant1/2$ for any $N\geqslant1$ it holds $K^{\prime}>0$ and $K^{\prime}\sqrt{n}=\lambda
K/(1+M)$.
We have shown that for any $N\geqslant 1$,%
\[
\mathbb{P}\left(  \left\Vert \alpha_{n}^{(N)}\right\Vert _{\mathcal{F}%
}\geqslant\lambda\right)  \leqslant2m_{N}\mathbb{P}\left(  \left\Vert
\alpha_{n}^{(N-1)}\right\Vert _{\mathcal{F}}\geqslant\frac{\lambda\beta
p_{N}}{1+M}\right)  .
\]
Applying (\ref{lambdaK}) again with $\lambda$ turned into the smaller $\lambda\beta p_{N}/(1+M)$ then iterating backward from $N_{0}$ we get, for $P_{N_{0}}=\prod_{N=1}^{N_{0}}p_{N}$ and $M_{N_{0}}%
=\prod_{N=1}^{N_{0}}m_{N}\leqslant1/P_{N_{0}}$,%

\begin{equation}
\mathbb{P}\left(  \left\Vert \alpha_{n}^{(N_{0})}\right\Vert _{\mathcal{F}%
}\geqslant\lambda\right)  \leqslant2^{N_{0}}M_{N_{0}}\mathbb{P}\left(
\left\Vert \alpha_{n}^{(0)}\right\Vert _{\mathcal{F}}\geqslant\frac
{\lambda P_{N_{0}}}{(1+M+\lambda/\sqrt{n})^{N_{0}}}\right).
\label{step1}
\end{equation}
The latter upper bound being increasing with $N_{0}$ we conclude that%
\[
\mathbb{P}\left(  \max_{0\leqslant N\leqslant N_{0}}\left\Vert \alpha
_{n}^{(N)}\right\Vert _{\mathcal{F}}\geqslant\lambda\right)  \leqslant
\sum_{N=1}^{N_{0}}\mathbb{P}\left(  \left\Vert \alpha_{n}^{(N)}\right\Vert
_{\mathcal{F}}\geqslant\lambda\right)  \leqslant N_{0}\mathbb{P}\left(
\left\Vert \alpha_{n}^{(N_{0})}\right\Vert _{\mathcal{F}}\geqslant
\lambda\right)  .
\]
\noindent\textbf{Step 2.} By Theorem 2.14.25 of van der Vaart
and Wellner~\cite{VanWell} or Corollary 2 of~\cite{BirgMass98}, for $n\geqslant1$, $t>0$, we have for some universal constants $D_{1}^{\prime}>0$, $D_{2}^{\prime}>0$, %
\[
\mathbb{P}\left(  ||\alpha_{n}^{(0)}||_{\mathcal{F}}>D_{1}^{\prime}\left(
\mathbb{\mu}_{n}+t\right)  \right)  \leqslant\exp\left(  -D_{2}^{\prime}%
\min\left(  \frac{t^{2}}{\sigma_{\mathcal{F}}^{2}},\frac{t\sqrt{n}}{M}\right)
\right),
\]
where, by the last maximal inequality in Theorem 2.14.2 of~\cite{VanWell} applied to $\mathcal{F}$ with envelop function constant to
$M$, it holds%
\[
\mathbb{\mu}_{n}=\mathbb{E}\left\Vert \frac{1}{\sqrt{n}}\sum_{i=1}^{n}%
f(X_{i})\right\Vert _{\mathcal{F}}\leqslant M\int_{0}^{1}%
\sqrt{1+\log N_{\left[ \ \right]  }(\mathcal{F},M\varepsilon,d_{P})}d\varepsilon.
\]
Under (BR), we have $\mathbb{\mu}_{n}<C$ with $C=M(1+b_0 / (1-r_0))$. 
For $\lambda_{0}=2D'_1C$ we get, for any $n>0$ and $\lambda_{0}<\lambda<2D_{1}^{\prime}\sigma_{\mathcal{F}}^{2}\sqrt{n}/M$,
\[
\mathbb{P}\left(  ||\alpha_{n}^{(0)}||_{\mathcal{F}}>\lambda\right)
\leqslant\mathbb{P}\left(  ||\alpha_{n}^{(0)}||_{\mathcal{F}}>D_{1}^{\prime
}\left(  \mathbb{\mu}_{n}+\frac{\lambda}{2D_{1}^{\prime}}\right)  \right)
\leqslant\exp\left(  -D_{2}^{\prime\prime}\lambda^{2}\right),
\]
where $D_{2}^{\prime\prime}=D_{2}^{\prime}/4(D_{1}^{\prime})^{2}%
\sigma_{\mathcal{F}}^{2}$. Therefore, according to (\ref{step1}), taking%
\begin{equation}
D_{0}=\frac{2D_{1}^{\prime}\sigma_{\mathcal{F}}^{2}}{M},\quad D_{1}%
=N_{0}2^{N_{0}}M_{N_{0}},\quad D_{2}=\frac{D_{2}^{\prime\prime}P_{N_{0}%
}^{2}}{(1+M+D_{0})^{2N_{0}}},\label{D012}%
\end{equation}
yields $\mathbb{P}(\max_{0\leqslant N\leqslant N_{0}}
||\alpha_{n}^{(N_{0})}||_{\mathcal{F}}\geqslant\lambda)\leqslant D_{1}\exp\left(  -D_{2}\lambda^{2}\right)  $ 
for $\lambda_{0}<\lambda<D_{0}\sqrt{n}$.\medskip

\noindent\textbf{Step 3.}\ By Theorem 2.14.9 in \cite{VanWell}, under (VC) there exists a constant
$D(c_{0})$ such that, for $t_{0}$ large enough and all $t\geqslant t_{0}$,%
\[
\mathbb{P}\left(  \left\Vert \alpha_{n}^{(0)}\right\Vert _{\mathcal{F}%
}\geqslant t\right)  \leqslant\left(  \frac{D(c_{0})t}{M\sqrt{v_{0}}}\right)
^{v_{0}}\exp\left(  -\frac{2t^{2}}{M^{2}}\right)  .
\]
Denote $\lambda_{1,n}$ and $\lambda_{2,n}$ the two solutions of $\lambda
P_{N_{0}}=t_{0}(1+M+\lambda/\sqrt{n})^{N_{0}}$. Notice that, for $n$
large, $\lambda_{1,n}$ is close to $t_{0}(1+M)^{N_{0}}/P_{N_{0}}$ and
$\lambda_{2,n}=O(n^{N_{0}/2(N_{0}-1)})$. Combined with (\ref{step1}) it ensues that for some $n_{0},\lambda_{0}$
it holds, for all $n>n_{0}$ and $\lambda_{0}<\lambda<2M\sqrt{n}$, $\mathbb{P}(\max_{0\leqslant N\leqslant N_{0}}\vert\vert \alpha
_{n}^{(N)}\vert\vert_{\mathcal{F}}\geqslant\lambda)
\leqslant D_{3}\lambda^{v_{0}}\exp(-D_{4}\lambda^{2})$
where
\begin{equation}
D_{3}=\frac{N_{0}2^{N_{0}}M_{N_{0}}}{(1+M)^{v_{0}N_{0}}}\left(  \frac{D(c_{0})P_{N_{0}}}%
{M\sqrt{v_{0}}}\right)  ^{v_{0}},\quad D_{4}=\frac{2P_{N_{0}}^{2}}%
{M^{2}(3M+1)^{2N_0}}.\label{D34}%
\end{equation}
Finally, at each step, add $S_{N_{0}}\left(1-p_{(N_{0})}\right)^{n}$ to take $B_{n,N_{0}}^{c}$ from (\ref{biendef}) into account.

\subsection{Proof of Proposition~\ref{ProCVProc} and Theorem
\ref{ProCVSpeed}.}

\noindent Theorem\ \ref{ProCVSpeed} implies Proposition \ref{ProCVProc}\ 
since the weak convergence on $(\ell^{\infty}(\mathcal{F}),\left\Vert
\cdot\right\Vert _{\mathcal{F}})$ is metrized by the L\'{e}vy-Prokhorov
distance between $\alpha_{n}^{(N)}$ and $\mathbb{G}_{n}^{(N)}$ which is%
\begin{equation}
\inf\left\{  \varepsilon
>0:\mathbb{P}^{\alpha_{n}^{(N)}}(A)\leqslant\mathbb{P}^{\mathbb{G}_{n}^{(N)}%
}(A^{\varepsilon})+\varepsilon,\mathbb{P}^{\mathbb{G}_{n}^{(N)}}%
(A)\leqslant\mathbb{P}^{\alpha_{n}^{(N)}}(A^{\varepsilon})+\varepsilon
\right\}  \leqslant d_0 v_{n}.\label{Levy}
\end{equation}
To see this, recall that we have $v_{n}>1/n^{\theta_{0}}$ for $\theta_{0}>1/2$ and $v_{n}\rightarrow0$ in
Theorem\ \ref{ProCVSpeed}, remind (\ref{PnA}) and (\ref{biendef}) then observe that%
\begin{align*}
\mathbb{P}\left(  \alpha_{n}^{(N)}\in A\right)   &  \leqslant\mathbb{P}\left(
\left\{  \alpha_{n}^{(N)}\in A\right\}  \cap\left\{  \left\Vert \alpha
_{n}^{(N)}-\mathbb{G}_{n}^{(N)}\right\Vert _{\mathcal{F}}\leqslant
d_0 v_{n}\right\}  \cap B_{n,N_{0}}\right)  \\
&  +\mathbb{P}\left(  \left\{  \left\Vert \alpha_{n}^{(N)}-\mathbb{G}%
_{n}^{(N)}\right\Vert _{\mathcal{F}}>d_0 v_{n}\right\}  \cap B_{n,N_{0}}\right)
+\mathbb{P}\left(  B_{n,N_{0}}^{c}\right)  \\
&  \leqslant\mathbb{P}(\mathbb{G}_{n}^{(N)}\in A^{v_{n}})+\frac{1}{n^{\theta_{0}}%
}+S_{N_{0}}\left(  1-p_{(N_{0})}\right)  ^{n},%
\end{align*}
which obviously remains true by exchanging $\alpha_{n}^{(N)}$ and $\mathbb{G}_{n}^{(N)}$. 
Since $v_{n}$ is the slowest sequence as $n\rightarrow +\infty$, if $n_{0}$ satisfies
$v_{n_{0}}>1/n_{0}^{\theta_{0}}+S_{N_{0}}\left(1-p_{(N_{0})}\right)^{n_{0}}$
then $v_{n}>1/n^{\theta_{0}}+S_{N_{0}}\left(1-p_{(N_{0})}\right)^{n}$
for all $n>n_{0}$. Whence (\ref{Levy}).\textit{\medskip}

\noindent We next establish Theorem\ \ref{ProCVSpeed}. 
Fix $N_{0}\in\mathbb{N}$.\textit{\medskip}

\noindent\textbf{Step 1.} Let introduce the transforms, for $f\in\mathcal{F}$,
$N\geqslant1$ and $1\leqslant j\leqslant m_{N}$,%
\begin{align*}
& \phi_{(j,N)}f=\left(  f-\mathbb{E}\left(  f\;|\;A_{j}^{(N)}\right)  \right)
1_{A_{j}^{(N)}},\\ 
& \phi_{(N)}f=\sum_{j=1}^{m_{N}}\phi_{(j,N)}f=f-\sum
_{j=1}^{m_{N}}\mathbb{E}(f\;|\;A_{j}^{(N)})1_{A_{j}^{(N)}}.
\end{align*}
It holds $P(\phi_{(N)}f)=P(\phi_{(j,N)}f)=0$ and, since $\mathcal{A}^{(N)}$ is
a partition of $\mathcal{X}$,
\begin{equation}
(\phi_{(j,N)}f)(\phi_{(j^{\prime},N)}g)=0,\quad1\leqslant j\neq j^{\prime
}\leqslant m_{N}.\label{phiortho}%
\end{equation}
Moreover, the $L_{2}(P)$ property of conditional expectations yields, with the notation (\ref{sigma}),
\begin{equation}
\sigma_{\phi_{(j,N)}f}^{2}=P(f_{(j,N)}^{2})\leqslant\sigma_{\phi_{(N)}f}%
^{2}=P(f_{(N)}^{2})=\sum_{j=1}^{m_{N}}\sigma_{\phi_{(j,N)}f}^{2}%
\leqslant\sigma_{f}^{2}.\label{var}%
\end{equation}
Next consider the class of backward iterated transforms%
\begin{align*}
& \mathcal{F}_{(N)}=\phi_{(1)}\circ...\circ\phi_{(N)}(\mathcal{F)}%
,\\
& \mathcal{H}_{(N)}=\bigcup\nolimits_{1\leqslant k\leqslant N}%
\bigcup\nolimits_{1\leqslant j\leqslant m_{k}}\phi_{(j,k)}\circ\phi
_{(k+1)}\circ...\circ\phi_{(N)}(\mathcal{F}),
\end{align*}
where $\phi_{(k+1)}\circ...\circ\phi_{(N)}=id$ if $k=N\geqslant1$ and
$\mathcal{F}_{(0)}=\mathcal{H}_{(0)}=\mathcal{F}$. Also write $\mathcal{F}%
_{0}=\bigcup\nolimits_{0\leqslant N\leqslant N_{0}}\mathcal{F}_{(N)}$ and
$\mathcal{H}_{0}=\bigcup\nolimits_{0\leqslant N\leqslant N_{0}}\mathcal{H}%
_{(N)}$. By iterating (\ref{var}) it comes $\sigma_{\mathcal{H}_{0}}%
^{2}\leqslant\sigma_{\mathcal{F}_{0}}^{2}\leqslant\sigma_{\mathcal{F}}^{2}$.
We first show that properties of $\mathcal{F}$ transfer to $\mathcal{F}_{(N)}%
$, $\mathcal{H}_{(N)}$ for $0\leqslant N\leqslant N_{0}$ and thus to
$\mathcal{F}_{0}$ and $\mathcal{H}_{0}$. Remind the constants defined at (\ref{deltaMN}).

\begin{lem}
\label{LemEntrop}Assume (\ref{deltaN}). If $\mathcal{F}$ is pointwise
measurable and bounded by $M$ then $\mathcal{F}_{(N)}$ and $\mathcal{H}_{(N)}$
(resp. $\mathcal{F}_{0}$ and $\mathcal{H}_{0}$) are pointwise measurable and
bounded by $(2M)^{N}/P_{N}$ (resp. $(2M)^{N_{0}}/P_{N_{0}}$). If 
$(VC)$ (resp. $(BR)$) holds then $\mathcal{F}_{0}$ and
$\mathcal{H}_{0}$ also satisfies $(VC)$ (resp. $(BR)$) with the same power
$\nu_{0}$ (resp. $r_{0}$) as $\mathcal{F}$.
\end{lem}

\noindent\textit{Proof.} If $\mathcal{F}$ is uniformly bounded by $M$ then for
$N\leqslant N_{0}$ we have%
\[
\sup_{\mathcal{F}}\sup_{\mathcal{X}}\left\vert \phi_{(N)}f\right\vert
=\sup_{\mathcal{F}}\max_{1\leqslant j\leqslant m_{N}}\sup_{\mathcal{X}%
}\left\vert \phi_{(j,N)}f\right\vert \leqslant M\left(  1+\frac{1}{p_{N}%
}\right)  \leqslant\frac{2M}{p_{N}},%
\]
thus, by backward induction from $N_{0}$ to $1$, $\mathcal{F}_{(N_{0})}$ and
$\mathcal{H}_{(N_{0})}$ are uniformly bounded by $(2M)^{N}/P_{N}$. It
readily follows that $\mathcal{F}_{0}$ and $\mathcal{H}_{0}$ are bounded by
$(2M)^{N_{0}}/P_{N_{0}}$. Assume that $f_{k}\in\mathcal{F}_{\ast}$
converges pointwise on $\mathcal{X}$ to $f\in\mathcal{F}$. From
\[
\lim_{k\rightarrow+\infty}1_{A_{j}^{(N)}}(X)f_{k}(X)=1_{A_{j}^{(N)}%
}(X)f(X)\quad\text{and}\quad P(1_{A_{j}^{(N)}}|f_{k}|)\leqslant P\left(
|f_{k}|\right)  \leqslant M,
\]
we deduce by dominated convergence that $\lim_{k\rightarrow+\infty
}\mathbb{E}(f_{k}\;|\;A_{j}^{(N)})=\mathbb{E}(f\;|\;A_{j}^{(N)})$. Thus
$\phi_{(j,N)}f_{k}$ converges pointwise to $\phi_{(j,N)}f$ and $\phi
_{(N)}f_{k}=\sum_{j=1}^{m_{N}}\phi_{(j,N)}f_{k}$ to $\phi_{(N)}%
f=\sum_{j=1}^{m_{N}}\phi_{(j,N)}f$. By iterating this reasoning backward
from $N$ to $1$ we obtain that $\mathcal{F}_{(N)}$ and $\mathcal{H}_{(N)}$ are
pointwise measurable, by using the countable classes $\phi_{(1)}\circ
...\circ\phi_{(N)}(\mathcal{F}_{\ast})$ and $%
{\textstyle\bigcup\nolimits_{1\leqslant k\leqslant N}}
{\textstyle\bigcup\nolimits_{1\leqslant j\leqslant m_{k}}}
\phi_{(j,k)}\circ\phi_{(k+1)}\circ...\circ\phi_{(N)}(\mathcal{F}_{\ast})$
respectively. Assume next that $\mathcal{F}$ satisfies $(VC)$. By
(\ref{phiortho}) we have%
\begin{align*}
d_{Q}^{2}(\phi_{(N)}f,\phi_{(N)}g)  & =\int_{\mathcal{X}}\left(
{\textstyle\sum\nolimits_{j=1}^{m_{N}}}
(\phi_{(j,N)}f-\phi_{(j,N)}g)\right)  ^{2}dQ\\
& =\sum_{j=1}^{m_{N}}d_{Q}^{2}(\phi_{(j,N)}f,\phi_{(j,N)}g)\\
& =\sum_{j=1}^{m_{N}}\int_{A_{j}^{(N)}}\left(  f-g-\mathbb{E(}f-g\;|\;A_{j}%
^{(N)})\right)  ^{2}dQ\\
& \leqslant\sum_{j=1}^{m_{N}}\int_{A_{j}^{(N)}}\left(  f-g-(Qf-Qg)\right)
^{2}dQ\\
& =d_{Q}^{2}(f,g)-(Qf-Qg)^{2},%
\end{align*}
thus $d_{Q}(f,g)<\varepsilon$ implies $d_{Q}^{2}(\phi_{(N)}f,\phi
_{(N)}g)\leqslant d_{Q}^{2}(f,g)<\varepsilon^{2}$. If $\mathcal{F}$ can be
covered by $N\left(  \mathcal{F},\varepsilon,d_{Q}\right)  $ balls of $d_{Q}%
$-radius $\varepsilon$ centered at some $g$ then $\phi_{(N)}(\mathcal{F})$ can
be covered by the same number of balls, centered at the corresponding
$\phi_{(N)}g$ and hence the same number of centers $\phi_{(1)}\circ
...\circ\phi_{(N)}g$ suffices to cover $\mathcal{F}_{(N)}$. All the
$\phi_{(j,k)}\circ\phi_{(k+1)}\circ...\circ\phi_{(N)}g$ are needed to cover
$\mathcal{H}_{(N)}$, that is $S_{N}N\left(  \mathcal{F},\varepsilon
,d_{Q}\right)  $. This shows that $\mathcal{F}_{0}$ (resp. $\mathcal{H}_{0}$)
obeys $(VC)$ with the same power $\nu_{0}$ and a constant $c_{0}(N_{0}+1)$
(resp. $c_{0}\sum_{N=0}^{N_{0}}S_{N}$). Assume now that $\mathcal{F}$
satisfies $(BR)$. If $g_{-}\leqslant f\leqslant g_{+}$ then we have
\begin{align*}
h_{(j,N)}^{-} &  =1_{A_{j}^{(N)}}g_{-}-1_{A_{j}^{(N)}}\mathbb{E}\left(
g_{+}\;|\;A_{j}^{(N)}\right)  \\
&  \leqslant\phi_{(j,N)}f\leqslant1_{A_{j}^{(N)}}g_{+}-1_{A_{j}^{(N)}%
}\mathbb{E}\left(  g_{-}\;|\;A_{j}^{(N)}\right)  =h_{(j,N)}^{+},%
\end{align*}
and the $L_{2}(P)$-size of the new bracket $[h_{(j,N)}^{-},h_{(j,N)}^{+}]$ is%
\begin{align*}
d_{P}^{2}(h_{(j,N)}^{-},h_{(j,N)}^{+}) &  =\int_{A_{j}^{(N)}}\left(
g_{+}-g_{-}+\mathbb{E}(g_{+}-g_{-}\;|\;A_{j}^{(N)})\right)  ^{2}dP\\
&  =P(1_{A_{j}^{(N)}}(g_{+}-g_{-})^{2})+P(A_{j}^{(N)})\mathbb{E}(g_{+}%
-g_{-}\;|\;A_{j}^{(N)})^{2}\\
&  +2\mathbb{E}(g_{+}-g_{-}\;|\;A_{j}^{(N)})P(1_{A_{j}^{(N)}}(g_{+}-g_{-})).
\end{align*}
If $d_{P}(g_{+},g_{-})<\varepsilon$ the H\"{o}lder
inequality yields $P(1_{A_{j}^{(N)}}(g_{+}-g_{-}))\leqslant\varepsilon\sqrt
{P(A_{j}^{(N)})}$ and $\mathbb{E}(g_{+}-g_{-}\;|\;A_{j}^{(N)})\leqslant
\varepsilon/\sqrt{P(A_{j}^{(N)})}$ hence%
\[
d_{P}^{2}(h_{(j,N)}^{-},h_{(j,N)}^{+})\leqslant P(1_{A_{j}^{(N)}}(g_{+}%
-g_{-})^{2})+3\varepsilon^{2},%
\]
so that $\phi_{(N)}f=\sum_{j=1}^{m_{N}}\phi_{(j,N)}f\in\lbrack h_{(N)}%
^{-},h_{(N)}^{+}]$ where $h_{(N)}^{\pm}=\sum_{j=1}^{m_{N}}h_{(j,N)}^{\pm}$
satisfies%
\[
d_{P}^{2}(h_{(N)}^{-},h_{(N)}^{+})=\sum_{j=1}^{m_{N}}d_{P}^{2}(h_{(j,N)}%
^{-},h_{(j,N)}^{+})\leqslant d_{P}^{2}(g_{-},g_{+})+3m_{N}\varepsilon
^{2}\leqslant4m_{N}\varepsilon^{2}.
\]
It ensues $N_{[\ ]}(\phi_{(N)}(\mathcal{F}),\varepsilon,d_{P})\leqslant
N_{[\ ]}(\mathcal{F},\varepsilon/2\sqrt{m_{N}},d_{P})$ and $N_{[\ ]}%
(\mathcal{F}_{(N)},\varepsilon,d_{P})\leqslant N_{[\ ]}(\mathcal{F}%
,\varepsilon/2^{N}\sqrt{M_{N}},d_{P})$. To cover $\phi_{(j,k)}\circ
\phi_{(k+1)}\circ...\circ\phi_{(N)}(\mathcal{F})$ one needs at most $m_{k}%
N_{[\ ]}(\mathcal{F},\varepsilon/2^{N-k}\sqrt{m_{k+1}...m_{N}},d_{P})$
brackets. We have proved that%
\begin{align*}
N_{[\ ]}(\mathcal{F}_{0},\varepsilon,d_{P}) &  \leqslant(N_{0}+1)N_{[\ ]}%
(\mathcal{F},\varepsilon/2^{N_{0}}\sqrt{M_{N_{0}}},d_{P}),\\
N_{[\ ]}(\mathcal{H}_{0},\varepsilon,d_{P}) &  \leqslant S_{N_{0}}%
N_{[\ ]}(\mathcal{F},\varepsilon/2^{N_{0}}\sqrt{M_{N_{0}}},d_{P}).
\end{align*}
Therefore $\mathcal{F}_{0}$, $\mathcal{H}_{0}$ satisfy $(BR)$ with power
$r_{0}$ and constant $2^{r_{0}N_{0}}M_{N_{0}}^{r_{0}}b_{0}$.$\quad\square
$\textit{\medskip}

\noindent\textbf{Step 2.} By (\ref{alphanN}) we have
\begin{align}
&  \alpha_{n}^{(N)}(f)=\sum_{j=1}^{m_{N}}\frac{P(A_{j}^{(N)})}{\mathbb{P}%
_{n}^{(N-1)}(A_{j}^{(N)})}\alpha_{n}^{(N-1)}(\phi_{(j,N)}f)=\alpha_{n}%
^{(N-1)}(\phi_{(N)}f)+\Gamma_{n}^{(N)}(f),\label{alphagamma}\\
&  \Gamma_{n}^{(N)}(f)=\sum_{j=1}^{m_{N}}q_{n}(j,N)\alpha_{n}^{(N-1)}\left(
\phi_{(j,N)}f\right)  ,\quad q_{n}(j,N)=\frac{P(A_{j}^{(N)})}{\mathbb{P}%
_{n}^{(N-1)}(A_{j}^{(N)})}-1.\nonumber
\end{align}
Under the convention that $\phi_{(N+1)}
\circ\phi_{(N)}=id$, iterating (\ref{alphagamma}) leads to%
\begin{align*}
& \alpha_{n}^{(N)}(f)=\alpha_{n}^{(0)}(\phi_{(1)}\circ...\circ\phi
_{(N)}f)+\digamma_{n}^{(N)}(f),\\ 
& \digamma_{n}^{(N)}(f)=\sum_{k=1}^{N}%
\Gamma_{n}^{(k)}(\phi_{(k+1)}\circ...\circ\phi_{(N)}f).
\end{align*}
Clearly the terms $\Gamma_{n}^{(k)}$ carry out some bias
and variance distortion. However the following lemma states that
$\alpha_{n}^{(0)}(\mathcal{F}_{(N)})$ is the main contribution to $\alpha
_{n}^{(N)}(\mathcal{F})$ and $\digamma_{n}^{(N)}(\mathcal{F})$ is an error
process.

\begin{lem}
\label{LemGamma}Consider the sequence
$v_{n}$ defined at Theorem\ \ref{ProCVSpeed}. 
If $\mathcal{F}$ satisfies $(VC)$ or $(BR)$ then there exists
$C_{0}<+\infty$ such that we almost surely have, for all $n$ large enough,
$\max_{0\leqslant N\leqslant N_{0}}\Vert\digamma_{n}^{(N)}\Vert_{\mathcal{F}%
}\leqslant C_{0}L\circ L(n)/\sqrt{n}$. Moreover, for any $\zeta>0$ and
$\theta>0$ there exists $n_3(\zeta,\theta)$ such that we have, for all
$n>n_3(\zeta,\theta)$,
\[
\mathbb{P}\left(  \max_{0\leqslant N\leqslant N_{0}}\Vert\digamma_{n}%
^{(N)}\Vert_{\mathcal{F}}>\zeta v_{n}\right)  \leqslant\frac{1}{2n^{\theta}}.
\]

\end{lem}

\noindent\textit{Proof.} \textbf{(i)} Let us apply Proposition~\ref{ProLLIMes}
to $\mathcal{F}$ and, thanks to Lemma \ref{LemEntrop}, to $\mathcal{H}_{0}$
and $\mathcal{H}_{(N)}$. So for all $ \varepsilon> 0 $ we have for all $ n $ large enough, $$ \max_{1\leqslant N\leqslant N_{0}}\max_{1\leqslant j\leqslant m_{N}}\left\vert
\alpha_{n}^{(N-1)}(A_{j}^{(N)})\right\vert \leqslant  \sigma_\mathcal{F}\kappa_{N_0}\sqrt{2L\circ L(n)} (1+\varepsilon). $$ The following statements are almost surely true, for
all $n$ large enough. On the one hand, for $ \varepsilon = \sqrt{2}-1>0$, %
\begin{equation}
\max_{1\leqslant N\leqslant N_{0}}\max_{1\leqslant j\leqslant m_{N}}\left\vert
\alpha_{n}^{(N-1)}(A_{j}^{(N)})\right\vert \leqslant b_{n}=2\sigma
_{\mathcal{F}}\kappa_{N_{0}}\sqrt{L\circ L(n)}.\label{supalpha}%
\end{equation}
On the other hand, having $\sigma_{\mathcal{H}_{0}}\leqslant\sigma
_{\mathcal{F}}$ by (\ref{var}),%
\[
\max_{1\leqslant N\leqslant N_{0}}\max_{1\leqslant k\leqslant N}%
\max_{1\leqslant j\leqslant m_{k}}\left\vert \alpha_{n}^{(k-1)}(\phi
_{(j,k)}\circ\phi_{(k+1)}\circ...\circ\phi_{(N)}f)\right\vert \leqslant
\max_{1\leqslant k\leqslant N_{0}}\left\Vert \alpha_{n}^{(k-1)}\right\Vert
_{\mathcal{H}_{0}}\leqslant b_{n}.
\]
By (\ref{supalpha}), $q_{n}(j,N)=1/\left(1+\alpha_{n}^{(N-1)}(A_{j}^{(N)}%
)/P(A_{j}^{(N)})\sqrt{n}\right)-1$ satisfies
\begin{equation}
\max_{1\leqslant N\leqslant N_{0}}\max_{1\leqslant j\leqslant m_{N}}\left\vert
q_{n}(j,N)\right\vert \frac{\sqrt{n}}{b_{n}}p_{(N_{0})}\leqslant
2\label{qn},%
\end{equation}
which implies%
\begin{align}
\left\Vert \digamma_{n}^{(N)}\right\Vert _{\mathcal{F}} &  \leqslant\sum
_{k=1}^{N}\max_{1\leqslant j\leqslant m_{k}}\left\vert q_{n}(j,k)\right\vert
\sum_{j=1}^{m_{k}}\left\vert \alpha_{n}^{(k-1)}\left(  \phi_{(j,k)}\circ
\phi_{(k+1)}\circ...\circ\phi_{(N)}f\right)  \right\vert, \label{FnN} \\
\max_{1\leqslant N\leqslant N_{0}}\left\Vert \digamma_{n}^{(N)}\right\Vert
_{\mathcal{F}} &  \leqslant\frac{2b_{n}}{\sqrt{n}p_{(N_{0})}}S_{N_{0}%
}\max_{1\leqslant k\leqslant N_{0}}\left\Vert \alpha_{n}^{(k-1)}\right\Vert
_{\mathcal{H}_{0}}\leqslant\frac{2b_{n}^{2}S_{N_{0}}}{\sqrt{n}p_{(N_{0}%
)}} \nonumber.
\end{align}
The almost sure result then holds with $C_{0}=8\sigma_{\mathcal{F}}^{2}%
\kappa_{N_{0}}^{2}S_{N_{0}}/p_{(N_{0})}$.\textit{\medskip}

\noindent\textbf{(ii)} We now work on the event $B_{n,N_{0}}$ of (\ref{biendef}). There 
obviously exists $n_{1}$ such that if $n>n_{1}$ then 
$S_{N_{0}}(1-p_{(N_{0})})^{n}\leqslant1/4n^{\theta}$. 
We can also find $\kappa>0$ so small that $n^{2\kappa}/\sqrt
{n}=o(v_{n})$ as $n\rightarrow+\infty$. Therefore, whatever $\zeta>0$ there exists
$n_{2}(\kappa,S_{N_{0}},\zeta,\mathcal{F},P)$ such that 
$\zeta v_{n}>2S_{N_{0}}n^{2\kappa}/p_{(N_{0})}\sqrt{n}$ for any 
$n\geqslant n_{2}$. Choosing $n\geqslant \max(n_{1},n_{2})$ 
we deduce as for (\ref{qn}) and (\ref{FnN}) that%
\begin{align*}
& \mathbb{P}\left(  \max_{0\leqslant N\leqslant N_{0}}\Vert\digamma_{n}%
^{(N)}\Vert_{\mathcal{F}}>\zeta v_{n}\right)  \\
& \leqslant\mathbb{P}\left(  S_{N_{0}}\max_{1\leqslant N\leqslant N_{0}%
}\left(  \left\Vert \alpha_{n}^{(N-1)}\right\Vert _{\mathcal{H}_{0}}%
\max_{1\leqslant j\leqslant m_{N}}\left\vert q_{n}(j,N)\right\vert \right)
>\zeta v_{n}\right)  \\
& \leqslant\mathbb{P}\left(  \left(  \max_{1\leqslant N\leqslant N_{0}}%
\max_{1\leqslant j\leqslant m_{N}}\left\vert q_{n}(j,N)\right\vert \right)
>\frac{2n^{\kappa}}{p_{(N_{0})}\sqrt{n}}\right)  +\mathbb{P}\left(
\max_{1\leqslant N\leqslant N_{0}}\left\Vert \alpha_{n}^{(N-1)}\right\Vert
_{\mathcal{H}_{0}}>n^{\kappa}\right)  \\
& \leqslant2\mathbb{P}\left(  \max_{1\leqslant N\leqslant N_{0}}\left\Vert
\alpha_{n}^{(N-1)}\right\Vert _{\mathcal{H}_{0}}>n^{\kappa}\right)  .
\end{align*}
By Proposition \ref{ProConcentr} we see that under (VC) or (BR) the latter
probability can be made less than $1/8n^{\theta}$ for any $n>n_3(\zeta,\theta)$ and 
$n_3(\zeta,\theta)$ large enough. Clearly $n_3(\zeta,\theta)$ depends on $\zeta, 
\theta, n_1, n_2$ and on the entropy of $\mathcal{H}_{0}$ thus all constants 
in Lemma \ref{LemEntrop} and Proposition \ref{ProConcentr} are involved.$\quad\square$\medskip

\noindent\textbf{Step 3.} Fix $\theta>0$. By Lemma \ref{LemEntrop} we can
apply Propositions 1 and 2 of Berthet and Mason~\cite{BerMas06} to $\mathcal{F}_{0}$, which ensures the following
Gaussian approximation. For some constant $c_{\theta}(\mathcal{F}_{0},P)>0$
and $n_{\theta}(\mathcal{F}_{0},P)>0$ we can build on a probability space
$(\Omega,\mathcal{T},\mathbb{P})$ a version of the sequence $\left\{
X_{n}\right\}$ of independent random variables with law $P$ and a sequence $\{\mathbb{G}_{n}^{(0)}\}$ of coupling
versions of $\mathbb{G}^{(0)}$ in such a way that, for all 
$n\geqslant n_{\theta}(\mathcal{F}_{0},P)$,%
\begin{equation}
\mathbb{P}\left(  \left\Vert \alpha_{n}^{(0)}-\mathbb{G}_{n}^{(0)}\right\Vert
_{\mathcal{F}_{0}}>c_{\theta}(\mathcal{F}_{0},P)v_{n}\right)  \leqslant
\frac{1}{2n^{\theta}}.\label{BM}%
\end{equation}
Keep in mind that constants $n_{\theta}$ and $c_{\theta}$ only depend on the entropy of $\mathcal{F}_{0}$ through the constants $ M, c_0, \nu_0, b_0, r_0 $. By choosing $\theta>1$, $d_{\theta}>c_{\theta}(\mathcal{F}_{0},P)$ then
applying Borel-Cantelli lemma to (\ref{BM}), it almost surely holds, for all
$n$ large enough,%
\begin{equation}
\left\Vert \alpha_{n}^{(0)}-\mathbb{G}_{n}^{(0)}\right\Vert _{\mathcal{F}_{0}%
}<d_{\theta}v_{n}.\label{BMps}%
\end{equation}

\noindent\textbf{Step 4.} 
Let $\theta_0>0$. We work under the event $B_{n,N_{0}}$ of (\ref{biendef}) with a
probability at least $1-1/4n^{\theta_0}$ provided that $n>n_{1}$. 
The process $\mathbb{G}^{(0)}$ being linear on $\mathcal{F}$ we see that the
recursive definition (\ref{GNf}) applied to the version $\mathbb{G}_{n}^{(0)}$
of $\mathbb{G}^{(0)}$ reads $\mathbb{G}_{n}^{(N)}(f)=\mathbb{G}_{n}%
^{(N-1)}(\phi_{(N)}f)$. This combined with (\ref{alphagamma}) readily gives%
\begin{align}
& \max_{1\leqslant N\leqslant N_{0}}\left\Vert \alpha_{n}^{(N)}-\mathbb{G}%
_{n}^{(N)}\right\Vert _{\mathcal{F}} \nonumber\\
& =\max_{1\leqslant N\leqslant N_{0}%
}\left\Vert \alpha_{n}^{(N-1)}(\phi_{(N)}f)-\mathbb{G}_{n}^{(N-1)}%
(\phi_{(N)}f)+\Gamma_{n}^{(N)}(f)\right\Vert _{\mathcal{F}}\nonumber\\
&  =\max_{1\leqslant N\leqslant N_{0}}\left\Vert \alpha_{n}^{(0)}(\phi
_{(1)}\circ...\circ\phi_{(N)}f)-\mathbb{G}_{n}^{(0)}(\phi_{(1)}\circ
...\circ\phi_{(N)}f)+\digamma_{n}^{(N)}(f)\right\Vert _{\mathcal{F}%
}\nonumber\\
&  \leqslant\left\Vert \alpha_{n}^{(0)}-\mathbb{G}_{n}^{(0)}\right\Vert
_{\mathcal{F}_{0}}+\max_{0\leqslant N\leqslant N_{0}}\Vert\digamma_{n}%
^{(N)}\Vert_{\mathcal{F}}.\label{approxN0}%
\end{align}
Remind that $v_{n}>L(n)/\sqrt{n}$ and Lemma \ref{LemGamma} holds. By
(\ref{BMps}) and (\ref{approxN0}) we almost surely have, for all $n$ large
enough and $d_0=2d_{\theta_0}$,%
\[
\max_{1\leqslant N\leqslant N_{0}}\left\Vert \alpha_{n}^{(N)}-\mathbb{G}%
_{n}^{(N)}\right\Vert _{\mathcal{F}}\leqslant d_{\theta_0}v_{n}+C_{0}%
\frac{L\circ L(n)}{\sqrt{n}}\leqslant d_{0}v_{n}.
\]
By Lemmas \ref{LemEntrop} and \ref{LemGamma}, (\ref{BM}) and (\ref{approxN0}), for $n_0>
\max(n_1,n_3(\zeta,\theta_0),n_{\theta_0}(\mathcal{F}_{0},P))$ and $d_{0}>c_{\theta_0}(\mathcal{F}_{0},P)+\zeta$ 
we have, for all $n\geqslant n_0$,%
\[
\mathbb{P}\left(  \max_{1\leqslant N\leqslant N_{0}}\left\Vert \alpha
_{n}^{(N)}-\mathbb{G}_{n}^{(N)}\right\Vert _{\mathcal{F}}>d_{0}v_{n}\right)
\leqslant\frac{1}{2n^{\theta_0}}+\mathbb{P}\left(  \max_{0\leqslant N\leqslant
N_{0}}\Vert\digamma_{n}^{(N)}\Vert_{\mathcal{F}}>\zeta v_{n}\right)
\leqslant\frac{1}{n^{\theta_0}}.
\]
To conclude observe that the parameters $N_{0}, M, M_{N_{0}}, S_{N_{0}}, P_{N_{0}},
p_{(N_{0})}, \theta_{0}, \nu_0, c_0, r_0, b_0$ have been used at one or several steps to 
finally define $n_0$ and $d_0$.

\subsection{Proof of Proposition~\ref{ProBias}}

\noindent Theorem \ref{ProCVSpeed} implies, for
$f\in\mathcal{F}$,%
\begin{equation}
\mathbb{P}_{n}^{(N)}(f)-P(f)=\frac{1}{\sqrt{n}}\mathbb{G}_{n}^{(N)}%
(f)+\frac{1}{\sqrt{n}}\mathbb{R}_{n}^{(N)}(f), \label{decompo}%
\end{equation}
where $\mathbb{G}_{n}^{(N)}$ is a sequence of versions of the centered
Gaussian process $\mathbb{G}^{(N)}$ from (\ref{GNf}) and the random sequence
$r_{n}^{(N)}=\Vert \mathbb{R}_{n}^{(N)}\Vert _{\mathcal{F}}$
satisfies%
\[
r_{n}^{(N)}\leqslant\left\Vert \mathbb{G}_{n}^{(N)}\right\Vert _{\mathcal{F}%
}+\left\Vert \alpha_{n}^{(N)}\right\Vert _{\mathcal{F}}\leqslant\left\Vert
\mathbb{G}_{n}^{(N)}\right\Vert _{\mathcal{F}}+2M\sqrt{n},\quad\lim
_{n\rightarrow+\infty}\frac{r_{n}^{(N)}}{v_{n}}\leqslant d_{0}\quad a.s.
\]
We have to be a little careful with the expectation, variance and 
covariance of the coupling error process $\mathbb{R}_{n}^{(N)}$.\medskip

\noindent\textbf{Step 1.}\ Since $\mathbb{G}_{n}^{(N)}(f)$ is centered the bias is controlled by%
\begin{equation}
\sup_{f\in\mathcal{F}}\frac{\sqrt{n}}{v_{n}}\left\vert \mathbb{E}\left(
\mathbb{P}_{n}^{(N)}(f)\right)  -P(f)\right\vert =\sup_{f\in\mathcal{F}%
}\left\vert \frac{1}{v_{n}}\mathbb{E}\left(  \mathbb{R}_{n}^{(N)}(f)\right)
\right\vert \leqslant\mathbb{E}\left(  \frac{r_{n}^{(N)}}{v_{n}}\right).
\label{bias}%
\end{equation}
Write $a_{n}=\sqrt{K\log n}$ where $K>0$ and $\theta_{0}>1$ from Theorem \ref{ProCVSpeed} can be chosen as large as needed. Then, for $\theta>1$, $\varepsilon>0$ and $k\in\mathbb{N}^{\ast
}$ consider the events%
\begin{align*}
A_{n} & =\left\{  r_{n}^{(N)}\leqslant(d_{0}+\varepsilon)v_{n}\right\}  ,\quad
B_{n}=\left\{  \left\Vert \mathbb{G}_{n}^{(N)}\right\Vert _{\mathcal{F}%
}\leqslant a_{n}\right\}  ,\\
C_{n,k} & =\left\{  \theta^{k-1}a_{n}<\left\Vert
\mathbb{G}_{n}^{(N)}\right\Vert _{\mathcal{F}}\leqslant\theta^{k}%
a_{n}\right\}  .
\end{align*}
By Theorem \ref{ProCVSpeed}, $\mathbb{P}(A_{n}^{c})<1/n^{\theta_{0}}$ and
$v_{n}>a_{n}/\sqrt{n}$ for all $n$ large enough, hence%
\begin{align*}
\frac{1}{v_{n}}\mathbb{E}\left(  r_{n}^{(N)}\right)   &  =\mathbb{E}\left(
\frac{r_{n}^{(N)}}{v_{n}}1_{A_{n}}\right)  +\mathbb{E}\left(  \frac
{r_{n}^{(N)}}{v_{n}}1_{A_{n}^{c}\cap B_{n}}\right)  +\mathbb{E}\left(
\frac{r_{n}^{(N)}}{v_{n}}1_{A_{n}^{c}\cap B_{n}^{c}}\right) \\
&  \leqslant d_{0}+\varepsilon+\frac{a_{n}+2M\sqrt{n}}{v_{n}}\mathbb{P}%
(A_{n}^{c})+\mathbb{E}\left(  \frac{r_{n}^{(N)}}{v_{n}}1_{B_{n}^{c}}\right) \\
&  \leqslant d_{0}+2\varepsilon+%
{\displaystyle\sum\limits_{k=1}^{+\infty}}
\mathbb{E}\left(  \frac{r_{n}^{(N)}}{v_{n}}1_{C_{n,k}}\right)  .
\end{align*}
By Propositions \ref{Thm2} and \ref{Thm2b}, $\mathbb{G}^{(N)}(f)$ is a centered
Gaussian process indexed by $\mathcal{F}$ such that, under (VC) or (BR), it holds
\begin{align}
\mathbb{E}\left(  \left\Vert \mathbb{G}^{(N)}\right\Vert _{\mathcal{F}%
}\right)  & <+\infty, \quad
\sup_{f\in\mathcal{F}}\mathbb{V}(\mathbb{G}^{(N)}(f))\leqslant
\sigma_{\mathcal{F}}^{2}<+\infty,\nonumber\\
\mathbb{E}\left(  \left\Vert \mathbb{G}^{(N)}\right\Vert _{\mathcal{F}}%
^{2}\right)   & \leqslant C_{\mathcal{F}}^{2}=\sigma_{\mathcal{F}}%
^{2}+\mathbb{E}\left(  \left\Vert \mathbb{G}^{(N)}\right\Vert _{\mathcal{F}%
}\right)  ^{2}<+\infty.
\label{SupVarGN}
\end{align}
Thus, by Borell's inequality -- see Appendix A.2 of \cite{VanWell} -- for any version $\mathbb{G}_{n}^{(N)}$ 
of $\mathbb{G}^{(N)}$, we have
\begin{equation}
\mathbb{P}\left(  \left\Vert \mathbb{G}_{n}^{(N)}\right\Vert _{\mathcal{F}%
}>\lambda\right)  \leqslant2\exp\left(  -\frac{\lambda^{2}}{8C_{\mathcal{F}%
}^{2}}\right)  . \label{Borell}%
\end{equation}
Therefore we have, since $\theta>1$ and $v_{n}>4M/\sqrt{n}>2a_{n}/n$ for $n$
large enough,%
\begin{align*}
\mathbb{E}\left(  \frac{r_{n}^{(N)}}{v_{n}}1_{C_{n,k}}\right)   &   \leqslant
\frac{\theta^{k}a_{n}+2M\sqrt{n}}{v_{n}}\mathbb{P}\left(  C_{n,k}\right)
\leqslant\theta^{k}n\mathbb{P}\left(  \left\Vert \mathbb{G}_{n}^{(N)}%
\right\Vert _{\mathcal{F}}>\theta^{k-1}a_{n}\right)\\
&   \leqslant2\theta^{k}%
n\exp\left(  -\frac{(\theta^{k-1}a_{n})^{2}}{8C_{\mathcal{F}}^{2}}\right),
\end{align*}
and the following series is converging to an arbitrarily small sum,%
\begin{align*}%
{\displaystyle\sum\limits_{k=1}^{+\infty}}
\mathbb{E}\left(  \frac{r_{n}^{(N)}}{v_{n}}1_{C_{n,k}}\right)   &
\leqslant2n\exp\left(  -\frac{a_{n}^{2}}{8C_{\mathcal{F}}^{2}}\right)
{\displaystyle\sum\limits_{k=1}^{+\infty}}
\theta^{k}\exp\left(  -\left(  \frac{\theta^{2(k-1)}-1}{8C_{\mathcal{F}}^{2}%
}\right)  a_{n}^{2}\right) \\
&  \leqslant n\exp\left(  -\frac{K\log n}{8C_{\mathcal{F}}^{2}}\right)
{\displaystyle\sum\limits_{k=1}^{+\infty}}
2e\theta^{k}\exp\left(  -\theta^{2(k-1)}\right)  \leqslant\frac{1}{n^{\delta}},%
\end{align*}
where $\delta<K/8C_{\mathcal{F}}^{2}-1$. It follows that (\ref{bias}) is
ultimately bounded by $d_{0}$.\medskip

\noindent\textbf{Step 2.}\ Starting from (\ref{decompo}) and the bias and
variance decomposition, the quadratic risk is in turn controlled by%
\begin{align}
&  \mathbb{E}\left(  ( \mathbb{P}_{n}^{(N)}(f)-P(f) ) ^{2}\right)
-\frac{1}{n}\mathbb{V}\left(  \mathbb{G}^{(N)}(f)\right)  \nonumber\\
&  =\left\vert \mathbb{E}\left(  \mathbb{P}_{n}^{(N)}(f)\right)
-P(f)\right\vert ^{2}+\frac{1}{n}\mathbb{V}\left(  \mathbb{R}_{n}^{(N)}(f)\right)
+\frac{2}{n}\mathrm{Cov}\left(  \mathbb{G}_{n}^{(N)}(f),\mathbb{R}_{n}^{(N)}(f)\right)
.\label{decompoRisk}%
\end{align}
\textbf{(i)} By Step 1, the first right-hand term is the squared bias, of
order $d_{0}^{2}v_{n}^{2}/n$. Concerning the second right-hand
term in (\ref{decompoRisk}), we bound 
$\mathbb{E}(\mathbb{R}_{n}^{(N)}(f)^{2})$.
Fix $\varepsilon>0$ and assume that $n$ is large enough for the following statements.
By setting $s_{n}^{(N)}=(r_{n}^{(N)})^{2}$ then using 
$v_{n}>a_{n}/\sqrt{n}$, $a_{n}=K\sqrt{\log n}<\sqrt{n}$ we
get, for $\theta_{0}=2$,%
\begin{align*}
\frac{1}{v_{n}^{2}}\sup_{f\in\mathcal{F}}\mathbb{E}\left(  \mathbb{R}_{n}%
^{(N)}(f)^{2}\right)   &  \leqslant\mathbb{E}\left(  \frac{s_{n}^{(N)}}%
{v_{n}^{2}}1_{A_{n}}\right)  +\mathbb{E}\left(  \frac{s_{n}^{(N)}}{v_{n}^{2}%
}1_{A_{n}^{c}\cap B_{n}}\right)  +\mathbb{E}\left(  \frac{s_{n}^{(N)}}%
{v_{n}^{2}}1_{A_{n}^{c}\cap B_{n}^{c}}\right)  \\
&  \leqslant(d_{0}+\varepsilon)^{2}+\left(  \frac{a_{n}+2M\sqrt{n}}{v_{n}%
}\right)  ^{2}\mathbb{P}(A_{n}^{c})+\mathbb{E}\left(  \frac{s_{n}^{(N)}}%
{v_{n}^{2}}1_{B_{n}^{c}}\right)  \\
&  \leqslant(d_{0}+2\varepsilon)^{2}+\left(  \frac{3M\sqrt{n}}{\log n}\right)  ^{2}\frac
{1}{n^{2}}+%
{\displaystyle\sum\limits_{k=1}^{+\infty}}
\mathbb{E}\left(  \frac{s_{n}^{(N)}}{v_{n}^{2}}1_{C_{n,k}}\right)  \\
&  \leqslant(d_{0}+3\varepsilon)^{2}+%
{\displaystyle\sum\limits_{k=1}^{+\infty}}
\theta^{2k}n^{2}\mathbb{P}\left(  \left\Vert \mathbb{G}_{n}^{(N)}\right\Vert
_{\mathcal{F}}>\theta^{k-1}a_{n}\right)  \leqslant(d_{0}+4\varepsilon)^{2},%
\end{align*}
where the series is equal to its first term $n^{2}\exp\left(  -a_{n}%
^{2}/8C_{\mathcal{F}}^{2}\right)  $ times a convergent series, by using
(\ref{Borell}) as for Step 1 with $K>16C_{\mathcal{F}}^{2}$.
We have shown that
\begin{equation}
\limsup_{n\rightarrow+\infty}\frac{1}{v_{n}^{2}}\sup_{f\in\mathcal{F}}%
\mathbb{V}\left(\mathbb{R}_{n}^{(N)}(f)\right) \leqslant
\limsup_{n\rightarrow+\infty}\frac{1}{v_{n}^{2}}
\mathbb{E}\left(s_{n}^{(N)}\right)
\leqslant d_{0}^{2}.\label{Varvn}%
\end{equation}

\noindent\textbf{(ii)} Concerning the covariance term in 
(\ref{decompoRisk}) it holds%
\begin{align*}
\frac{1}{v_{n}}\left\vert \mathrm{Cov}\left(\mathbb{G}_{n}^{(N)}(f),\mathbb{R}_{n}%
^{(N)}(f)\right)  \right\vert  &  =\frac{1}{v_{n}}\left\vert \mathbb{E}\left(
\mathbb{G}_{n}^{(N)}(f)\mathbb{R}_{n}^{(N)}(f)\right)  \right\vert \\
&  \leqslant \frac{1}{v_{n}} \mathbb{E}\left(\left\vert
\mathbb{G}_{n}^{(N)}(f)\right\vert r_{n}^{(N)}\right)    \\
&  = T_{A_{n}}(f)+T_{A_{n}^{c}\cap B_{n}}(f)+T_{A_{n}^{c}\cap
B_{n}^{c}}(f),
\end{align*}
where $
T_{D}(f)=\mathbb{E}(  1_{D} \vert \mathbb{G}_{n}^{(N)}(f)\vert
r_{n}^{(N)}/v_{n}) $ for $ D\in\left\{  A_{n},A_{n}^{c}\cap
B_{n},A_{n}^{c}\cap B_{n}^{c}\right\}$.
We have, by Proposition~\ref{Thm2} and~\ref{Thm2b},
\begin{align*}
    T_{A_{n}}(f)&\leqslant\mathbb{E}\left(  \left\vert \mathbb{G}_{n}^{(N)}%
(f)\right\vert (d_{0}+\varepsilon)1_{A_{n}}\right) \\ 
&=(d_{0}%
+\varepsilon) \sqrt{\mathbb{V}\left(\mathbb{G}_n^{(N)}(f)\right)} \mathbb{E}\left(  \left\vert \mathcal{N}%
(0,1)\right\vert \right) \leqslant \sqrt{\frac{2}{\pi}}(d_{0}+\varepsilon)\sigma_{f} .
\end{align*}
By using again $\mathbb{P}(A_{n}^{c})<1/n^{2}$ we see that%
\[
T_{A_{n}^{c}\cap B_{n}}(f)\leqslant\mathbb{E}\left(  a_{n}\left(
\frac{2M\sqrt{n}+a_{n}}{v_{n}}\right)  1_{A_{n}^{c}\cap B_{n}}\right)
\leqslant a_{n}\left(  \frac{3M\sqrt{n}}{v_{n}}\right)  \frac{1}{n^{2}%
}\leqslant\varepsilon.
\]
Lastly, for $g_{n}^{(N)}=\left\Vert \mathbb{G}_{n}^{(N)}\right\Vert
_{\mathcal{F}}$, $K$ large and all $n$ large enough it holds, by
(\ref{SupVarGN}) and (\ref{Borell}),%
\begin{align*}
T_{A_{n}^{c}\cap B_{n}^{c}}(f) &  \leqslant\mathbb{E}\left(  g_{n}%
^{(N)}\left(  \frac{2M\sqrt{n}+g_{n}^{(N)}}{v_{n}}\right)  1_{B_{n}^{c}%
}\right)  \\
&  =%
{\displaystyle\sum\limits_{k=1}^{+\infty}}
\mathbb{E}\left(  g_{n}^{(N)}\left(  \frac{2M\sqrt{n}+g_{n}^{(N)}}{v_{n}%
}\right)  1_{C_{n,k}}\right)  \\
&  \leqslant%
{\displaystyle\sum\limits_{k=1}^{+\infty}}
\theta^{2k}a_{n}^{2}n\mathbb{P}\left(  \left\Vert \mathbb{G}_{n}%
^{(N)}\right\Vert _{\mathcal{F}}>\theta^{k-1}a_{n}\right)  \leqslant
\varepsilon.
\end{align*}
The above upper bounds imply, by (\ref{decompoRisk}),
(\ref{Varvn}) and (\ref{sigma}),
\[
\limsup_{n\rightarrow+\infty}\frac{n}{v_{n}}\sup_{f\in\mathcal{F}}\left\vert
\mathbb{E}\left( ( \mathbb{P}_{n}^{(N)}(f)-P(f) ) ^{2}\right)
-\frac{1}{n}\mathbb{V}\left(  \mathbb{G}^{(N)}(f)\right)  \right\vert
\leqslant\sqrt{\frac{8}{\pi}}d_{0}\sigma_{\mathcal{F}}.
\]

\noindent\textbf{Step 3.}\ Let extend Step 2 to the covariance. By Step 1 we have, for all $n$ large,%
\begin{align*}
& \left\vert \mathrm{Cov}\left(  \mathbb{P}_{n}^{(N)}(f),\mathbb{P}_{n}^{(N)}(g)\right)
-\mathbb{E}\left(  \left(  \mathbb{P}_{n}^{(N)}(f)-P(f)\right)  \left(
\mathbb{P}_{n}^{(N)}(g)-P(g)\right)  \right)  \right\vert \\
& =\left\vert \left( \mathbb{E}(  \mathbb{P}_{n}^{(N)}(f))
-P(f) \right) \left( \mathbb{E}(  \mathbb{P}_{n}^{(N)}(g))  -P(g)\right)\right\vert
<d_{0}^2\frac{v_{n}^{2}}{n}.%
\end{align*}
Now, by the upper bounds computed at (i) and (ii) of Step 2,%
\begin{align*}
&  \left\vert \mathbb{E}\left( (\mathbb{P}_{n}^{(N)}(f)-P(f))
(\mathbb{P}_{n}^{(N)}(g)-P(g)) \right) -\frac{1}{n}\mathrm{Cov}\left(
\mathbb{G}^{(N)}(f),\mathbb{G}^{(N)}(g)\right)  \right\vert \\
&  \leqslant\frac{1}{n}\mathbb{E}\left(  \left\vert \mathbb{G}_{n}%
^{(N)}(f)\mathbb{R}_{n}^{(N)}(g)\right\vert \right)  +\frac{1}{n}%
\mathbb{E}\left(  \left\vert \mathbb{G}_{n}^{(N)}(g)\mathbb{R}_{n}%
^{(N)}(f)\right\vert \right)  +\frac{1}{n}\mathbb{E}\left(  \left\vert
\mathbb{R}_{n}^{(N)}(f)\mathbb{R}_{n}^{(N)}(g)\right\vert \right)  \\
&  \leqslant\frac{2}{n}\sup_{f\in\mathcal{F}}\mathbb{E}\left(  \left\vert
\mathbb{G}_{n}^{(N)}(f)\right\vert r_{n}^{(N)}\right)  +\frac{1}{n}%
\mathbb{E}\left(  s_{n}^{(N)}\right)  \\
&  \leqslant\frac{2}{n}(d_{0}+\varepsilon)\sigma_{\mathcal{F}}\sqrt{\frac
{2}{\pi}}v_{n}+\frac{1}{n}(d_{0}+\varepsilon)^{2}v_{n}^{2}.
\end{align*}

\subsection{Proof of Proposition~\ref{ProBerryEsseen}.}

Fix $N_{0}\in\mathbb{N}$. Let apply Theorem \ref{ProCVSpeed}
with $\theta_0=2$, from which we also use $n_{0}$, $d_{0}$ and $v_{n}$. We have, for all $0\leqslant
N\leqslant N_{0}$, $\varphi\in\mathcal{L}$, $x\in\mathbb{R}$ and $n\geqslant
n_{0}$,%
\begin{align*}
\mathbb{P}\left(  \varphi(\alpha_{n}^{(N)})\leqslant x\right)   &
\leqslant\frac{1}{n^{2}}+\mathbb{P}\left(  \left\{  \varphi(\alpha_{n}%
^{(N)})\leqslant x\right\}  \cap\left\{  \left\Vert \alpha_{n}^{(N)}%
-\mathbb{G}_{n}^{(N)}\right\Vert _{\mathcal{F}}<d_{0}v_{n}\right\}  \right)
\\
&  \leqslant\frac{1}{n^{2}}+\mathbb{P}\left(  \varphi(\mathbb{G}_{n}%
^{(N)})\leqslant x+d_{0}C_{1}v_{n}\right)  \\
&  \leqslant\frac{1}{n^{2}}+\mathbb{P}\left(  \varphi(\mathbb{G}_{n}%
^{(N)})\leqslant x\right)  +d_{0}C_{1}C_{2}v_{n},%
\end{align*}
and%
\begin{align*}
& \mathbb{P}\left(  \varphi(\mathbb{G}_{n}^{(N)})\leqslant x-d_{0}C_{1}%
v_{n}\right)  \\
& \leqslant\frac{1}{n^{2}}+\mathbb{P}\left(  \left\{
\varphi(\mathbb{G}_{n}^{(N)})\leqslant x-d_{0}C_{1}v_{n}\right\}  \cap\left\{
\left\Vert \alpha_{n}^{(N)}-\mathbb{G}_{n}^{(N)}\right\Vert _{\mathcal{F}%
}<d_{0}v_{n}\right\}  \right)  \\
&  \leqslant\frac{1}{n^{2}}+\mathbb{P}\left(  \varphi(\alpha_{n}%
^{(N)})\leqslant\varphi(\mathbb{G}_{n}^{(N)})+d_{0}C_{1}v_{n}\leqslant
x\right),
\end{align*}
so that%
\[
\mathbb{P}\left(  \varphi(\alpha_{n}^{(N)})\leqslant x\right)  \geqslant
\mathbb{P}\left(  \varphi(\mathbb{G}_{n}^{(N)})\leqslant x\right)  -d_{0}%
C_{1}C_{2}v_{n}-\frac{1}{n^{2}}.
\]
This establishes the second statement of Proposition \ref{ProBerryEsseen}
provided $d_{1}>d_{0}$ and $n\geqslant n_{1}\geqslant n_{0}$ where $n_{1}$ is
large enough to have $(d_{1}-d_{0} C_1 C_2)v_{n}>n^{-2}$. The first statement
coincides with the special case $\mathcal{L}=\left\{  \varphi_{f}%
:f\in\mathcal{F}_{0}\right\}  $ where $\varphi_{f}(g)=g(f)$ are pointwise
projectors and we then have a Lipshitz constant $C_{1}=1$ whereas $\varphi
_{f}(\mathbb{G}_{n}^{(N)})=\mathbb{G}_{n}^{(N)}(f)$ has a Gaussian density
bounded by
\[
\frac{1}{\sqrt{2\pi\mathbb{V}(\mathbb{G}_{n}^{(N)}(f))}}\leqslant C_{2}%
=\frac{1}{\sqrt{2\pi}\sigma_{0}}<+\infty.
\]

\section{Proofs concerning the limiting process}\label{proofsLP}

\subsection{Proof of Proposition~\ref{Thm2}}

\textbf{Step 1.} Let us first relate $\mathbb{G}^{(N)}(\mathcal{F})$ from (\ref{GNf})
to $\mathbb{G}(\mathcal{F})=\mathbb{G}^{(0)}(\mathcal{F})$ from (\ref{GF}) by means 
of the vectors $\Phi_{k}^{(N)}(f)$ introduced at (\ref{phikNf}) before
Proposition~\ref{Thm2}.

\begin{lem}
\label{Thm1}For all $N\in\mathbb{N}^{\ast}$ and $f\in\mathcal{F}$ it holds
\[
\mathbb{G}^{(N)}(f)=\mathbb{G}(f)-\sum_{k=1}^{N}\Phi_{k}%
^{(N)}(f)^{t}\cdot\mathbb{G}\left[\mathcal{A}^{(k)}\right].
\]
\end{lem}

\noindent\textit{Proof.} The formula is true for $N=0$. Assume that it is the
case for $N\geqslant0$. Recall that sets $A\in\mathcal{A}_{\cup}^{(N_{0})}$
from (1.9) are identified to $f=1_{A}$. By (\ref{GNf}),%
\begin{align*}
\mathbb{G}^{(N+1)}(f) &  =\mathbb{G}^{(N)}(f)-\mathbb{E}\left[f\mid\mathcal{A}%
^{(N+1)}\right]^{t}\cdot\mathbb{G}^{(N)}\left[\mathcal{A}^{(N+1)}\right]\\
&  =\mathbb{G}(f)-\sum_{k=1}^{N}\Phi_{k}^{(N)}(f)^{t}\cdot\mathbb{G}%
\left[\mathcal{A}^{(k)}\right]\\
&  -\sum_{j=1}^{m_{N+1}}\mathbb{E}(f\mid A_{j}^{(N+1)})  \left(
\mathbb{G}(A_{j}^{(N+1)})-\sum_{k=1}^{N}\Phi_{k}^{(N)}(A_{j}^{(N+1)}%
)^{t}\cdot\mathbb{G}\left[\mathcal{A}^{(k)}\right]\right)  \\
&  =\mathbb{G}(f)-\mathbb{E}\left[f\mid\mathcal{A}^{(N+1)}\right]^{t}\cdot\mathbb{G}%
\left[\mathcal{A}^{(N+1)}\right]\\
&  -\sum_{k=1}^{N}\left[  \Phi_{k}^{(N)}(f)-\Phi_{k}^{(N)}[\mathcal{A}%
^{(N+1)}]\cdot\mathbb{E}[f\mid\mathcal{A}^{(N+1)}]\right]  ^{t}\cdot
\mathbb{G}\left[\mathcal{A}^{(k)}\right],
\end{align*}
where the $m_{k}\times
m_{N+1}$ matrix $\Phi_{k}^{(N)}[\mathcal{A}^{(N+1)}]=(\Phi_{k}^{(N)}%
(A_{1}^{(N+1)}),\dots,\Phi_{k}^{(N)}(A_{m_{N+1}}^{(N+1)}))$ satisfies%
\[
\sum_{j=1}^{m_{N+1}}\mathbb{E}(  f\mid A_{j}^{(N+1)}) \Phi
_{k}^{(N)}(A_{j}^{(N+1)})^{t}=\left[  \Phi_{k}^{(N)}[\mathcal{A}%
^{(N+1)}]\cdot\mathbb{E}[f\mid\mathcal{A}^{(N+1)}]\right] ^{t}.
\]
Now observe that $\Phi_{N+1}^{(N+1)}(f)=\mathbb{E}[f\mid\mathcal{A}^{(N+1)}]$ by the definition of $ \Phi_{k}^{(N+1)} $ given by (\ref{phikNf}) when $ k=N+1$. 
It remains to show that%
\begin{equation}
\Phi_{k}^{(N)}(f)-\Phi_{k}^{(N)}[\mathcal{A}^{(N+1)}]\cdot\mathbb{E}%
[f\mid\mathcal{A}^{(N+1)}]=\Phi_{k}^{(N+1)}(f).\label{phikN}%
\end{equation}
For $1\leqslant k\leqslant N$ and $1\leqslant j\leqslant m_{N+1}$ we have%
\begin{align*}
& \Phi_{k}^{(N)}(A_{j}^{(N+1)}) =\mathbb{E}\left[A_{j}^{(N+1)}\mid\mathcal{A}%
^{(k)}\right]+ \\
& \sum_{\substack{1\leqslant L\leqslant N-k\\k<l_{1}%
<...<l_{L}\leqslant N}}(-1)^{L}\mathbf{P}_{\mathcal{A}^{(l_{1})}%
|\mathcal{A}^{(k)}}\mathbf{P}_{\mathcal{A}^{(l_{2})}|\mathcal{A}%
^{(l_{1})}}\dots\mathbf{P}_{\mathcal{A}^{(l_{L})}|\mathcal{A}%
^{(l_{L-1})}}\cdot\mathbb{E}\left[A_{j}^{(N+1)}|\mathcal{A}^{(l_{L})}\right],
\end{align*}
where, for $l=k,k+1,...,N$ the vector%
\[
\mathbb{E}\left[A_{j}^{(N+1)}\mid\mathcal{A}^{(l)}\right]=\left(  \frac{P(
A_{j}^{(N+1)}\cap A_{1}^{(l)})  }{P(  A_{1}^{(l)})
},...,\frac{P( A_{j}^{(N+1)}\cap A_{m_{l}}^{(l)}) }{P(
A_{m_{l}}^{(l)}) }\right) ^{t},%
\]
is also the $j$-th column of $\mathbf{P}_{\mathcal{A}^{(N+1)}|\mathcal{A}%
^{(l)}}$. Therefore, by turning $L$ into $L^{\prime}=L+1$,%
\begin{align*}
& -\Phi_{k}^{(N)}\left[\mathcal{A}^{(N+1)}\right]\cdot\mathbb{E}\left[  f\mid
\mathcal{A}^{(N+1)}\right]  \\
& =-\sum_{j=1}^{m_{N+1}}\mathbb{E}(f\mid A_{j}^{(N+1)})  \Phi
_{k}^{(N)}(A_{j}^{(N+1)})\\
& =-\sum_{j=1}^{m_{N+1}}\mathbb{E}(f\mid A_{j}^{(N+1)})
\mathbb{E}\left[A_{j}^{(N+1)}\mid\mathcal{A}^{(k)}\right]\\
& +\sum_{j=1}^{m_{N+1}}\mathbb{E}(f\mid A_{j}^{(N+1)})
\sum_{\substack{1\leqslant L\leqslant N-k\\k<l_{1}<...<l_{L}\leqslant
N}}(-1)^{L+1}\mathbf{P}_{\mathcal{A}^{(l_{1})}|\mathcal{A}^{(k)}}%
\dots\mathbf{P}_{\mathcal{A}^{(l_{L})}|\mathcal{A}^{(l_{L-1})}}\cdot
\mathbb{E}\left[A_{j}^{(N+1)}|\mathcal{A}^{(l_{L})}\right]\\
& =(-1)^{1}\mathbf{P}_{\mathcal{A}^{(N+1)}|\mathcal{A}^{(k)}}\mathbb{E}%
\left[f|\mathcal{A}^{(N+1)}\right]\\
& +\sum_{\substack{1\leqslant L^{\prime}\leqslant N+1-k\\k<l_{1}%
<...<l_{L^{\prime}}=N+1}}(-1)^{L^{\prime}}\mathbf{P}_{\mathcal{A}^{(l_{1}%
)}|\mathcal{A}^{(k)}}\dots\mathbf{P}_{\mathcal{A}^{(l_{L^{\prime}-1}%
)}|\mathcal{A}^{(l_{L^{\prime}-2})}}\mathbf{P}_{\mathcal{A}%
^{(N+1)}|\mathcal{A}^{(l_{L^{\prime}-1})}}\cdot\mathbb{E}\left[f|\mathcal{A}%
^{(N+1)}\right],
\end{align*}
where all terms are different from those in
\[
\Phi_{k}^{(N)}(f)=\mathbb{E}\left[  f\mid\mathcal{A}^{(k)}\right]
+\sum_{\substack{1\leqslant L^{\prime}\leqslant N+1-k\\k<l_{1}%
<...<l_{L^{\prime}}<N+1}}(-1)^{L}\mathbf{P}_{\mathcal{A}^{(l_{1})}%
|\mathcal{A}^{(k)}}\dots\mathbf{P}_{\mathcal{A}^{(l_{L^{\prime}}%
)}|\mathcal{A}^{(l_{L^{\prime}-1})}}\cdot\mathbb{E}\left[f|\mathcal{A}%
^{(l_{L^{\prime}})}\right].
\]
Having collected all terms of $\Phi_{k}^{(N)}(f)$ in (\ref{phikNf}), this establishes (\ref{phikN}). The proof is completed by induction.$\quad\square$\medskip

\noindent The functions $\Phi_{k}^{(N)}$ and the process $\mathbb{G}$ are linear, hence
Lemma \ref{Thm1} implies that $\mathbb{G}^{(N)}$ is a linear process. Moreover
$\mathbb{G}(f)$ and $\mathbb{G}[\mathcal{A}^{(k)}]$ being centered Gaussian,
Lemma \ref{Thm1} proves that $\mathbb{G}^{(N)}(f)$ is a centered Gaussian
random variable.\medskip

\noindent\textbf{Step 2.} To compute the covariance of 
$\mathbb{G}^{(N)}(\mathcal{F})$ we need the following properties. Recall 
that $\mathbf{P}_{\mathcal{A}^{(k)}|\mathcal{A}^{(k)}%
}=\mathrm{Id}_{m_k}$ is the identity matrix of $\mathbb{R}^{m_k}$.

\begin{lem}
\label{covariance} For $1\leqslant k,l\leqslant N$ and $f\in\mathcal{F}$ we have%
\begin{align}
\mathrm{Cov}\left(  \mathbb{G}[\mathcal{A}^{(k)}],\mathbb{G}(f)\right)    &
=\mathbb{V}\left(  \mathbb{G}[\mathcal{A}^{(k)}]\right)  \cdot\mathbb{E}%
[f\mid\mathcal{A}^{(k)}],\label{VE}\\
\mathrm{Cov}\left(  \mathbb{G}[\mathcal{A}^{(k)}],\mathbb{G}[\mathcal{A}%
^{(l)}]\right)    & =\mathbb{V}\left(  \mathbb{G}[\mathcal{A}^{(k)}]\right)
\mathbf{P}_{\mathcal{A}^{(l)}|\mathcal{A}^{(k)}},\label{VP}\\
\Phi_{k}^{(N)}(f)  & =\mathbb{E}\left[  f\mid\mathcal{A}^{(k)}\right]
-\sum_{k<l\leqslant N}\mathbf{P}_{\mathcal{A}^{(l)}|\mathcal{A}^{(k)}}\cdot\Phi
_{l}^{(N)}(f).\label{recphikNf}%
\end{align}

\end{lem}

\noindent\textit{Proof.} The $j$-th coordinate of the vector $\mathbb{V}\left(
\mathbb{G}[\mathcal{A}^{(k)}]\right)  \cdot\mathbb{E}[f\mid\mathcal{A}^{(k)}]$
is%
\begin{align*}
&  P(A_{j}^{(k)})(1-P(A_{j}^{(k)}))\mathbb{E}(f\mid A_{j}^{(k)})-\sum_{j\neq
i\leqslant m_{k}}P(A_{i}^{(k)})P(A_{j}^{(k)})\mathbb{E}(f\mid A_{i}^{(k)})\\
&  =\mathbb{E}(1_{A_{j}^{(k)}}f)-P(A_{j}^{(k)})\sum_{1\leqslant i\leqslant
m_{k}}\mathbb{E}(1_{A_{i}^{(k)}}f)\\
&  =\mathrm{Cov}\left(  \mathbb{G}(A_{j}^{(k)}),\mathbb{G}(f)\right)  .
\end{align*}
Likewise the $(i,j)$-th coordinate of the matrix $\mathbb{V}\left(
\mathbb{G}[\mathcal{A}^{(k)}]\right)  \cdot\mathbf{P}_{\mathcal{A}%
^{(l)}|\mathcal{A}^{(k)}}$ is%
\begin{align*}
&  P(A_{i}^{(k)})(1-P(A_{i}^{(k)}))P(A_{j}^{(l)}\mid A_{i}^{(k)})-\sum_{j\neq
m\leqslant m_{k}}P(A_{i}^{(k)})P(A_{m}^{(k)})P(A_{j}^{(l)}\mid A_{m}^{(k)})\\
&  =P(A_{j}^{(l)}\cap A_{i}^{(k)})-P(A_{i}^{(k)})\sum_{1\leqslant m\leqslant
m_{k}}P(A_{j}^{(l)}\cap A_{m}^{(k)})\\
&  =\mathrm{Cov}\left(  \mathbb{G}(A_{i}^{(k)}),\mathbb{G}(A_{j}%
^{(l)})\right)  .
\end{align*}
By the definition (\ref{phikNf}) of the vectors $\Phi_{l}^{(N)}(f)$ we get%
\begin{align*}
&  \sum_{k<l\leqslant N}\mathbf{P}_{\mathcal{A}^{(l)}|\mathcal{A}^{(k)}}%
\cdot\Phi_{l}^{(N)}(f)\\
&  =\sum_{k<l\leqslant N}\mathbf{P}_{\mathcal{A}^{(l)}|\mathcal{A}^{(k)}%
}\cdot\mathbb{E}\left[  f\mid\mathcal{A}^{(l)}\right]  \\
&  +\sum_{\substack{k<l\leqslant N,1\leqslant L\leqslant N-l\\l<l_{1}%
<...<l_{L}\leqslant N}}(-1)^{L}\mathbf{P}_{\mathcal{A}^{(l)}|\mathcal{A}%
^{(k)}}\mathbf{P}_{\mathcal{A}^{(l_{1})}|\mathcal{A}^{(l)}}\dots
\mathbf{P}_{\mathcal{A}^{(l_{L})}|\mathcal{A}^{(l_{L-1})}}\cdot\mathbb{E}\left[
f\mid\mathcal{A}^{(l_{L})}\right]  \\
&  =\sum_{\substack{1\leqslant L\leqslant N-k\\k<l_{1}<...<l_{L}\leqslant
N}}(-1)^{L+1}\mathbf{P}_{\mathcal{A}^{(l_{1})}|\mathcal{A}^{(k)}}%
\dots\mathbf{P}_{\mathcal{A}^{(l_{L})}|\mathcal{A}^{(l_{L-1})}}\cdot\mathbb{E}%
\left[  f\mid\mathcal{A}^{(l_{L})}\right]  \\
&  =\mathbb{E}\left[  f\mid\mathcal{A}^{(k)}\right]  -\Phi_{k}^{(N)}(f),
\end{align*}
which yields (\ref{recphikNf}).$\quad\square$\medskip

\noindent\textbf{Step 3.} Let us first compute the variance of $\mathbb{G}^{(N)}(f)$. By Lemma \ref{Thm1} we have
\begin{align*}
&  \mathbb{V}\left(  \mathbb{G}^{(N)}(f)\right)  -\mathbb{V}\left(
\mathbb{G}(f)\right)  \\
&  =\sum_{k=1}^{N}\Phi_{k}^{(N)}(f)^{t}\cdot\mathbb{V}\left(  \mathbb{G}%
[\mathcal{A}^{(k)}]\right)  \cdot\Phi_{k}^{(N)}(f)-2\sum_{k=1}^{N}\Phi
_{k}^{(N)}(f)^{t}\cdot\mathrm{Cov}\left(  \mathbb{G}[\mathcal{A}%
^{(k)}],\mathbb{G}(f)\right)  \\
&  +2\sum_{1\leqslant k<l\leqslant N}\Phi_{k}^{(N)}(f)^{t}\cdot\mathrm{Cov}%
\left(  \mathbb{G}[\mathcal{A}^{(k)}],\mathbb{G}[\mathcal{A}^{(l)}]\right)
\cdot\Phi_{l}^{(N)}(f),
\end{align*}
hence Lemma \ref{covariance} gives, through (\ref{VE}) and (\ref{VP}),%
\[
\mathbb{V}\left(  \mathbb{G}^{(N)}(f)\right)  -\mathbb{V}\left(
\mathbb{G}(f)\right)  =\sum_{k=1}^{N}\Phi_{k}^{(N)}(f)^{t}\cdot\mathbb{V}%
\left(  \mathbb{G}[\mathcal{A}^{(k)}]\right)  \cdot\Psi_{k}^{(N)}(f),
\]
where, by (\ref{recphikNf}),%
\begin{equation}
\Psi_{k}^{(N)}(f)=\Phi_{k}^{(N)}(f)-2\mathbb{E}\left[  f\mid\mathcal{A}%
^{(k)}\right]  +2\sum_{k<l\leqslant N}\mathbf{P}_{\mathcal{A}^{(l)}%
|\mathcal{A}^{(k)}}\cdot\Phi_{l}^{(N)}(f)=-\Phi_{k}^{(N)}(f).\label{psi-phi}
\end{equation}
The formula (\ref{VGNf}) is proved. It extends to the
covariance since, by Lemma \ref{Thm1},%
\begin{align*}
&  \mathrm{Cov}(\mathbb{G}^{(N)}(f),\mathbb{G}^{(N)}(g))-\mathrm{Cov}%
(\mathbb{G}(f),\mathbb{G}(g))\\
&  =\frac{1}{2}\left(  \Upsilon_{N}(f,g)-2\Gamma_{N}(f,g)\right)  +\frac{1}%
{2}\left(  \Upsilon_{N}(g,f)-2\Gamma_{N}(g,f)\right),
\end{align*}
where, by (\ref{VE}) and (\ref{VP}) again,%
\begin{align*}
\Upsilon_{N}(f,g) &  =\sum_{1\leqslant k\leqslant l\leqslant N}\Phi_{k}^{(N)}(f)^{t}%
\cdot\mathbb{V}\left(  \mathbb{G}[\mathcal{A}^{(k)}]\right)
\mathbf{P}_{\mathcal{A}^{(l)}|\mathcal{A}^{(k)}}\cdot\Phi_{l}^{(N)}(g),\\
\Gamma_{N}(f,g) & =\sum_{l=1}^{N}\Phi_{k}^{(N)}(f)^{t}\cdot\mathbb{V}\left(
\mathbb{G}[\mathcal{A}^{(k)}]\right)\cdot\mathbb{E}
\left[g\mid\mathcal{A}^{(k)}\right].
\end{align*}
By replacing $\Psi_{k}^{(N)}(g)$ with $-\Phi_{k}^{(N)}(g)$
according to (\ref{psi-phi}), we obtain%
\[
\frac{1}{2}\left(  \Upsilon_{N}(f,g)-2\Gamma_{N}(f,g)\right)  =-\frac{1}%
{2}\sum_{k=1}^{N}\Phi_{k}^{(N)}(f)^{t}\cdot\mathbb{V}\left(  \mathbb{G}%
[\mathcal{A}^{(k)}]\right)  \cdot\Phi_{k}^{(N)}(g),
\]
which is symmetric in $f$ and $g$. The covariance formula of
Proposition~\ref{Thm2} is proved.

\subsection{Proof of Propositions~\ref{Thm2b} and \ref{Thm3}}

Since $\mathbb{V}\left(  \mathbb{G}[\mathcal{A}^{(k)}]\right)  $ is
semi-definite positive, for all $1\leqslant k\leqslant N$ and $f\in
\mathcal{F}$ we have
\[
\Phi_{k}^{(N)}(f)^{t}\cdot\mathbb{V}(\mathbb{G}%
[\mathcal{A}^{(k)}])\cdot\Phi_{k}^{(N)}(f)\geqslant0,
\]
and the variance part
(\ref{VGNf}) of Proposition \ref{Thm2b} follows from Proposition
\ref{Thm2}. For any $m\in\mathbb{N}_{\ast}$, $(f_{1},...,f_{m})\in
\mathcal{F}^{m}$ and $u\in\mathbb{R}^{m}$, it further holds, by Proposition \ref{Thm2}
again,%
\begin{align*}
& u^{t}\left(  \Sigma_{m}^{(0)}-\Sigma_{m}^{(N)}\right)  u \\
& =\sum_{1\leqslant
i,j\leqslant m}u_{i}u_{j}\left(  \mathrm{Cov}(\mathbb{G}(f_{i}),\mathbb{G}%
(f_{j}))-\mathrm{Cov}(\mathbb{G}^{(N)}(f_{i}),\mathbb{G}^{(N)}(f_{j}))\right)
\\
& =\sum_{k=1}^{N}\sum_{1\leqslant i,j\leqslant m}\left(  u_{i}\Phi_{k}%
^{(N)}(f_{i})\right)  ^{t}\cdot\mathbb{V}\left(  \mathbb{G}[\mathcal{A}%
^{(k)}]\right)  \cdot\left(  u_{j}\Phi_{k}^{(N)}(f_{j})\right)  \\
& =\sum_{k=1}^{N}\left(  \sum_{1\leqslant i\leqslant m}u_{i}\Phi_{k}%
^{(N)}(f_{i})\right)  ^{t}\cdot\mathbb{V}\left(  \mathbb{G}[\mathcal{A}%
^{(k)}]\right)  \cdot\left(  \sum_{1\leqslant j\leqslant m}u_{j}\Phi_{k}%
^{(N)}(f_{j})\right)  \geqslant0.
\end{align*}
Under the wrapping hypothesis of Proposition~\ref{Thm3} we have%
\[
\Phi_{N_{0}-k}^{(N_{0})}=\Phi_{N_{1}-k}^{(N_{1})},\quad0\leqslant k<N_{0},
\]
since the corresponding (\ref{phikNf}) only involves 
$\mathcal{A}^{(N_{0}-k)}=\mathcal{A}^{(N_{1}-k)}$ for
$0\leqslant k<N_{0}$.
Assuming moreover $N_{1}\geqslant2N_{0}$ we get, by Proposition~\ref{Thm2},
\begin{align*}
&  \mathbb{V}\left(  \mathbb{G}^{(N_{0})}(f)\right)  -\mathbb{V}\left(
\mathbb{G}^{(N_{1})}(f)\right)  \\
&  =\sum_{k=1}^{N_{1}}\Phi_{k}^{(N_{1})}(f)^{t}\cdot\mathbb{V}\left(
\mathbb{G}[\mathcal{A}^{(k)}]\right)  \cdot\Phi_{k}^{(N_{1})}(f)-\sum
_{k=1}^{N_{0}}\Phi_{k}^{(N_{0})}(f)^{t}\cdot\mathbb{V}\left(  \mathbb{G}%
[\mathcal{A}^{(k)}]\right)  \cdot\Phi_{k}^{(N_{0})}(f)\\
&  =\sum_{k=N_{0}+1}^{N_{1}-N_{0}}\Phi_{k}^{(N_{1})}(f)^{t}\cdot
\mathbb{V}\left(  \mathbb{G}[\mathcal{A}^{(k)}]\right)  \cdot\Phi_{k}%
^{(N_{1})}(f)\geqslant0,
\end{align*}
and, for $m\in\mathbb{N}_{\ast}$, $(f_{1},...,f_{m})\in\mathcal{F}^{m}$ and
$u\in\mathbb{R}^{m}$,%
\begin{align*}
&  u^{t} \cdot\left(  \Sigma_{m}^{(N_{0})}-\Sigma_{m}^{(N_{1})}\right) \cdot u\\
&  =u^{t}\cdot \left(  (\Sigma_{m}^{(0)}-\Sigma_{m}^{(N_{1})})-(\Sigma_{m}%
^{(0)}-\Sigma_{m}^{(N_{0})})\right) \cdot u\\
&  =\sum_{k=N_{0}+1}^{N_{1}-N_{0}}\left(  \sum_{1\leqslant i\leqslant m}%
u_{i}\Phi_{k}^{(N)}(f_{i})\right)  ^{t}\cdot\mathbb{V}\left(  \mathbb{G}%
[\mathcal{A}^{(k)}]\right)  \cdot\left(  \sum_{1\leqslant j\leqslant m}%
u_{j}\Phi_{k}^{(N)}(f_{j})\right)  \geqslant0.
\end{align*}

\subsection{Proof of Proposition~\ref{Pro3}}

We show the result by double induction. For $m=0$ we have $\mathbb{G}%
^{(0)}(f)=\mathbb{G}(f)$ and, by (\ref{GNf}), $\mathbb{G}^{(1)}(f)=\mathbb{G}%
(f)-\mathbb{E}[f|\mathcal{A}]^{t}\mathbb{G}[\mathcal{A}]$. Assume that
(\ref{G2m}) and (\ref{G2m1}) are true for $m\in\mathbb{N}$. For $m+1$ we have,
by the raking ratio transform (\ref{GNf}),%
\begin{align}
\mathbb{G}^{(2m+2)}(f) &  =\mathbb{G}^{(2m+1)}(f)-\mathbb{E}[f|\mathcal{B}%
]^{t}\cdot\mathbb{G}^{(2m+1)}[\mathcal{B}],\label{G2m2}\\
\mathbb{G}^{(2m+3)}(f) &  =\mathbb{G}^{(2m+2)}(f)-\mathbb{E}[f|\mathcal{A}%
]^{t}\cdot\mathbb{G}^{(2m+2)}[\mathcal{A}].\label{G2m3}%
\end{align}
For $1\leqslant j\leqslant m_{2}$ and $f=1_{B_{j}}$ we get, by (\ref{G2m1}), 
(\ref{RakingRatioS1p}), (\ref{RakingRatioS2i}), (\ref{RakingRatioS1i})
and (\ref{RakingRatioS2p}),%
\begin{align*}
& \mathbb{G}^{(2m+1)}(B_{j}) \\ 
& =\mathbb{G}(B_{j})-S_{1,\text{odd}}^{(m-1)}(1_{B_{j}})^{t}\cdot\mathbb{G}[\mathcal{A}]
-S_{2,\text{odd}}^{(m-1)}(1_{B_{j}})^{t}\cdot\mathbb{G}[\mathcal{B}] \\ 
&=\mathbb{G}(B_{j})-\left(  \sum_{k=0}^{m-1}\left(  \mathbf{P}_{\mathcal{A}|\mathcal{B}%
}\mathbf{P}_{\mathcal{B}|\mathcal{A}}\right)  ^{k}\cdot\left(  \mathbb{E}%
[1_{B_{j}}|\mathcal{B}]-\mathbf{P}_{\mathcal{A}|\mathcal{B}}\cdot\mathbb{E}%
[1_{B_{j}}|\mathcal{A}]\right)  \right)^{t}\cdot\mathbb{G}[\mathcal{B}]
\\
&-\left(  \sum_{k=0}^{m-1}\left(  \mathbf{P}_{\mathcal{B}%
|\mathcal{A}}\mathbf{P}_{\mathcal{A}|\mathcal{B}}\right)  ^{k}\cdot\left(
\mathbb{E}[1_{B_{j}}|\mathcal{A}]-\mathbf{P}_{\mathcal{B}|\mathcal{A}%
}\cdot\mathbb{E}[1_{B_{j}}|\mathcal{B}]\right)  +\left(  \mathbf{P}_{\mathcal{B}%
|\mathcal{A}}\mathbf{P}_{\mathcal{A}|\mathcal{B}}\right)^{m}\cdot\mathbb{E}%
[1_{B_{j}}|\mathcal{A}]\right)^{t}\cdot\mathbb{G}[\mathcal{A}],
\end{align*}
where $\mathbb{E}[1_{B_{j}}|\mathcal{A}]$ is the $j$-th column of
$\mathbf{P}_{\mathcal{B}|\mathcal{A}}$ and $\mathbb{E}[1_{B_{j}}|\mathcal{B}]$
is the $j$-th unit vector of $\mathbb{R}^{m_{2}}$. Therefore%
\begin{align*}
& \mathbb{G}^{(2m+1)}[\mathcal{B}] \\ 
& =\mathbb{G}[\mathcal{B}]-\left(\sum_{k=0}^{m-1}(\mathbf{P}_{\mathcal{A}|\mathcal{B}%
}\mathbf{P}_{\mathcal{B}|\mathcal{A}})^{k}(\mathrm{Id}_{m_{2}}-\mathbf{P}%
_{\mathcal{A}|\mathcal{B}}\mathbf{P}_{\mathcal{B}|\mathcal{A}})\right)
^{t}\cdot\mathbb{G}[\mathcal{B}] \\
& \quad -\left(\sum_{k=0}^{m-1}(\mathbf{P}_{\mathcal{B}|\mathcal{A}}\mathbf{P}_{\mathcal{A}%
|\mathcal{B}})^{k}(\mathbf{P}_{\mathcal{B}|\mathcal{A}}-\mathbf{P}%
_{\mathcal{B}|\mathcal{A}}\mathrm{Id}_{m_{2}})+(\mathbf{P}_{\mathcal{B}%
|\mathcal{A}}\mathbf{P}_{\mathcal{A}|\mathcal{B}})^{m}\mathbf{P}%
_{\mathcal{B}|\mathcal{A}}\right)^{t}\cdot\mathbb{G}[\mathcal{A}] \\
& =\left( (\mathbf{P}_{\mathcal{A}|\mathcal{B}%
}\mathbf{P}_{\mathcal{B}|\mathcal{A}})^{m}\right) ^{t}\cdot\mathbb{G}%
[\mathcal{B}] -\left( (\mathbf{P}_{\mathcal{B}|\mathcal{A}}\mathbf{P}_{\mathcal{A}%
|\mathcal{B}})^{m}\mathbf{P}_{\mathcal{B}|\mathcal{A}}\right) ^{t}%
\cdot\mathbb{G}[\mathcal{A}],
\end{align*}
Finally (\ref{G2m}) and (\ref{G2m2}) then again (\ref{RakingRatioS1p}),
(\ref{RakingRatioS2i}), (\ref{RakingRatioS1i}) and (\ref{RakingRatioS2p})
together imply%
\begin{align*}
& \mathbb{G}^{(2m+2)}(f) \\ 
& =\mathbb{G}(f)- S_{1,\text{odd}}^{(m-1)}(f) ^{t}\cdot\mathbb{G}[\mathcal{A}]
- S_{2,\text{odd}}^{(m-1)}(f) ^{t}\cdot\mathbb{G}[\mathcal{B}]
-\mathbb{E}[f|\mathcal{B}]^{t}\cdot\mathbb{G}^{(2m+1)}[\mathcal{B}]\\
& =\mathbb{G}(f)-\left( S_{1,\text{odd}}^{(m-1)}(f)-(\mathbf{P}%
_{\mathcal{B}|\mathcal{A}}\mathbf{P}_{\mathcal{A}|\mathcal{B}})^{m}%
\mathbf{P}_{\mathcal{B}|\mathcal{A}}\cdot\mathbb{E}[f|\mathcal{B}]\right)^{t}
\cdot\mathbb{G}[\mathcal{A}]\\
& \quad-\left(  S_{2,\text{odd}}^{(m-1)}(f)+(\mathbf{P}_{\mathcal{A}|\mathcal{B}%
}\mathbf{P}_{\mathcal{B}|\mathcal{A}})^{m}\cdot\mathbb{E}[f|\mathcal{B}]\right)^{t}
\cdot\mathbb{G}[\mathcal{B}]\\
& =\mathbb{G}(f)-S_{1,\text{even}}^{(m)}(f) ^{t}\cdot\mathbb{G}[\mathcal{A}]
- S_{2,\text{even}}^{(m-1)}(f) ^{t}\cdot\mathbb{G}[\mathcal{B}],
\end{align*}
and (\ref{G2m}) is valid for $m+1$. If $1\leqslant i\leqslant m_{1}$ then
$\mathbb{E}[1_{A_{i}}|\mathcal{B}]$ is the $i$-th column of $\mathbf{P}%
_{\mathcal{A}|\mathcal{B}}$ and $\mathbb{E}[1_{A_{i}}|\mathcal{A}]$ is the
$i$-th unit vector of $\mathbb{R}^{m_{1}}$ thus (\ref{G2m}) for $m+1$ and
$f=1_{B_{i}}$ in turn entails%
\begin{align*}
& \mathbb{G}^{(2m+2)}[\mathcal{A}] \\ 
& =\mathbb{G}[\mathcal{A}]-\left(\sum_{k=0}^{m}\left( \mathbf{P}_{\mathcal{B}|\mathcal{A}}
\mathbf{P}_{\mathcal{A}|\mathcal{B}}\right) ^{k}\left(  \mathrm{Id}_{m_{1}}
-\mathbf{P}_{\mathcal{B}|\mathcal{A}}\mathbf{P}_{\mathcal{A}|\mathcal{B}%
}\right) \right) ^{t}\cdot\mathbb{G}[\mathcal{A}]\\
& \quad -\left( \sum_{k=0}^{m-1}\left(  \mathbf{P}_{\mathcal{A}|\mathcal{B}%
}\mathbf{P}_{\mathcal{B}|\mathcal{A}}\right)  ^{k}\left(  \mathbf{P}%
_{\mathcal{A}|\mathcal{B}}-\mathbf{P}_{\mathcal{A}|\mathcal{B}}\mathrm{Id}_{m_{1}}\right)
+\left(  \mathbf{P}_{\mathcal{A}|\mathcal{B}}\mathbf{P}%
_{\mathcal{B}|\mathcal{A}}\right)  ^{m}\mathbf{P}_{\mathcal{A}|\mathcal{B}%
}\right) ^{t}\cdot\mathbb{G}[\mathcal{B}]\\
& =\left( (\mathbf{P}_{\mathcal{B}|\mathcal{A}}\mathbf{P}_{\mathcal{A}%
|\mathcal{B}})^{m+1}\right) ^{t}\cdot\mathbb{G}[\mathcal{A}]
-\left( (\mathbf{P}_{\mathcal{A}|\mathcal{B}}\mathbf{P}_{\mathcal{B}|\mathcal{A}})  ^{m}\mathbf{P}_{\mathcal{A}|\mathcal{B}}\right) ^{t}\cdot\mathbb{G}%
[\mathcal{B}],
\end{align*}
and also, thanks to (\ref{G2m1}) and (\ref{G2m3}),%
\begin{align*}
& \mathbb{G}^{(2m+3)}(f) \\ 
& =\mathbb{G}(f)-\left( S_{1,\text{odd}}^{(m)}(f)- (\mathbf{P}_{\mathcal{B}|\mathcal{A}}
\mathbf{P}_{\mathcal{A}|\mathcal{B}})^{m+1}\cdot\mathbb{E}[f|\mathcal{A}]\right)
^{t}\cdot\mathbb{G}[\mathcal{A}] \\
& \quad -\left( S_{2,\text{odd}}^{(m-1)}(f)+\left(  \mathbf{P}_{\mathcal{A}%
|\mathcal{B}}\mathbf{P}_{\mathcal{B}|\mathcal{A}}\right)  ^{m}\cdot\mathbb{E}%
[f|\mathcal{B}]\right)^{t}\cdot\mathbb{G}[\mathcal{B}]-\mathbb{E}[f|\mathcal{A}%
]^{t}\cdot\mathbb{G}^{(2m+2)}[\mathcal{A}]\\
& =\mathbb{G}(f)- S_{1,\text{odd}}^{(m)}(f) ^{t}\cdot\mathbb{G}[\mathcal{A}]
- S_{2,\text{odd}}^{(m)}(f) ^{t}\cdot\mathbb{G}[\mathcal{B}],
\end{align*}
which is (\ref{G2m1}) for $m+1$.%

\subsection{Proof of Proposition~\ref{Pro4}}

\noindent\textbf{Step 1.} For $m\geqslant1$ let $\mathbf{0}_{m,m}$ be the $m\times m$ 
null matrix. Also recall the vectors $P(\mathcal{A})=(P(A_{1}),\dots,P(A_{m_{1}}))$ and 
$P(\mathcal{B})=(P(B_{1}),\dots,P(B_{m_{2}}))$.

\begin{lem}
\label{LemUlAl}Assume (ER). For $l=1,2$ there exists an invertible $m_{l}\times m_{l}$
matrix $U_{l}$ and an upper triangular $(m_{l}-1)\times(m_{l}-1)$
matrix $T_{l}$ such that%
\begin{align*}
\mathbf{P}_{\mathcal{B}|\mathcal{A}}\mathbf{P}_{\mathcal{A}|\mathcal{B}}  &
=U_{1}\left(
\begin{matrix}
1 & \mathbf{0}_{1,m_{1}-1}\\
\mathbf{0}_{m_{1}-1,1} & T_{1}%
\end{matrix}
\right)  U_{1}^{-1},\quad\lim_{k\rightarrow+\infty}T_{1}^{k}=\mathbf{0}%
_{m_{1}-1,m_{1}-1},\\
\mathbf{P}_{\mathcal{A}|\mathcal{B}}\mathbf{P}_{\mathcal{B}|\mathcal{A}}  &
=U_{2}\left(
\begin{matrix}
1 & \mathbf{0}_{1,m_{2}-1}\\
\mathbf{0}_{m_{2}-1,1} & T_{2}%
\end{matrix}
\right)  U_{2}^{-1},\quad\lim_{k\rightarrow+\infty}T_{2}^{k}=\mathbf{0}%
_{m_{2}-1,m_{2}-1}.
\end{align*}

\end{lem}

\noindent\textit{Proof.} Since $\mathcal{A}$ is a partition, for $1\leqslant 
i\leqslant m_{2}$ the sum of the
$m_{1}$ terms of row $i$ of $\mathbf{P}_{\mathcal{A}|\mathcal{B}}$ is
$\sum_{j=1}^{m_{1}}P(A_{j}\mid B_{i})=1$
hence $\mathbf{P}_{\mathcal{A}|\mathcal{B}}$ is stochastic. Likewise
$\mathbf{P}_{\mathcal{B}|\mathcal{A}}$ is stochastic and, by stability, so are
$\mathbf{P}_{\mathcal{A}|\mathcal{B}}\mathbf{P}_{\mathcal{B}|\mathcal{A}}$ and
$\mathbf{P}_{\mathcal{B}|\mathcal{A}}\mathbf{P}_{\mathcal{A}|\mathcal{B}}$.
Let the column of $1$'s associated to their eigenvalue $1$ be in first position 
in their respective matrix $U_{1},U_{2}$ of eigenvectors. 
The announced decomposition is always true with some 
upper triangular matrices $T_{l}$ having Jordan decomposition $T_{l}=D_{l}+N_{l}$
where $D_{l}=Q_{l}\Delta_{l}Q_{l}^{-1}$, $\Delta_{l}$
is a diagonal $(m_{l}-1)\times(m_{l}-1)$ matrix, $Q_{l}$ is an
invertible $(m_{l}-1)\times(m_{l}-1)$ matrix and $N_{l}$ is a
nilpotent $(m_{l}-1)\times(m_{l}-1)$ matrix of order $n_{l}\geqslant1$
that commute with $D_{l}$. Next observe that
\begin{align*}
(P(\mathcal{A})\cdot\mathbf{P}_{\mathcal{B}|\mathcal{A}}\mathbf{P}%
_{\mathcal{A}|\mathcal{B}})_{k} &  =\sum_{i=1}^{m_{1}}P(A_{i})(\mathbf{P}%
_{\mathcal{B}|\mathcal{A}}\mathbf{P}_{\mathcal{A}|\mathcal{B}})_{i,k} \\
& =\sum_{i=1}^{m_{1}}P(A_{i}) \sum_{j=1}^{m_{2}}(\mathbf{P}%
_{\mathcal{B}|\mathcal{A}})_{i,j}(\mathbf{P}_{\mathcal{A}|\mathcal{B}})_{j,k} \\
&  =\sum_{j=1}^{m_{2}}\sum_{i=1}^{m_{1}}P(A_{i})P(B_{j}|A_{i})P(A_{k}%
|B_{j})=P(A_{k}),
\end{align*}
which proves that $\mathbf{P}_{\mathcal{B}|\mathcal{A}}\mathbf{P}%
_{\mathcal{A}|\mathcal{B}}$ has invariant probability $P(\mathcal{A})$.
Similarly, $P(\mathcal{B})$ is invariant for $\mathbf{P}_{\mathcal{A}%
|\mathcal{B}}\mathbf{P}_{\mathcal{B}|\mathcal{A}}$, and the first line of $U_{1}^{-1}$ and $U_{2}^{-1}$ is $P[\mathcal{A}]$ and $P[\mathcal{B}]$ respectively. 
Under (ER) these matrices are ergodic, which ensures that the eigenvalues of $\Delta_{l}$ 
have moduli strictly less than the dominant $1$ since it is the case of 
eigenvalues of $T_{l}$ hence $D_{l}$. It follows that%
\[
\lim_{k\rightarrow+\infty}\Delta_{l}^{k}=\mathbf{0}_{m_{l}-1,m_{l}-1},\quad
l=1,2.
\]
Furthermore, since $N_{l}$ and $D_{l}$ commute it holds
\[
T_{l}^{k}=\sum_{j=0}^{n_{l}-1}\binom{k}{j}N_{l}^{j}D_{l}^{k-j},\quad
l=1,2,\quad k\geqslant n_{l}.
\]
We conclude that $\lim_{k\rightarrow+\infty}T_{l}^{k}=\mathbf{0}%
_{m_{l}-1,m_{l}-1}$.$\quad\square$\medskip

\noindent \textbf{Step 2.} Let $V_{1}(f)=(\mathbb{E}[f|\mathcal{A}]
-\mathbf{P}_{\mathcal{B}|\mathcal{A}} \cdot \mathbb{E}[f|\mathcal{B}])$ and
$V_{2}(f)=(\mathbb{E}[f|\mathcal{B}]-\mathbf{P}_{\mathcal{A}|\mathcal{B}}\cdot \mathbb{E}[f|\mathcal{A}]).$

\begin{lem}
\label{Lem3}Under (ER) we have
\[
\lim_{k\rightarrow+\infty}(\mathbf{P}_{\mathcal{B}|\mathcal{A}}\mathbf{P}%
_{\mathcal{A}|\mathcal{B}})^{k}\cdot V_{1}(f)=\mathbf{0}_{m_{1},1},\;
\lim_{k\rightarrow+\infty}(\mathbf{P}_{\mathcal{A}|\mathcal{B}}\mathbf{P}%
_{\mathcal{B}|\mathcal{A}})^{k}\cdot V_{2}(f)=\mathbf{0}_{m_{2},1}.
\]

\end{lem}

\noindent\textit{Proof.} By Lemma~\ref{LemUlAl} we have%
\begin{equation}
\lim_{k\rightarrow+\infty}(\mathbf{P}_{\mathcal{B}|\mathcal{A}}\mathbf{P}%
_{\mathcal{A}|\mathcal{B}})^{k}=\left(
\begin{smallmatrix}
P(\mathcal{A})\\
\vdots\\P(\mathcal{A})\end{smallmatrix}\right), \;
\lim_{k\rightarrow+\infty}(\mathbf{P}_{\mathcal{A}|\mathcal{B}}
\mathbf{P}_{\mathcal{B}|\mathcal{A}})^{k}=\left(
\begin{smallmatrix}
P(\mathcal{B})\\
\vdots\\
P(\mathcal{B})
\end{smallmatrix}
\right).\label{PAPB}%
\end{equation}
The scalar product of $P(\mathcal{A})$ by $V_{1}(f)$ is null since
$P(\mathcal{A})\cdot\mathbb{E}[f|\mathcal{A}]=P\mathbb{(}f)$ and
\begin{align*}
P(\mathcal{A})\cdot\mathbf{P}_{\mathcal{B}|\mathcal{A}}\mathbb{E}%
[f|\mathcal{B}]) &=\sum_{j=1}^{m_{1}}P(A_{j})\sum_{k=1}^{m_{2}}P(B_{k}%
\;|\;A_{j})\mathbb{E}(f|B_{k})\\
& =\sum_{j=1}^{m_{1}}\sum_{k=1}^{m_{2}}P(A_{j}\cap
B_{k})\mathbb{E}[f|B_{k}]=P(f).
\end{align*}
Likewise we get $P(\mathcal{B})\cdot V_{2}(f)=0.\quad\square$\medskip

\noindent The following convergences hold for any matrix norm. By Lemma~\ref{LemUlAl} we have%
\begin{align*}
\sum_{k=0}^{N}(\mathbf{P}_{\mathcal{B}|\mathcal{A}}\mathbf{P}_{\mathcal{A}%
|\mathcal{B}})^{k} &  =U_{1}\left(
\begin{matrix}
N+1 & \mathbf{0}_{1,m_{1}-1}\\
\mathbf{0}_{m_{1}-1,1} & \sum_{k=0}^{N}T_{1}^{k}%
\end{matrix}
\right)  U_{1}^{-1},\\
\sum_{k=0}^{N}(\mathbf{P}_{\mathcal{A}|\mathcal{B}}\mathbf{P}_{\mathcal{B}%
|\mathcal{A}})^{k} &  =U_{2}\left(
\begin{matrix}
N+1 & \mathbf{0}_{1,m_{2}-1}\\
\mathbf{0}_{m_{2}-1,1} & \sum_{k=0}^{N}T_{2}^{k}%
\end{matrix}
\right)  U_{2}^{-1}.
\end{align*}
Now, the matrices $\mathrm{Id}_{m_{1}-1}-T_{1}$ and $\mathrm{Id}_{m_{2}%
-1}-T_{2}$ are nonsingular since $1$ is a dominant eigenvalue of
$\mathbf{P}_{\mathcal{B}|\mathcal{A}}\mathbf{P}_{\mathcal{A}|\mathcal{B}}$ and
$\mathbf{P}_{\mathcal{A}|\mathcal{B}}\mathbf{P}_{\mathcal{B}|\mathcal{A}}%
$.\ Recalling~(\ref{RakingRatioS1p}) and (\ref{RakingRatioS2i}), by Lemma~\ref{LemUlAl} and~\ref{Lem3} it follows that
\begin{align*}
S_{1,\text{even}}^{(N)}(f) &  =U_{1}\left(
\begin{matrix}
0 & \mathbf{0}_{1,m_{1}-1}\\
\mathbf{0}_{m_{1}-1,1} & (\mathrm{Id}_{m_{1}-1}-T_{1})^{-1}(\mathrm{Id}%
_{m_{1}-1}-T_{1}^{N+1})
\end{matrix}
\right)  U_{1}^{-1}\cdot V_{1}(f),\\
S_{2,\text{odd}}^{(N)}(f) &  =U_{2}\left(
\begin{matrix}
0 & \mathbf{0}_{1,m_{2}-1}\\
\mathbf{0}_{m_{2}-1,1} & (\mathrm{Id}_{m_{2}-1}-T_{2})^{-1}(\mathrm{Id}%
_{m_{2}-1}-T_{2}^{N+1})
\end{matrix}
\right)  U_{2}^{-1}\cdot V_{2}(f),
\end{align*}
which, by Lemma \ref{LemUlAl}, converge respectively to
\begin{align*}
S_{1,\text{even}}(f) &  =U_{1}\left(
\begin{matrix}
0 & \mathbf{0}_{1,m_{1}-1}\\
\mathbf{0}_{m_{1}-1,1} & (\mathrm{Id}_{m_{1}-1}-T_{1})^{-1}%
\end{matrix}
\right)  U_{1}^{-1}\cdot V_{1}(f),\\
S_{2,\text{odd}}(f) &  =U_{2}\left(
\begin{matrix}
0 & \mathbf{0}_{1,m_{2}-1}\\
\mathbf{0}_{m_{2}-1,1} & (\mathrm{Id}_{m_{2}-1}-T_{2})^{-1}%
\end{matrix}
\right)  U_{2}^{-1}\cdot V_{2}(f).
\end{align*}
Since we have already seen by using (\ref{PAPB}) and the notations of Proposition~\ref{Pro4} that%
\[
\lim_{k\rightarrow+\infty}(\mathbf{P}_{\mathcal{B}|\mathcal{A}}\mathbf{P}%
_{\mathcal{A}|\mathcal{B}})^{k}\cdot\mathbb{E}[f|\mathcal{A}]=P_{1}[f],\text{
}\lim_{k\rightarrow+\infty}(\mathbf{P}_{\mathcal{A}|\mathcal{B}}%
\mathbf{P}_{\mathcal{B}|\mathcal{A}})^{k}\cdot\mathbb{E}[f|\mathcal{B}%
]=P_{2}[f],
\]
we conclude by (\ref{RakingRatioS1i}) and (\ref{RakingRatioS2p}) that
$S_{1,\text{odd}}^{(N)}(f)$, $S_{2,\text{even}}^{(N)}(f)$ converge to
the vectors $S_{1,\text{odd}}(f)=S_{1,\text{even}}(f)+P_{1}[f]$,
$S_{2,\text{even}}(f)=S_{1,\text{odd}}(f)+P_{1}[f]$ respectively.\medskip

\noindent \textbf{Step 3.} Given the spectral radius $\rho(T_{l})<1$ of $T_{l}$ let $\lambda
_{l}=\rho(T_{l})+\varepsilon<1$, $l=1,2$ for any $ \varepsilon>0$. Then there exists a vector norm
$\left\Vert \cdot\right\Vert _{l}$ on $\mathbb{C}^{m_{l}-1}$ such that its
induced matrix norm $\left\vert \left\Vert \cdot\right\Vert \right\vert _{l}$
on matrices $(m_{l}-1)\times(m_{l}-1)$ satisfies $\left\vert \left\Vert
T_{l}\right\Vert \right\vert _{l}\leqslant\lambda_{l}$. Introduce the vector
norm $\left\Vert (x_{1},...,x_{m_{l}})^{t}\right\Vert _{l}^{\prime
}=\left\vert x_{1}\right\vert +\left\Vert (x_{2},...,x_{m_{l}})^{t}\right\Vert
_{l}$ on $\mathbb{C}^{m_{l}}$ and the induced operator norm $\left\vert \left\Vert \cdot\right\Vert
\right\vert _{l}^{\prime}$ for $m_{l}\times m_{l}$ matrices. Then we have%
\[
\left\vert \left\Vert \left(
\begin{array}
[c]{cc}%
0 & \mathbf{0}_{1,m_{l}-1}\\
\mathbf{0}_{m_{l}-1,1} & T
\end{array}
\right)  \right\Vert \right\vert _{l}^{\prime}=\sup\left\{  x \in \mathbb{C}^{m_{l}}:\frac
{0+\left\Vert T(x_{2},...,x_{m_{l}})^{t}\right\Vert _{l}}{\left\vert
x_{1}\right\vert +\left\Vert (x_{2},...,x_{m_{l}})^{t}\right\Vert _{l}%
}\right\}  =\left\vert \left\Vert T\right\Vert \right\vert _{l},%
\]
for any $m_{l}\times m_{l}$ matrix $T$. Let $K_{l}=\left\vert \left\Vert (\mathrm{Id}%
_{m_{l}-1}-T_{l})^{-1}\right\Vert \right\vert _{l} $, $ \widetilde{K}_l = |||U_l|||_l' \ |||U_l^{-1}|||_l'$ and $ K_l^{\prime}>0$ be such that
$\left\Vert \cdot\right\Vert _{l}^{\prime}\leqslant K_{l}^{\prime}\left\Vert
\cdot\right\Vert _{\infty}$.\ By using Lemmas~\ref{LemUlAl} and \ref{Lem3} we get%
\begin{align*}
&\left\vert \left\Vert S_{l,even}^{(N)}(f)-S_{l,even}(f)\right\Vert \right\vert
_{l}^{\prime} \\
&\leqslant\left\vert \left\Vert U_{l}\left(
\begin{array}
[c]{cc}%
0 & \mathbf{0}_{1,m_{l}-1}\\
\mathbf{0}_{m_{l}-1,1} & -(\mathrm{Id}_{m_{l}-1}-T_{l})^{-1}T_{l}^{N+1}%
\end{array}
\right)  U_{l}^{-1}\right\Vert \right\vert _{l}^{\prime}\left\Vert
V_{l}(f)\right\Vert _{l}^{\prime}\\
& \leqslant \widetilde{K}_l \left\vert \left\Vert (\mathrm{Id}_{m_{l}-1}-T_{l})^{-1}T_{l}%
^{N+1}\right\Vert \right\vert _{l}^{\prime
}\left\Vert V_{l}(f)\right\Vert _{l}^{\prime
}\\
& \leqslant \widetilde{K}_l K_{l}\lambda_{l}^{N+1}K_{l}^{\prime}\left\Vert V_{l}(f)\right\Vert
_{\infty} \leqslant c_{l}\lambda_{l}^{N+1},%
\end{align*}
where $c_{l}=\widetilde{K}_l K_{l}K_{l}^{\prime}M$. Similar constants show up for $ \left\vert \left\Vert S_{l,odd}^{(N)}(f)-S_{l,odd}(f)\right\Vert \right\vert
_{l}^{\prime} $. The final constants $c_1 $ and $c_2$ depend on $ \lambda_1,\lambda_2,\varepsilon$, both matrices (ER) but also the two implicit constants relating the norms $||\cdot||_{m_1},||\cdot||_{m_2}$ of Proposition~\ref{Pro4} to the equivalent norms $ ||\cdot||_1',||\cdot||_2' $ .

\subsection{Proof of Theorem~\ref{Thm4}}

Write $Z_1=\max_{1\leqslant j\leqslant m_{1}}\left\vert \mathbb{G}(A_{j})\right\vert$ 
and $Z_2=\max_{1\leqslant j\leqslant m_{2}}\left\vert \mathbb{G}(B_{j})\right\vert$. 
According to Proposition~\ref{Pro4} the random variables appearing on the
right hand side of the following formulae, for $\ast\in\left\{  \text{even,odd}\right\}  $,%
\begin{align*}
\sup_{f\in\mathcal{F}}\left\vert (S_{1,\ast}^{(N)}(f)-S_{1,\ast}(f))^{t}%
\cdot\mathbb{G}[\mathcal{A}]\right\vert  & \leqslant c_{1}\sup_{f\in
\mathcal{F}}\left\Vert S_{1,\ast}^{(N)}(f)-S_{1,\ast}(f)\right\Vert _{m_{1}}Z_1, \\
\sup_{f\in\mathcal{F}}\left\vert (S_{2,\ast}^{(N)}(f)-S_{2,\ast}(f))^{t}%
\cdot\mathbb{G}[\mathcal{B}]\right\vert  & \leqslant c_{2}\sup_{f\in
\mathcal{F}}\left\Vert S_{2,\ast}^{(N)}(f)-S_{2,\ast}(f)\right\Vert _{m_{2}}Z_2,
\end{align*}
almost surely converge to $0$ since $\mathbb{P}\left(  ||\mathbb{G}%
||_{\mathcal{F}}<+\infty\right)  =1$. Hence the processes $\mathbb{G}%
^{(2m)},\mathbb{G}^{(2m+1)}$ converge almost surely in $\ell^{\infty
}(\mathcal{F})$ to $\mathbb{G}_{\text{even}}^{(\infty)},\mathbb{G}%
_{\text{odd}}^{(\infty)}$ defined by
\begin{align*}
\mathbb{G}_{\text{even}}^{(\infty)}(f) &  =\mathbb{G}(f)-S_{1,\text{even}%
}(f)^{t}\cdot\mathbb{G}[\mathcal{A}]-S_{2,\text{even}}(f)^{t}\cdot
\mathbb{G}[\mathcal{B}]\\
&  =\mathbb{G}^{(\infty)}(f)-P_{2}[f]^{t}\cdot\mathbb{G}[\mathcal{B}],\\
\mathbb{G}_{\text{odd}}^{(\infty)}(f) &  =\mathbb{G}(f)-S_{1,\text{odd}%
}(f)^{t}\cdot\mathbb{G}[\mathcal{A}]-S_{2,\text{odd}}(f)^{t}\cdot
\mathbb{G}[\mathcal{B}]\\
&  =\mathbb{G}^{(\infty)}(f)-P_{1}[f]^{t}\cdot\mathbb{G}[\mathcal{A}],
\end{align*}
with $\mathbb{G}^{(\infty)}(f)=\mathbb{G}(f)-S_{1,\text{even}}(f)^{t}%
\cdot\mathbb{G}[\mathcal{A}]-S_{2,\text{odd}}(f)^{t}\cdot\mathbb{G}%
[\mathcal{B}]$ and using (\ref{RakingRatioS2p}), (\ref{RakingRatioS1i}).\ Since $P_{1}[f]^{t}\cdot\mathbb{G}%
[\mathcal{A}]=P(f)\sum_{j=1}^{m_{1}}\mathbb{G}(A_{j})=P(f)\mathbb{G}%
(1)=0$ and $P_{2}[f]^{t}\cdot\mathbb{G}[\mathcal{B}]=0$ almost surely, we see that
$\mathbb{G}_{\text{even}}^{(\infty)}(\mathcal{F})=\mathbb{G}_{\text{odd}%
}^{(\infty)}(\mathcal{F})=\mathbb{G}^{(\infty)}(\mathcal{F})$. Applying Proposition~\ref{Pro4}
with the supremum norms $||\cdot||_{m_1}$ and $||\cdot||_{m_2}$ further yields, for any $m\geqslant0$ and $c_{0}=m_1 c_{1}+m_2 c_{2}$,%
\begin{align}
\left\Vert\mathbb{G}^{(2m)}-\mathbb{G}_{\text{even}}^{(\infty)}\right\Vert_{\mathcal{F}} &
=\left\Vert (S_{1,\text{even}}^{(m-1)}-S_{1,\text{even}})^{t}\cdot\mathbb{G}%
[\mathcal{A}]+(S_{2,\text{even}}^{(m-2)}-S_{2,\text{even}})^{t}\cdot\mathbb{G}%
[\mathcal{B}]\right\Vert _{\mathcal{F}}\nonumber \\
&  \leqslant c_{0}\max(\lambda_{1},\lambda_{2})^{m-2}Z, \label{Zpair}
\end{align}
where $Z=\max(Z_1,Z_2)$, and%
\begin{equation}
\left\Vert\mathbb{G}^{(2m+1)}-\mathbb{G}_{\text{odd}}^{(\infty)}\right\Vert_{\mathcal{F}%
}\leqslant c_{0}\max(\lambda_{1},\lambda_{2})^{m-1}Z.\label{Zimpair}
\end{equation}
Let $\varepsilon_{N}=q_{N}%
c_{0}\max(\lambda_{1},\lambda_{2})^{N/2}$ and $q_{N}=F_{Z}^{-1}(c_{0}%
\max(\lambda_{1},\lambda_{2})^{N/2})$, which is well defined for $N$ large enough. From (\ref{Zpair}) and (\ref{Zimpair}) we deduce that
\[
\mathbb{P}\left(  \left\Vert \mathbb{G}^{(N)}-\mathbb{G}^{(\infty)}\right\Vert
_{\mathcal{F}}>\varepsilon_{N}\right)  \leqslant\mathbb{P}\left(
Z>q_{N}\right)  \leqslant c_{0}\max(\lambda_{1},\lambda_{2})^{N/2},%
\]
whence an upper bound for the L\'{e}vy-Prokhorov distance%
\[
d_{LP}(\mathbb{G}^{(N)},\mathbb{G}^{(\infty)})\leqslant\max\left(
\mathbb{P}\left(  \left\Vert \mathbb{G}^{(N)}-\mathbb{G}^{(\infty)}\right\Vert
_{\mathcal{F}}>\varepsilon_{N}\right)  ,\varepsilon_{N}\right)  \leqslant
c_{0}q_{N}\max(\lambda_{1},\lambda_{2})^{N/2}.
\]
Let $\Phi$ denote the standard Gaussian distribution function, $c_{5}=m_{1}+m_{2}$ and $c_{4}^{2}%
=\max_{D\in\mathcal{A}\cup\mathcal{B}}\left\{  P(D)(1-P(D))\right\}$. The union bound%
\[
\mathbb{P}\left(  Z>\lambda\right)  \leqslant c_{5}\left(  1-\Phi\left(
\frac{\lambda}{c_{4}}\right)  \right)  \leqslant\frac{c_{5}c_{4}}{\sqrt{2\pi
}\lambda}\exp\left(  -\frac{\lambda^{2}}{2c_{4}^{2}}\right),
\]
shows that $q_{N}= c_{6}c_{4}\sqrt{N\log(1/c_{0}\max(\lambda_{1},\lambda_{2}))}$ for some $c_{6}>0$.

\appendix\section*{Appendix}\label{app}

\subsection{Elementary example}\label{RRExemple}

The Raking-Ratio algorithm changes the weights of cells of a contingency
table in such a way that given margins are respected, just as if the sample should
have respected the expected values of known probabilities. Let us illustrate the method from the following basic two-way contingency table.
\begin{center}
  \begin{tabular}{|c|c|c|c|c|c|}
      \hline
      $\mathbb{P}_n(1_{A_i^{(1)} \cap A_j^{(2)}})$ & $A_1^{(2)}$ &$ A_2^{(2)}$ & $A_3^{(2)}$ 
      &Total&{\color{black} Excepted total}\\ \hline
      $A_1^{(1)}$&0.2&0.25&0.1&0.55&{\color{black}0.52}\\ \hline
      $A_2^{(1)}$&0.1&0.2&0.15&0.45&{\color{black}0.48}\\ \hline
      Total&0.3&0.45&0.25&1&\\ \hline
      Excepted total&{\color{black}0.31}&{\color{black}0.4}&{\color{black}0.29}&&$N=0$\\ \hline
  \end{tabular}
  \end{center}
  The margins of this sample differ from the known margins, here called expected total. Firstly the weights of lines are corrected, hence each cell is multiplied by the ratio of the expected total and the actual one, this is step $N=1$.
  \begin{center}
  \begin{tabular}{|c|c|c|c|c|c|}
      \hline
      $\mathbb{P}_n^{(1)}(1_{A_i^{(1)} \cap A_j^{(2)}})$ & $A_1^{(2)}$ &$ A_2^{(2)}$ & $A_3^{(2)}$ 
      &Total&{\color{black} Excepted total}\\ \hline
      $A_1^{(1)}$&0.189&0.236&0.095&0.52&{\color{black}0.52}\\ \hline
      $A_2^{(1)}$&0.11&0.21&0.16&0.48&{\color{black}0.48}\\ \hline
      Total&0.299&0.446&0.255&1&\\ \hline
      Excepted total&{\color{black}0.31}&{\color{black}0.4}&{\color{black}0.29}&&$N=1$\\ \hline
  \end{tabular}
  \end{center}
  The totals for each column are similarly corrected at step $N=2$. Typically the margins of the lines no longer match the expected frequencies. Here they move in the right direction. Some estimators based on $\mathbb{P}_n^{(2)}$ may be improved.
  \begin{center}
  \begin{tabular}{|c|c|c|c|c|c|}
      \hline
      $\mathbb{P}_n^{(2)}(1_{A_i^{(1)} \cap A_j^{(2)}})$ & $A_1^{(2)}$ &$ A_2^{(2)}$ & $A_3^{(2)}$ 
      &Total&{\color{black} Excepted total}\\ \hline
      $A_1^{(1)}$&0.196&0.212&0.108&0.516&{\color{black}0.52}\\ \hline
      $A_2^{(1)}$&0.114&0.188&0.182&0.484&{\color{black}0.48}\\ \hline
      Total&0.31&0.4&0.29&1&\\ \hline
      Excepted total&{\color{black}0.31}&{\color{black}0.4}&{\color{black}0.29}&&$N=2$\\ \hline
  \end{tabular}
  \end{center} The last two operations are repeated until stabilization. The 
  algorithm converges to the Kullback projection of the initial joint law. The rate depends only
  on the initial table compared to the desired marginals. It
  takes only 7 iterations in our case to match the expected margins. 
  \begin{center}
  \begin{tabular}{|c|c|c|c|c|c|}
      \hline
      $\mathbb{P}_n^{(7)}(1_{A_i^{(1)} \cap A_j^{(2)}})$ & $A_1^{(2)}$ &$ A_2^{(2)}$ 
      & $A_3^{(2)}$ &Total&{\color{black} Excepted total}\\ \hline
      $A_1^{(1)}$&0.199&0.212&0.109&0.52&{\color{black}0.52}\\ \hline
      $A_2^{(1)}$&0.111&0.188&0.181&0.48&{\color{black}0.48}\\ \hline
      Total&0.31&0.4&0.29&1&\\ \hline
      Excepted total&{\color{black}0.31}&{\color{black}0.4}&{\color{black}0.29}&&$N=7$\\ \hline
  \end{tabular}
  \end{center} The final raked frequencies are slightly moved away the initial ones,
  however this has to be compared with the natural sampling oscillation order 
  $1/\sqrt{n}$ -- insidiously $n$ was not mentioned. 
  For small samples such changes
  are likely to occur that may improve a large class of estimators, and worsen others. Our theoretical results showed that the improvement is uniform over a large class as $ n \to +\infty $ and $ N $ is fixed.

\subsection{Counterexample of Remark J}\label{AppCounterEx}

Let assume that $ P $ satisfies the following probability values\begin{center}
\begin{tabular}{|c|c|c|c|c|}
	\hline$ P(A_i \cap B_j) $ & $A_1$ &$A_2$ &$A_3$ & $ P(B_j) $ \\ \hline
    $ B_1 $ & 0.2 & 0.25 & 0.1 & 0.55 \\ \hline
    $ B_2 $ & 0.25 & 0.1 & 0.1 & 0.45 \\ \hline
    $ P(A_i) $ & 0.45 & 0.35 & 0.2 & \\ \hline
\end{tabular} 
\end{center} and that $f$ has the following conditional expectations \begin{center}
\begin{tabular}{|c|c|c|c|c|}
	\hline$ \mathbb{E}(f|A_i \cap B_j) $ & $A_1$ &$A_2$ &$A_3$ & $ \mathbb{E}^{(2)}(f) \simeq $ \\ \hline
    $ B_1 $ &0.75 & -0.5 & 0.5 & 0.136 \\ \hline
    $ B_2 $ & 0.5 & 0.25 & 0.5 & 0.444 \\ \hline
    $ \mathbb{E}^{(1)}(f) \simeq $ & 0.611 & -0.286 & 0.5 & \\ \hline
\end{tabular} 
\end{center} By supposing also that $ \mathbb{V}(f|A_i \cap B_j)=0.5 $ for all 
$i=1,2,3$ and $j=1,2 $ we can compute the theoretical limiting variances 
from Proposition \ref{Thm2}. We get 
$ \mathbb{V}(\mathbb{G}^{(0)}(f)) \simeq 0.734 $; 
$ \mathbb{V}(\mathbb{G}^{(1)}(f)) \simeq 0.563 $;
$ \mathbb{V}(\mathbb{G}^{(2)}(f)) \simeq 0.569 $;
$ \mathbb{V}(\mathbb{G}^{(3)}(f)) \simeq 0.402 $.
The fact that $ \mathbb{V}(\mathbb{G}^{(2)}(f)) > \mathbb{V}(\mathbb{G}^{(1)}(f)) $ 
shows that the variance doesn't decrease necessarily at each step. As 
predicted by Propositions \ref{Thm2b} and \ref{Thm3} we have 
$ \mathbb{V}(\mathbb{G}^{(N)}(f)) < \mathbb{V}(\mathbb{G}^{(0)}(f)) $ 
for $N=1,2,3$ and
$ \mathbb{V}(\mathbb{G}^{(3)}(f)) < \mathbb{V}(\mathbb{G}^{(1)}(f)) $.

\subsection{Raked empirical means over a class}\label{AppEx}

\noindent\textbf{General framework.} Many specific settings in statistics may be
modeled through $\mathcal{F}$. Typically $\mathcal{X}$ is of very large or
infinite dimension and each $f(X)$ is one variable with mean $P(f)$ in the
population. To control correlations between such variables one needs to extend
$\mathcal{F}$ into $\mathcal{F}_{\times}=\left\{  fg:f,g\in\mathcal{F}%
\right\}  $ and consider the covariance process $\alpha_{n}^{(N)}%
(\mathcal{F}_{\times})$. Random vectors $(Y_{1},...,Y_{k})=(f_{1}%
(X),...,f_{k}(X))$ can in turn be combined into real valued random variables
$g_{\theta}(X)=\varphi_{\theta}(Y_{1},...,Y_{k})$ through parameters $\theta$
and functions $g_{\theta}$ that should be included in $\mathcal{F}$ and so on.
Consider for instance $g_{\theta}(X)=\theta_{1}Y_{1}+...+\theta_{k}%
Y_{k}+\varepsilon_{\sigma}(X)$ with a collection of possible residual
functions $\varepsilon_{\sigma}$ turning part of the randomness of $X$ into a noise with variance $\sigma^{2}$. The (VC) or (BR) entropy of $\mathcal{F}$ rules the variety and
complexity of models or statistics one can simultaneously deal with. 
We refer to Pollard~\cite{Poll90}, Shorack and Wellner~\cite{ShoWell} and Wellner~\cite{Well92} for classical
statistical models where an empirical process indexed by functions is easily identified.\smallskip

\noindent\textbf{Direct applications.} Since the limiting process 
$\mathbb{G}^{(N)}$ of $\alpha_{n}^{(N)}$ has less variance than 
$\mathbb{G}^{(0)}$, Theorem \ref{ProCVSpeed} can be applied to 
revisit the limiting behavior of
classical estimators or tests by using $\mathbb{P}_{n}^{(N)}$ instead of
$\mathbb{P}_{n}$ and prove that the induced asymptotic variances or risk
decrease. For instance, in the case of goodness of fit tests, the threshold
decreases at any given test level while the power increases against
any alternative distribution $Q$ that do not satisfy the margin conditions. As a matter of
fact, enforcing $\mathbb{P}_{n}^{(N)}$ to look like $P$ instead of the true $Q$ over
all $\mathcal{A}^{(N)}$ makes $\mathbb{P}_{n}^{(N)}$ go very far from $P$\ on
sets where $Q$ was already far from $P$.\smallskip

\noindent\textbf{Example: two raked distribution functions.} Let $(X,Y)$ be a real centered Gaussian vector with covariance matrix $ \left( \begin{smallmatrix}
	3&-1 \\ -1 & 1
\end{smallmatrix} \right)$. We consider the raked joint estimation of the two distribution functions $F_X, F_Y$. An auxiliary information provides their values at points $-2$ to $2$, every $0.5$. The class $ \mathcal{F} $ we need contains for all $ t \in \mathbb{R} $ the functions $ f_t^X(x,y) = \mathbf{1}_{]-\infty,t]}(x), f_t^Y(x,y) = \mathbf{1}_{]-\infty,t]}(y) $ thus (VC) holds.  For $ Z = X,Y $ let $ F_{Z,n}^{(N)}(t) = \sum_{Z_i \leqslant t} \mathbb{P}_n^{(N)}(\{Z_i\}) $ be the $N$-th raked empirical distribution function and write $Z_{(1)}\leqslant \dots \leqslant Z_{(n)}$ the order statistics. To exploit at best the information we use $N=2m$, $F_{X,n}^{(2m-1)}$ and $F_{Y,n}^{(2m)}$. Consider $ d_{Z,n}^{(N)} = \sum_{i=1}^{n-1} (Z_{(i+1)}-Z_{(i)}) \vert F_{Z,n}^{(N)}(Z_{(i+1)})-F_Z(Z_{(i+1)})\vert $ which approximates on $[Z_{(1)},Z_{(n)}]$ the $ L^1 $-distance between $ F_{Z,n}^{(N)}$ and $ F_{Z} $. Denote $ \#_{Z,n}^{(N)} $ the random proportion of sample points where $ F_{Z,n}^{(N)} $ is closer to $ F_Z $ than $ F_{Z,n}^{(0)} $. The table below provides Monte-Carlo estimates of $ D_{Z,n}^{(N)} = \mathbb{E}(d_{Z,n}^{(N)}) $ and $ p_{Z,n}^{(N)} = \mathbb{E}(\#_{Z,n}^{(N)}) $ from 1000 simulations based on samples of size $ n=200$: \begin{center}
\begin{tabular}{|c|c|c|c|c|c|}
	\hline $Z$ & $D_{Z,n}^{(0)}$ & $ D_{Z,n}^{(10)}$ &$ D_{Z,n}^{(\infty)}$ & $ p_{Z,n}^{(10)} $ & $p_{Z,n}^{(\infty)} $  \\ \hline
    $X$ & 0.084  & 0.058& 0.065 & 0.752 & 0.724 \\ \hline
    $Y$ & 0.085 & 0.043 & 0.053 & 0.731 & 0.681 \\ \hline 
\end{tabular}
\end{center}
This shows some improvement, especially for $N=10$. For $n$ rather small it seems not always relevant to wait for the stabilization of the algorithm -- here denoted $N=
\infty$. Our theoretical results provide guaranties only for $N\leqslant N_0$, $N_0$ fixed and $n \geqslant n_0$ for $ n_0 > 0 $. We also observed on graphical representations that the way $F_{Z,n}^{(N)}$ leaves $F_{Z,n}^{(0)}$ to cross $F_Z$ at the known points tends to accentuate the error at a few short intervals where $F_{Z,n}^{(0)}$ is far from $F_Z$. This is less the case as the auxiliary information partition is refined or the sample size increases.
\smallskip

\noindent\textbf{Example: raked covariance matrices.} Given $d\in\mathbb{N}_{\ast}$ and
$f_{1},...,f_{d}$ let $\mathbb{V}(Y)$ denote the
covariance matrix of the random vector $Y=(f_{1}(X),...,f_{d}(X))$ which we
assume to be centered for simplicity. Instead of the empirical
covariance $\mathbb{V}_{n}^{(0)}(Y)=n^{-1}%
{\textstyle\sum\nolimits_{i=1}^{n}}
Y_{i}^{t}Y_{i}$ consider its raked version%
\[
\mathbb{V}_{n}^{(N)}(Y)=\left(  \left(  \mathbb{P}_{n}^{(N)}(f_{i}%
f_{j})\right)  _{i,j}\right)  .
\]
Let $\left\Vert \mathbb{\cdot
}\right\Vert $ denote the Froebenius norm and define
\[
\varphi_{Y}(\alpha_{n}^{(N)})=\sqrt{n}\left\Vert \mathbb{V}
_{n}^{(N)}(Y)-\mathbb{V}(Y)\right\Vert.
\]
In other words,%
\[
\varphi_{Y}^{2}(\alpha_{n}^{(N)})=%
{\displaystyle\sum\limits_{i=1}^{d}}
{\displaystyle\sum\limits_{j=1}^{d}}
\left(  \alpha_{n}^{(N)}(f_{i}f_{j})\right)  ^{2}, \quad
\varphi_{Y}^{2}(\mathbb{G}^{(N)})=
{\displaystyle\sum\limits_{i=1}^{d}}
{\displaystyle\sum\limits_{j=1}^{d}}
\left(  \mathbb{G}^{(N)}(f_{i}f_{j})\right)  ^{2}.
\]
In the context of Proposition \ref{ProBerryEsseen} observe that $\varphi_{Y}$ is $\left(  \left\Vert \mathbb{\cdot}\right\Vert
_{\mathcal{F}},\left\Vert \mathbb{\cdot}\right\Vert \right)  $-Lipshitz with
parameter $C_{1}=d$. Clearly $\varphi_{Y}(\mathbb{G}^{(N)})$ has a bounded
density since $\varphi_{Y}^{2}(\mathbb{G}^{(N)})$ is a quadratic form with Gaussian
components and has a modified $\mathcal{X}^{2}$ distribution. Choosing a finite collection of such $\varphi_{Y}$ ensures that $C_{2}<+\infty
$. More generally by letting $(f_{1},...,f_{d})$ vary among a small entropy infinite subset
$\mathcal{L}_{d}$ of $\mathcal{F}^{d}$ and imposing some regularity or localization constraints to the $f_{i}$ one may have $C_{2}<+\infty$ while $\{f_i f_j:f_i,f_j\in\mathcal{F}\}$ satisfies (BR). The largest $C_{2}$ still works for
$\mathcal{L}=%
{\textstyle\bigcup\nolimits_{d\leqslant d_{0}}}
\mathcal{L}_{d}$. Therefore Proposition \ref{ProBerryEsseen} guaranties
that%
\[
\max_{\begin{subarray}{c}
	0\leqslant N\leqslant N_{0} \\
    d\leqslant d_{0}
\end{subarray}}
\sup_{\begin{subarray}{c}
    (f_{1},...,f_{d})\in\mathcal{L}_{d} \\
    x>0
\end{subarray}}
\left\vert \mathbb{P}\left(
\varphi_{Y}(\alpha_{n}^{(N)})\leqslant x\right)
-\mathbb{P}\left(\varphi_{Y}(\mathbb{G}^{(N)})\leqslant x\right)
\right\vert \leqslant d_{1}C_{1}C_{2}v_{n},
\]
where it holds, for all $N\leqslant N_{0}$, $d_0\leqslant d_{1}$, $(f_{1},...,f_{d})\in\mathcal{L}_{d}$
and $x>0$,
\[
\mathbb{P}\left(\varphi_{Y}(\mathbb{G}^{(N)})\leqslant x\right) \leqslant
\mathbb{P}\left(\varphi_{Y}(\mathbb{G}^{(0)})\leqslant x\right),
\]
by the variance reduction property of Proposition \ref{Thm2b}. Hence we asymptotically have $\mathbb{P}(
\varphi_{Y}(\alpha_{n}^{(N)})\leqslant x)<
\mathbb{P}(
\varphi_{Y}(\alpha_{n}^{(0)})\leqslant x)-\varepsilon$ uniformly among $Y$ such that 
$\mathbb{P}\left(\varphi_{Y}(\mathbb{G}^{(N)})\leqslant x\right)<
\mathbb{P}\left(\varphi_{Y}(\mathbb{G}^{(0)})\leqslant x\right)
-2\varepsilon$, for any fixed $\varepsilon>0$.

\end{document}